\documentclass[11pt]{article}

\usepackage{lineno,hyperref}
\usepackage{graphicx,dcolumn,amsmath,latexsym,amssymb}
\usepackage{amsfonts,subfigure,caption}
\usepackage{amsthm}
\usepackage{color,bm}
\usepackage{placeins}
\usepackage{trimclip}
\usepackage{placeins}

\usepackage[ruled,vlined]{algorithm2e}

\SetCommentSty{mycommfont}
\SetKwInput{KwInput}{Input}                
\SetKwInput{KwOutput}{Output} 

\usepackage{tikz}
\usepackage{pgfplots, capt-of}
\pgfplotsset{compat = newest}

\theoremstyle{definition} \newtheorem{definition}{Definition}
\theoremstyle{plain} 
\theoremstyle{plain} 
\theoremstyle{plain} \newtheorem{theorem}[definition]{Theorem}
\theoremstyle{plain}

\newcommand{\cbb}{\mathbb{C}}

\newcommand{\Ecal}{\mathcal{E}}

\def\K{{\mathcal{K}}}
\def\S{{\mathcal{S}}}

\def\E{{\mathcal{E}}}

\def\Em{{\mathcal{E}^{\rho}_{\text{max}}}}
\def\Tm{{T^{\rho}_{\text{max}}}}
\newcommand{\argmax}{\operatorname{argmax}}

\def\bbeta{\boldsymbol{\beta}}
\def\bxi{\boldsymbol{\xi}}

\def\RR{{\mathbb{R}}}
\def\EE{{\mathbb{E}}}
\def\CC{{\mathbb{C}}}
\def\ZZ{{\mathbb{Z}}}

\def\hu{\widehat{u}}
\def\w{\mathbf{w}}

\def\tg0{\tilde{g}_{\Ecal_0,T}}

\makeatletter
\providecommand*{\upY}{%
  \mathbin{%
    \mathpalette\@updownY{0}%
  }%
}
\providecommand*{\downY}{%
  \mathbin{%
    \mathpalette\@updownY{1}%
  }%
}
\providecommand*{\UpDownYFactor}{1}
\newcommand*{\@updownY}[2]{%
  \sbox0{$#1+\m@th$}%
  \dimen2=.5\dimexpr\wd0-\ht0-\dp0\relax
  \sbox2{$#1\vcenter{}$}%
  \dimen4=\dimexpr\ht0-\ht2\relax
  \setbox0=\hbox to 0pt{%
    \hss
    \clipbox{%
      0pt %
      {\dimexpr\totalheight-\UpDownYFactor\dimen4\relax} %
      0pt %
      -\dimen2%
    }{$#1|$}%
    \hss
  }%
  \ht0=\dimexpr\ht0-\dimen2\relax
  \kern\dimen2 %
  \raise\ht2\hbox{%
    \ifnum#2=0 %
      {\rotatebox{120}{\copy0}}%
      \copy0 %
      {\rotatebox{-120}{\copy0}}%
    \else
      {\rotatebox{60}{\copy0}}%
      {\rotatebox{180}{\copy0}}%
      {\rotatebox{-60}{\copy0}}%
    \fi
  }%
  \kern\dimen2 %
}
\makeatother

\begin{document}

\title{Singularity Formation in the Deterministic and Stochastic Fractional Burgers Equation}
\author{Elkin Ram\'irez and Bartosz Protas
\\ 
Department of Mathematics \& Statistics, \\
McMaster University  \\
Hamilton, Ontario L8S4K1, CANADA }

\date{\today}
\maketitle

\begin{abstract}
  This study is motivated by the question of how singularity formation
  and other forms of extreme behavior in nonlinear dissipative partial
  differential equations are affected by stochastic excitations. To
  address this question we consider the 1D fractional Burgers equation
  with additive colored noise as a model problem. This system is
  interesting, because in the deterministic setting it exhibits
  finite-time blow-up or a globally well-posed behavior depending on
  the value of the fractional dissipation exponent. The problem is
  studied by performing a series of accurate numerical computations
  combining spectrally-accurate spatial discretization with a
  Monte-Carlo approach. First, we carefully document the singularity
  formation in the deterministic system in the supercritical regime
  where the blow-up time is shown to be a decreasing function of the
  fractional dissipation exponent. Our main result for the stochastic
  problem is that there is no evidence for the noise to regularize the
  evolution by suppressing blow-up in the supercritical regime, or for
  the noise to trigger blow-up in the subcritical regime. However, as
  the noise amplitude becomes large, the blow-up times in the
  supercritical regime are shown to exhibit an increasingly
  non-Gaussian behavior. Analogous observations are also made for the
  maximum attained values of the enstrophy and the times when the
  maxima occur in the subcritical regime.
\end{abstract}

\begin{flushleft}
\end{flushleft}


\section{Introduction}
\label{sec:intro}

A standard mathematical model describing the motion of viscous
incompressible fluids is the Navier-Stokes system governing the
evolution of a velocity vector field and a scalar pressure field
resulting from the conservation of mass and momentum in such fluids.
Despite the ubiquity of this model in diverse applications spanning
many different areas of science and engineering, our understanding of
some of its key mathematical properties is still far from
satisfactory. Most importantly, it is not known whether the
Navier-Stokes system in three dimensions (3D) is globally well-posed
in the classical sense. More precisely, for arbitrary smooth initial
data smooth (classical) solutions of the 3D Navier-Stokes system have
been shown to exist for finite times only and formation of
singularities in finite time has not been ruled out
\cite{d09,Robinson2020}. On the other hand, suitable weak solutions,
which are possibly non-smooth, have been shown to exist globally in
time. The significance of this regularity problem for the 3D
Navier-Stokes system has been recognized by the Clay Mathematics
Institute which included it among its seven ``Millennium Problems''
posed as challenges for the mathematical community \cite{f00}.

Given the difficulties in investigating the Navier-Stokes system in
3D, a lot of research has been focused on the study of various
simplified models, often in one dimension (1D). Among many such
models, the fractional Burgers equation stands out, because depending
on the value of the fractional dissipation exponent $\alpha$, this
system is either globally well-posed (in the critical and subcritical
regime when $\alpha \ge 1/2$), or exhibits finite-time blow-up (in the
supercritical regime when $\alpha < 1/2$) \cite{kns08}. Thus, there is
a certain analogy with the 3D Navier-Stokes system which is also known
to be globally well-posed in the classical sense in the presence of
fractional dissipation with exponent $\alpha \ge 5/4$ \cite{kp12}.
Moreover, when $\alpha = 1$ the fractional system becomes the
well-known viscous Burgers equation which is arguably the most
commonly studied 1D model of fluid flow \cite{bk07}.

An interesting problem is how the behavior of various hydrodynamic
models, especially in regard to possible singularity formation, may be
affected by the presence of noise represented by a suitably-defined
stochastic forcing.  More specifically, the key question is whether
via some interaction with the nonlinearity {and dissipation present in
  the system} such stochastic forcing may accelerate or delay the
formation of a singularity, or perhaps even prevent it entirely
\cite{f15}.  These questions are of course nuanced by the fact that
they may be considered either for individual trajectories or in
suitable statistical terms.  The idea that stochastic excitation could
act to re-establish global well-posedness in a system exhibiting a
finite-time blow-up in the deterministic setting has been considered
for some time, although more progress has been made on the related
problem of restoring uniqueness.  There are in fact some model
problems, including certain transport equations \cite{f15} and some
versions of the Schr\"odinger equation \cite{ddm02a}, whose behavior
is indeed regularized by noise. While there are a few related results
available for the 3D Navier-Stokes and Euler equations \cite{f15b},
here we mention the studies \cite{ar09,ar10} where it was shown that
singularity formation (gradient blow-up) in the inviscid Burgers
equation can be prevented by a certain stochastic excitation of the
associated Lagrangian particle trajectories.

There exists a large body of literature devoted to investigations of
stochastically forced Burgers equation used as a model for
three-dimensional (3D) turbulence. Below we mention a few landmark
studies and refer the reader to the survey paper \cite{bk07} for
additional details and references. The majority of these
investigations aimed to characterize the solutions obtained in
statistical equilibrium, attained by averaging over sufficiently long
times, in terms of properties of the stochastic forcing. Given the
motivation to obtain insights about actual turbulent flows, the main
quantities of interest in these studies were the scaling of the energy
spectrum, evidence for intermittency in the anomalous scaling of the
structure functions and the statistics of $\partial_x u$, such as the
tails (exponential vs.~algebraic) of its probability density function
\cite{cy95a,cy95b,ztg97}. Remarkably, some of these results were also
established with mathematical rigor \cite{b14}.  The aforementioned
quantities were also studied in flows evolving from stochastic initial
data \cite{gk93}. In this context we mention the investigations
\cite{s92,saf92} which focused on the statistics of shock waves in the
limit of vanishing viscosity $\nu$.  As regards technical
developments, a number of interesting results were obtained using
optimization-based instanton formulations {\cite{mkv16,bfkl97,ggs15}}.

Questions concerning extreme behavior in stochastic Burgers flows, as
quantified by the growth of certain Sobolev norms of the solutions,
and how it relates to the deterministic case \cite{ap11a} were
recently investigated in \cite{PocasProtas2018}.  Problems related to
singularity formation in the dispersive Burgers equation were studied
numerically in \cite{ks15a}.  However, to the best of our knowledge,
questions concerning the effect of stochastic forcing on the solutions
of the fractional Burgers equation, especially in the supercritical
regime where finite-time blow-up is known to occur, have not been
considered. In this study we use carefully designed, highly accurate
numerical computations to provide new insights about two related
problems concerning the deterministic and stochastic fractional
Burgers equation.

The first goal is to provide a precise quantitative characterization
of the singularity formation, in terms of the blow-up times, in the
deterministic fractional Burgers equation in the supercritical regime.
These results will then serve as a point of departure for our
investigation of the second problem where we will be interested in how
this singular behavior is affected by stochastic forcing. The key
finding here is that while in the presence of a stochastic forcing
gradient singularity still occurs in supercritical fractional Burgers
flows, mean blow-up times become shorter as the magnitude of the
stochastic forcing increases. Moreover, for large amplitudes of the
stochastic forcing the probability distribution of the blow-up times
becomes increasingly non-Gaussian and develops algebraic tails
corresponding to delayed blow-up times. This means that while on
average the effect of the stochastic forcing is to accelerate the
singularity formation, situations where blow-up is delayed as compared
to the deterministic case also become more likely.

The structure of the paper is as follows: in the next section we
introduce the deterministic and stochastic variants of the fractional
Burgers equation, state some basic facts about these systems and
define some diagnostic quantities; then, in Section \ref{sec:numer} we
discuss the numerical approaches used to solve these problems and
evaluate the diagnostic quantities; our computational results for the
deterministic and stochastic problem are presented in Sections
\ref{sec:results_deter} and \ref{sec:results_stoch}, respectively,
whereas discussion and conclusions are deferred to Section
\ref{sec:final}.

\section{Fractional Burgers Equation}
\label{sec:Burgers}

In this section we introduce the fractional Burgers equation in its
deterministic and stochastic versions. For simplicity, these problems
are considered on a periodic domain $(0,2\pi)$.  We will also discuss
some of their key properties and will introduce diagnostic quantities
useful for characterizing the regularity of their solutions.
Hereafter, we will use the Sobolev space $H^{s}_{p}(0, 2\pi)$ defined
for $s \in \RR^+$ as \cite{af05}
\begin{equation}
H^{s}_{p}(0,2\pi) := \left\{v\in L^{2}_{p}(0, 2\pi): ||v||_{H^{s}_{p}(0, 2\pi)}^{2}=\sum_{k=-\infty}^{\infty}\left(1+|k|^2\right)^{s}|\widehat{v}_{k}|^2<\infty\right\},
\label{eq:Hs}
\end{equation}
where $L^{2}_{p}(0, 2\pi)$ stands for the space of square-integrable
$2\pi$-periodic functions and $[\,\widehat{\cdot}\,]_k$ represents the
Fourier coefficient corresponding to the wavenumber $k \in \ZZ$
(``:='' means ``equal to by definition'').

\subsection{Deterministic Problem}
\label{sec:BurgersDet}

We consider the 1D fractional Burgers equation
\begin{subequations}
\label{FBE}
\begin{align}
\partial_{t}u+\frac{1}{2}\,\partial_{x}u^{2}+\nu\,(-\Delta)^{\alpha} u=0 \quad & \mbox{in}\ (0,T]\times (0,2\pi),\label{FBEa}\\
\text{periodic boundary conditions} \quad& \mbox{for}\ t\in(0,T],\label{FBEb}\\
u(0,x)=g(x)\quad & \mbox{for}\ x\in(0,2\pi),\label{FBEc}
\end{align}
\end{subequations}
where $\nu>0$ is the viscosity coefficient and $(-\Delta)^{\alpha}$ is
the fractional Laplacian with $\alpha\in[0,1]$ the fractional
dissipation exponent.  For a sufficiently smooth function $v \; : \;
[0,2\pi] \rightarrow \RR$, the fractional Laplacian is defined in
terms of the relation
\begin{equation}
\left[ \widehat{(-\Delta)^\alpha v} \right]_k := |k|^{2\alpha} [ \widehat{v} ]_k, \qquad k \in \ZZ.
\label{eq:fraclap}
\end{equation}
In system \eqref{FBE} $T>0$ represents the length of the time window
and $g\in H_{p}^{1}(0, 2\pi)$ is the initial condition assumed to have
zero mean. This last assumption reflects the fact that the mean of the
solution is preserved by the evolution governed by system \eqref{FBE}.
Some key results addressing questions about existence of solutions of
system \eqref{FBE} are found in \cite{kns08} and are briefly
summarized below.

\begin{theorem}[subcritical case]
  \label{thm:subcrit}
  Assume that $\alpha>1/2$ and the initial data $g \in H^{s}$, $s>
  3/2-2\alpha,$ $s\geq0.$ Then, there exists a unique global solution
  of problem \eqref{FBE} that belongs to
  $C(\left[0,\infty\right),H^{s})$ and is real analytic in $x$ for
  $t>0$.
\end{theorem} 

\begin{theorem}[critical case]
\label{thm:crit}
Assume that $\alpha = 1/2$ and $g\in H^{s},$ $s > 1/2$. Then, there
exists a global solution of the system \eqref{FBE} which is real
analytic in $x$ for any $t > 0$.
\end{theorem} 

\begin{theorem}[supercritical case]
\label{thm:supercrit}
Assume that $0 <\alpha < 1/2$. Then, there exists smooth periodic
initial data $g$ such that the solution $u$ of \eqref{FBE} blows up in
$H^{s}$ for each $s > 3/2 -2\alpha$ in a finite time.
\end{theorem}

\begin{figure}
\begin{center}
\begin{tikzpicture}
\begin{axis}[
xmin = 0, xmax = 1/2,
ymin = 0.5, ymax = 3/2,
width = 8cm,
height = 6.5cm,
xtick={0,0.25,0.5},
xticklabels={$0$, $\frac{1}{4}$,  $\frac{1}{2}$},
ytick={0.5,1,1.5},
yticklabels={$\frac{1}{2}$, $1$,  $\frac{3}{2}$},
smooth,
xlabel=$\alpha$,
ylabel=$s$,
]
\addplot[color=red,domain = 0:1/2, ultra thick]{(3/2)-2*x};
\addplot [fill=gray!90] (1/2,1/2) -- (1/2,3/2) -- (0,3/2) -- cycle;
\addplot[domain = 0:1/4, ultra thick,dash dot]{1};
\addplot[domain = 1/4:1/2, ultra thick]{1};
\addplot[color=red,domain = 0:1/2, ultra thick]{(3/2)-2*x};
\addlegendentry{$s=\frac{3}{2}-2\alpha$};

\end{axis}
\end{tikzpicture}
\captionof{figure}{The shaded region represents the set of values
  $(\alpha,s)$, cf.~\eqref{eq:fraclap} and \eqref{eq:Hs}, for which
  finite-time blow-up in system \eqref{FBE} is predicted by Theorem
  \ref{thm:supercrit}. The horizontal line represents the space
  $H^1_p$.}
\label{fig:blowup}
\end{center}
\end{figure}
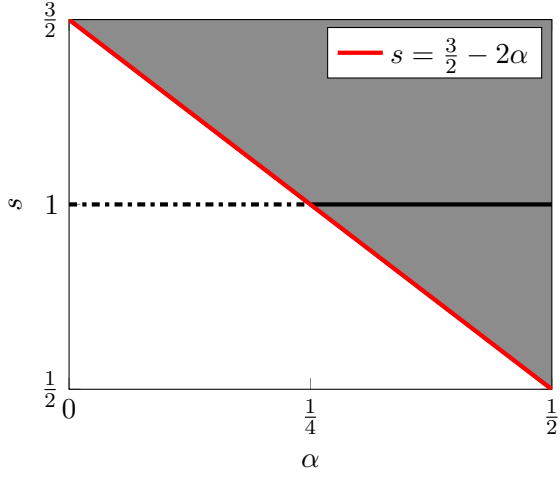

Our main interest is in the existence of solutions in the Sobolev
space $H^1_p(0,2\pi)$ and from Theorems \ref{thm:subcrit} and
\ref{thm:crit} it is clear that system \eqref{FBE} is globally
well-posed in this space when $\alpha \ge 1/2$. However, as regards
blow-up in the supercritical regime, the situation is more nuanced,
since Theorem \ref{thm:supercrit} predicts that blow-up in
$H^1_p(0,2\pi)$ occurs only when $\alpha \in (1/4,1/2)$, cf.~Figure
\ref{fig:blowup}. Thus, Theorem \ref{thm:supercrit} is inconclusive as
regards what happens when $\alpha \in [0,1/4]$ and in Section
\ref{sec:results_deter} we present computational evidence indicating
that finite-time blow-up occurs in this regime as well. Interestingly,
as was conjectured in \cite{Yun2018}, if $\alpha \in [0,1/4]$, system
\eqref{FBE} may not be even locally well-posed in $H^{1}_{p}(0,2\pi)$
for {\it some} initial data. Other results concerning the regularity
of the fractional Burgers equation were obtained \cite{adv07,ccs10},
whereas in \cite{ddl09} it was shown that under certain conditions on
the initial data blow-up occurs in $W^{1,\infty}$ for all $\alpha <
1/2$.

\subsubsection{Limiting cases}
\label{sec:limiting}

Here we briefly discuss the forms of system \eqref{FBE} corresponding
to the limits $\nu \rightarrow 0$ and $\alpha \rightarrow 0$, the
latter with $\nu > 0$ fixed. As regards the first limit, we recall
that \eqref{FBE} turns into the inviscid Burgers equation
\begin{subequations}
\label{eq:iBE}
\begin{align}
\partial_{t}u+\frac{1}{2}\,\partial_{x}u^{2}=0 \quad & \mbox{in}\ (0,T^*)\times (0,2\pi),\label{inviscida}\\
\text{periodic boundary conditions} \quad& \mbox{for}\ t\in(0,T^*), \\
u(0,x)=g(x)\quad & \mbox{for}\ x\in(0,2\pi)
\end{align}
\end{subequations}
which is well understood \cite{kl04}. In particular, it is known that
in solutions of \eqref{eq:iBE} corresponding to smooth nonzero initial
conditions a shock-type singularity forms at the time
\begin{equation}
T^{*}_0= -\frac{1}{\inf_x g'(x)}.
\label{eq:TsiBE}
\end{equation}
On the other hand, in the second case corresponding to the limit
$\alpha \rightarrow 0$ with $\nu > 0$, system \eqref{FBE} turns into
\begin{subequations}
\label{eq:fBE}
\begin{align}
\partial_{t}u+\frac{1}{2}\,\partial_{x}u^{2}+\nu\,u=0 \quad & \mbox{in}\,(0,T^*)\times (0,2\pi),\\
\text{periodic boundary conditions} \quad& \mbox{for}\ t\in(0,T^*), \\
u(0,x)=g(x)\quad & \mbox{for}\,x\in(0,2\pi)
\end{align}
\end{subequations}
which can be solved using elementary techniques \cite{Ramirez2020}. It
can be shown that if $\inf_x g'(x)+\nu<0$ a singularity forms in its
solutions at the time
\begin{equation}
T^{*}_1= -\frac{1}{\nu}\ln\left[\frac{\nu}{\inf_x g'(x)}+1\right].
\label{eq:TsfBE}
\end{equation}

\subsection{Stochastic version}
\label{sec:BurgersStoch}

As regards the stochastic model, in the context of dissipative systems
such as \eqref{FBE} one typically studies additive noise. The reason
is that, as argued in \cite[Section 5.5.2]{f15}, multiplicative noise
tends to have similar effect to dissipative terms, so if the equation
already involves such a term, then no major qualitative changes in the
solution behavior can be expected. We will thus consider the following
stochastic form of the fractional Burgers system \eqref{FBE}
\begin{subequations}
\label{SBE}
\begin{align}
\partial_{t}u+\frac{1}{2}\,\partial_{x}u^{2}+\nu\,(-\Delta)^{\alpha} u=\zeta(t,x) \quad & \mbox{in}\ (0,T]\times (0,2\pi), \label{SBEa}\\
\text{periodic boundary conditions} \quad& \mbox{for}\ t\in(0,T], \\
u(0,x)=g(x)\quad & \mbox{for}\ x\in(0,2\pi),\label{SBEd}
\end{align}
\end{subequations}
where $\zeta$ is a random field. Therefore, at any point $(t,x) \in
(0,T]\times [0,2\pi]$ our solution becomes a random variable
$u={u(t,x;\omega)}$ for $\omega$ in some probability space $\Omega$.
Inputs to stochastic systems are often modelled in terms of Gaussian
noise which is white both in space and in time, and is thus associated
with an infinite-variance Wiener process \cite{lps14}.  However, as
shown in \cite{PocasProtas2018}, such a noise model acts uniformly on
the entire wavenumber spectrum of the solution and therefore does not
ensure that its individual realizations are well defined in the
Sobolev space $H^1_p$. White noise is thus not suitable for the
problem considered here and we shall therefore adopt a formulation
where $\zeta$ is the derivative of a Wiener process with finite
variance, which is the most ``aggressive'' stochastic excitation still
leaving problem \eqref{SBE} with $\alpha \in (1/2,1]$ well-posed in
$H^1_p$.

We will therefore consider a  colored-in-space Gaussian noise
\begin{equation}
\label{eq:zeta}
\zeta(t,x) := \rho\frac{dW(t,x)}{dt}, \qquad x \in [0,2\pi], \ t \ge 0,
\end{equation}
where $\rho>0$ is a constant and $W(t,x)$ is a cylindrical Wiener
process given by the expression

\begin{equation}
\label{eq:W}
W(t,x) := \sum_{k=1}^N\gamma_{k}\,\beta_{k}(t)\,\chi_{k}(x)
\end{equation}
in which $\{\beta_{k}(t)\}_{k=1}^N$ are independent and identically
distributed (i.i.d.) standard Brownian motions,
$\{\chi_{k}\}_{k\in\mathbb{N}}$ are elements of a trigonometric
orthonormal basis of $L^{2}_{p}(0,2\pi)$ and
$\{\gamma_{k}\}_{k\in\mathbb{N}}$ are scaling coefficients. We will
set $\chi_{0}(x) = 1$, $\chi_{2k}(x) = \sqrt{2}\cos(kx)$ and
$\chi_{2k-1}(x) = \sqrt{2}\sin(kx)$ with $x \in [0,2\pi]$ and
$k=1,2,\dots$. For technical reasons, ansatz \eqref{eq:W} involves a
finite (but large) number $N$ of Fourier modes which in practice will
be chosen equal to the numerical resolution (such that noise will be
acting on all Fourier components in the spatial discretization of
system \eqref{SBE}).  As regards the scaling coefficients, we choose
them to be $l^2$-summable
\begin{equation}
\label{eq:gammas}
\gamma_{0}=0,\quad \gamma_{2k-1}=\gamma_{2k}=\frac{1}{k},\qquad k>0,
\end{equation}
such that $W(t,x)$ has finite variance in the limit $N \rightarrow
\infty$.  

In order to define solutions of the stochastic problem more precisely,
it is necessary to rewrite system \eqref{SBE} in the corresponding
differential form \cite{lps14}
\begin{subequations}
\label{SBE2}
\begin{align}
du=\left(-\frac{1}{2}\,\partial_{x}u^{2}-\nu\,(-\Delta)^{\alpha} u\right)dt+\rho\,dW\quad & \mbox{in}\ (0,T]\times (0,2\pi),\label{SBE2a}\\
\text{periodic boundary conditions} \quad& \mbox{for}\ t\in(0,T], \\
u(0,x)=g(x)\quad & \mbox{for}\ x\in(0,2\pi).\label{SBE2d}
\end{align}
\end{subequations}
Then, the relevant concept of a solution is the mild solution defined
as
\begin{equation}
\label{eq:SBE2mild}
u(t)=e^{-tA}g-\frac{1}{2}\int_{0}^{t} e^{-(t-s)A}\partial_{x}u^2 ds + \rho\int_{0}^{t} e^{-(t-s)A}dW(s),
\end{equation}
where $A:=\nu\,(-\Delta)^{\alpha}$ and the semigroup $e^{-tA}$ is
defined in terms of its action on the elements of the basis
$\{\phi_{k}\}_{k\in\mathbb{Z}}=\{e^{ikx}\}_{k\in\mathbb{Z}}$ of
$L^{2}_{p}(0,2\pi)$ as $e^{-tA}e^{ikx}=e^{-\nu
  t|k|^{2\alpha}}e^{ikx}$, whereas the second integral is understood
in It\^{o}'s sense \cite{lps14}. With the stochastic excitation
defined in \eqref{eq:W}, mild solutions \eqref{eq:SBE2mild} in the
subcritical regime ($\alpha \in (1/2,1]$) remain well-defined in the
space $H^1_p(0,2\pi)$ also in the limit $N \rightarrow \infty$
\cite{PocasProtas2018,Ramirez2020}.

\subsection{Diagnostic quantities}
\label{sec:diagn}

We now introduce two diagnostic quantities which will be used to
monitor the regularity of solutions to systems \eqref{FBE} and
\eqref{SBE} in computations. The first one is the enstrophy which is
proportional to the $H^1$ seminorm of the solution
\begin{equation}
\mathcal{E}(t) := \pi\int_{0}^{2\pi}\left| \partial_{x} u(t,x)\right|^2\,dx.
\label{eq:E}
\end{equation} 
This quantity is useful because its unbounded growth signals
singularity formation in the solutions of system \eqref{FBE}, a
property which also holds for 3D Navier-Stokes flows
\cite{d09,Robinson2020}.

The second quantity we will consider is the width of the analyticity
strip $\delta \; : \; [0,T] \rightarrow \RR^+ \cup \{ \infty \}$
characterizing the distance $\delta(t)$ from the real axis to the
nearest singularity in the complex extension of the solution
$u(t,\cdot)$ at time $t$. Clearly, the width of the analyticity strip
$\delta(t)$ vanishes as time $t$ approaches the blow-up time $T^*$
when the solution ceases to be analytic.  Since solutions of the
stochastic system \eqref{SBE} are not, in general, analytic in $x$,
the width of the analyticity strip cannot be used to characterized
their regularity.

Numerical approaches we use to approximate solutions of the
deterministic and stochastic problems \eqref{FBE} and \eqref{SBE}, and
the diagnostic quantities described above are discussed in the next
section.

\section{Numerical Approaches}
\label{sec:numer}

\subsection{Deterministic Fractional Burgers Equation}
\label{sec:numerdet}

The deterministic fractional Burgers system \eqref{FBE} is solved
numerically using a standard pseudo-spectral Fourier-Galerkin approach
\cite{b01,canuto:SpecMthd} in which the solution is approximated as
\begin{equation}
u(t,x) \approx u_{N}(t,x) := \sum_{k=-N/2+1}^{N/2} \widehat{u}_{k}(t)e^{ikx}, 
\label{eq:un}
\end{equation}
where $\widehat{u}_{k}(t)$, $k=-N/2,\dots,N/2$, are the Fourier
coefficients and $N = 2^{n}$ for some $n\in\mathbb{N}$ is the
resolution. We note that since the solution $u(t,x)$ is real-valued,
the Fourier coefficients satisfy the conjugate symmetry, i.e.,
$\widehat{u}_{-k}=\overline{\widehat{u}}_{k}$, $k=-N/2,\dots,N/2$,
where the bar denotes complex conjugation, such that only the Fourier
coefficients with $k>0$ need to be computed. Moreover, since we
consider initial data with zero mean and the mean of the solution is
preserved by system \eqref{FBE}, we have $\widehat{u}_0(t) = 0$,
$\forall t \ge 0$. Plugging ansatz \eqref{eq:un} into \eqref{FBE}
leads to the following system of ordinary differential equations
(ODEs)
\begin{subequations}
\label{FBEhat}
\begin{align}
\frac{d\widehat{{\bf u}}(t)}{dt}&={\bf r}(\widehat{{\bf u}}(t))+{\bf A}\widehat{{\bf u}}(t),\label{fouriersystema}\\
\widehat{{\bf u}}(0)&=\widehat{{\bf g}},
\label{fouriersystemb}
\end{align}
\end{subequations} 
where $\widehat{{\bf
    u}}(t)=\left[\widehat{u}_{1}(t),\widehat{u}_{2}(t),\dots
  ,\widehat{u}_{N/2}(t)\right]^{T} \in \CC^{N/2}$ and $\widehat{{\bf
    g}}=[\widehat{g}_{1},\widehat{g}_{2},\dots\widehat{g}_{N/2}]^{T} \in
\CC^{N/2}$, whereas ${\bf r}$ and ${\bf A}$ represent the nonlinear
and diagonal linear operators defined as
\begin{equation}
\left[{\bf r}(\widehat{{\bf u}}(t))\right]_{k}:=-\frac{1}{2}\,i\,k\,\widehat{\left[u^2(t)\right]}_{k},\quad 
\left[{\bf A}\widehat{{\bf u}}(t)\right]_{k}:=-\nu\,k^{2\alpha}\,\widehat{u}_{k}(t), \quad k=1,\dots, N/2.
\label{eq:rA}
\end{equation}
The nonlinear term in \eqref{eq:rA} is evaluated in the physical space
with discrete Fourier transforms computed using the FFT and dealiasing
performed based on the ``3/2 rule'' \cite{b01,canuto:SpecMthd}.
Integration in time is carried out using a hybrid approach combining
the Crank-Nicolson method with a three-step Runge-Kutta method applied
to the linear and nonlinear terms in \eqref{eq:rA}, respectively
\cite{NumRenaissance}.  

Given that we are interested in computing solutions of system
\eqref{FBE} as they approach the blow-up time $T^*$ in the
supercritical regime, cf.~Theorem \ref{thm:supercrit}, an important
element of the numerical approach is automatic grid refinement. The
spatial resolution is refined (by doubling $N$) each time the solution
$u(t,x)$ becomes marginally resolved which is signaled by the
magnitude of its Fourier coefficients corresponding to the largest
wavenumbers $k$ becoming larger than the machine precision, i.e., when
$|\widehat{u}_k| > \mathcal{O}(10^{-16})$ for $k \lessapprox N/2$.
Grid refinement in space then triggers a corresponding refinement of
the time step. Since resolution refinement significantly increases the
computational cost, time integration must be stopped at some point
before reaching the blow-up time $T^*$. Further details of the
approach presented here are described together with validation results
in \cite{Ramirez2020}.

\subsection{Stochastic Fractional Burgers Equation}
\label{sec:numerstoch}

A pseudo-spectral Fourier-Galerkin method with dealiasing is also used
to discretize problem \eqref{SBE} in space. Substituting
representation \eqref{eq:un} into \eqref{SBE2}, we obtain the
following system of stochastic ODEs
\begin{subequations}
\label{SBEhat}
\begin{align}
d\,\widehat{\bf u}& = \left({\bf r}(\widehat{{\bf u}}(t))+{\bf A}\widehat{{\bf u}}(t)\right)dt+\rho\,d{\bf W}(t),\\
\widehat{\bf u}(0)&=\widehat{\bf g}
\end{align}
\end{subequations}
in which $\widehat{{\bf u}}(t),$ ${\bf r}(\widehat{{\bf u}}(t))$ and
${\bf A}\widehat{{\bf u}}$ are as defined as above in Section
\ref{sec:numerdet} and ${\bf W}(t)=[W_{1}(t),\dots,W_{N/2}(t)]^{T}$,
where
\begin{equation}
W_{k}(t) := \frac{\sqrt{2}}{2k}\left(\beta_{2k}(t)-i\beta_{2k-1}(t)\right)
\label{eq:Wk}
\end{equation} 
and $\beta_{1},\dots,\beta_{N}$ are i.i.d standard Brownian motions.
Discretization in the probability space is carried out using a
Monte-Carlo approach where system \eqref{SBEhat} is solved $M$ times,
each time using a different random realization of the Brownian motions
in \eqref{eq:Wk}. Each such realization of system \eqref{SBEhat} is
discretized in time using a stochastic Runge-Kutta method of order
one-and-half introduced in \cite{Chang1987}. It is based on the
  following scheme
\begin{equation}
\begin{split}
{\bf Q}^{(n)}&:={\bf \widehat{u}}^{n}+\frac{1}{2}\Delta t\,{\bf f}({\bf \widehat{u}}^{n}),\\
{\bf Q}^{*(n)}&:={\bf \widehat{u}}^{n}+\frac{1}{2}\Delta t\,{\bf f}({\bf \widehat{u}}^{n}) +\frac{3}{2}\rho\sqrt{\Delta t}\,{\boldsymbol \beta},\\
{\bf \widehat{u}}^{n+1}&={\bf \widehat{u}}^{n} + \rho\Delta {\bf W}^{n} +\frac{1}{3}\Delta t\left[{\bf f}({\bf Q}^{n})+2{\bf f}({\bf Q}^{*(n)})\right], \qquad n = 1,2,\dots,
\label{RKstochastic}
\end{split}
\end{equation}
where $\Delta t$ is the time step, ${\bf \widehat{u}}^{n}$ is the
solution of system \eqref{SBEhat} at time $t^n := n\Delta t$, ${\bf
  f(\widehat{u})}:={\bf r}(\widehat{{\bf u}}(t))+ {\bf A}\widehat{{\bf
    u}}(t)$, whereas $\bbeta$ and $\Delta {\bf W}^{n} ={\bf
  W}^{n+1}-{\bf W}^{n}$ are defined as $\bbeta:=\frac{1}{2} \bxi +
\frac{\sqrt{3}}{6}{\boldsymbol{\eta}}$ and $\Delta {\bf
  W}^{n}:=\sqrt{\Delta t} \, \bxi$, where $\bxi$ and
$\boldsymbol{\eta}$ are two independent $\cbb^{N/2}$-valued Gaussian
random variables with a joint distribution $\mathcal{N}(0,I_{N/2})$.

Our choice of the Monte Carlo approach to noise sampling is motivated
by its well-understood convergence properties and straightforward
implementation. While more modern approaches, such as polynomial chaos
expansions, may in principle achieve faster convergence, they suffer
from much higher computational complexity.  Moreover, the nonlinear
term will have a rather complicated expression in the polynomial
orthonormal basis \cite{lps14}. Further details of the approach
presented here are described together with validation results in
\cite{Ramirez2020}. We add that given the structure of the stochastic
forcing in \eqref{eq:zeta}--\eqref{eq:gammas}, the resolution required
to accurately approximate the solutions is primarily determined by the
properties of the noise which do not change in time.  Thus, grid
refinement is not performed in the solutions of the stochastic problem.

\subsection{Evaluation of Diagnostic Quantities}
\label{sec:numerdiag}

Using Parceval's identity and representation \eqref{eq:un}, the
enstrophy \eqref{eq:E} is approximated as
\begin{equation}
\mathcal{E}(t) \approx 4\,\pi^2\sum_{k=1}^{N/2} k^2 |\widehat{u}_{k}(t)|^2.
\label{eq:En}
\end{equation} 
As regards the width of the analyticity strip, we use the method
introduced in \cite{ssf83}. It assumes that the Fourier spectrum of
the solution $u(t,x)$ can be expressed as
\begin{equation}
|\widehat{u}_{k}(t)|\sim C(t)\, |k|^{\tilde{\alpha}(t)}e^{-\delta(t)\,k},
\label{eq:d}
\end{equation}
where $\delta(t)$ is the width of the analyticity strip of $u(t,x)$,
$\tilde{\alpha}(t)$ is the order of the nearest complex singularity
and $C(t)$ is an adjustable parameter. An estimate of $\delta(t)$ can be
then obtained by minimizing the least-squares error between ansatz
\eqref{eq:d} and the amplitudes of the Fourier coefficients
$\hu_1(t),\dots,\hu_{N/2}(t)$.

\subsection{Estimates of the blow-up time}
\label{sec:numerTs}

One of the main aims of this study is to understand how the blow-up
time $T^{*}$ in the supercritical regime depends on the fractional
dissipation exponent $\alpha$ for a certain fixed initial condition.
We thus need to estimate $T^*$ based on the properties of the solution
$u(x,t)$ as $t \rightarrow \left(T^*\right)^-$.  From the discussion
in Section \ref{sec:diagn} we know that $\E(t) \rightarrow \infty$ and
$\delta(t) \rightarrow 0$ as $t \rightarrow \left(T^*\right)^-$.
Hence, following the ideas discussed in \cite{bk08,bb12}, the
functions of time $\E(t)$ and $\delta(t)$ can be locally approximated
for $t \in I$, where $I \subset [0,T^{*})$ is a closed interval, in
terms of the relation $c\left(T^{*}-t\right)^{\gamma}$, where $c>0$ is
a constant and the exponent $\gamma < 0$ for $\mathcal{E}(t)$ and
$\gamma > 0$ for $\delta(t)$. To fix attention, we will focus here on
estimating the blow-up time based on the time evolution of the
enstrophy $\E(t)$.  We will denote this quantity $T^{*}_{\mathcal{E}}$
and the estimate $T^{*}_{\delta}$ based on the time evolution of the
width of the analyticity strip $\delta(t)$ can be determined in a
analogous manner.  We assume that the numerical solution is available
at discrete times $0 \le t_i < T^*$ such that $t_{i}<t_{i+1}$.  A
local estimate of $T^{*}_{\mathcal{E}}$ can be then obtained by
solving the minimization problem
\begin{equation}
\min_{(c,\gamma,T^{*})\in\mathbb{R}^{3}} \sum_{t_{i}\in I} \left[\ln \frac{c(T^* - t_i)^\gamma}{\E(t_i)} \right]^2,
\label{eq:minf}
\end{equation}
in which, due to a large range of values attained by $\E(t)$, a
logarithmic objective function is preferred to the more standard
least-squares formulation. While problem \eqref{eq:minf} can in
principle be solved easily using standard tools of numerical
optimization \cite{nw00}, it is non-convex and the solution will be a
local minimizer depending on the choice of the initial guess.
Therefore, finding a suitable initial guess for the solution of
problem \eqref{eq:minf} is key to obtaining an accurate estimate
$T^{*}_{\mathcal{E}}$. To address this issue, we introduce a family of
$K$ ``sliding'' time windows $I_{j}\subset [0,T^{*})$, $j=1,\dots,K$,
centered at $t_{j P}$, where the positive integer $P$ represents the
number of discrete points the window is shifted as $j$ is incremented
(i.e., the window $I_{j}$ slides towards longer times as index $j$
increases until the last window $I_K$ coincides with the end of the
time interval on which a numerical solution of \eqref{FBE} has been
obtained). We then solve optimization problem \eqref{eq:minf} over $I
= I_j$ to obtain estimates for $c,\,T^{*}$ and $\gamma$ which we
denote $c_{\mathcal{E}}(t_j),\,T^{*}_{\mathcal{E}}(t_j)$ and
$\gamma_{\mathcal{E}}(t_j)$.  These estimates are then used as the
initial guess to solve problem \eqref{eq:minf} over $I = I_{j+1}$ and
so on until we finally solve this problem on $I = I_K$ which provides
our estimate of $T^{*}_{\mathcal{E}}$.  The idea of solving
optimization problem \eqref{eq:minf} repeatedly using such a family of
time windows is to obtain a good initial guess at early times $t$ when
the enstrophy $\E(t)$ is still varying slowly before solving these
problems close to the blow-up time.  A summary of this approach is
provided in Algorithm \ref{alg:Ts}. An estimate of $T^{*}_{\delta}$
can be obtained in the same way replacing $\E(t)$ with $\delta(t)$ in
\eqref{eq:minf}. Computational algorithms introduced in this section
have been implemented in {\tt julia} with post-processing performed in
{\tt MATLAB}.  Solutions of the deterministic and stochastic
fractional Burgers problems \eqref{FBE} and \eqref{SBE2} are discussed
next.

\bigskip

\begin{algorithm}[H]
{\bf Step 1:} Set $j=1$.

{\bf Step 2:} Do while $j \le K$:
\begin{itemize}
\item Obtain estimates $c_\mathcal{E}(t_{j})$, $\gamma_\mathcal{E}(t_{j})$ and  $T^{*}_{\mathcal{E}}(t_{j})$ by solving problem \eqref{eq:minf} over $I = I_j$ using initial guess $\w$.
\item Update the initial guess $\w =\left[c_\mathcal{E}(t_{j}),\gamma_\mathcal{E}(t_{j}),T^{*}_{\mathcal{E}}(t_{j})\right]$.
\item Set $j=j+1$ (slide the window).
\end{itemize}

{\bf Step 3:} Return $T^{*}_{\mathcal{E}}(t_{M})$. 

\caption{Approximation of the blow-up time $T^{*}$ based on the time evolution of $\E(t)$. \newline
{\bf Input:} $\mathcal{E}(t_{i})$, $i=1,\dots,Q$, where $Q$ is the number of discrete time steps. \newline
\hspace*{1.1cm} A family of sliding time windows $I_j,$ $j=1,\dots,K$.  \newline 
\hspace*{1.1cm} An initial guess $\w$ for the parameters $[c_\E, \gamma_\E, T^*_\E ]$.  \newline 
{\bf Output:} An approximation of $T^{*}_{\mathcal{E}}$.
}
\label{alg:Ts}
\end{algorithm}

\section{Results --- Deterministic Case}
\label{sec:results_deter}

In this section we analyze solutions of system \eqref{FBE} obtained in
the subcritical and supercritical regimes, focusing on how the blow-up
time in the latter case varies with problem parameters.  Unless
indicated otherwise, the initial condition and the viscosity are given
by $g(x)=\sin(x)$ and $\nu=0.11$.  System \eqref{FBE} is solved using
the numerical approach described in Section \ref{sec:numerdet} with
adaptive resolution varying from $N=2^{9}$ to $N=2^{18}$ in the
supercritical case and from $N=2^{9}$ to $N=2^{22}$ for the
subcritical case. However, in this latter case it is not always
necessary to use the highest resolution.

\begin{figure}[h]
\centering
\mbox{
\subfigure[]{\includegraphics[width=0.482\textwidth]{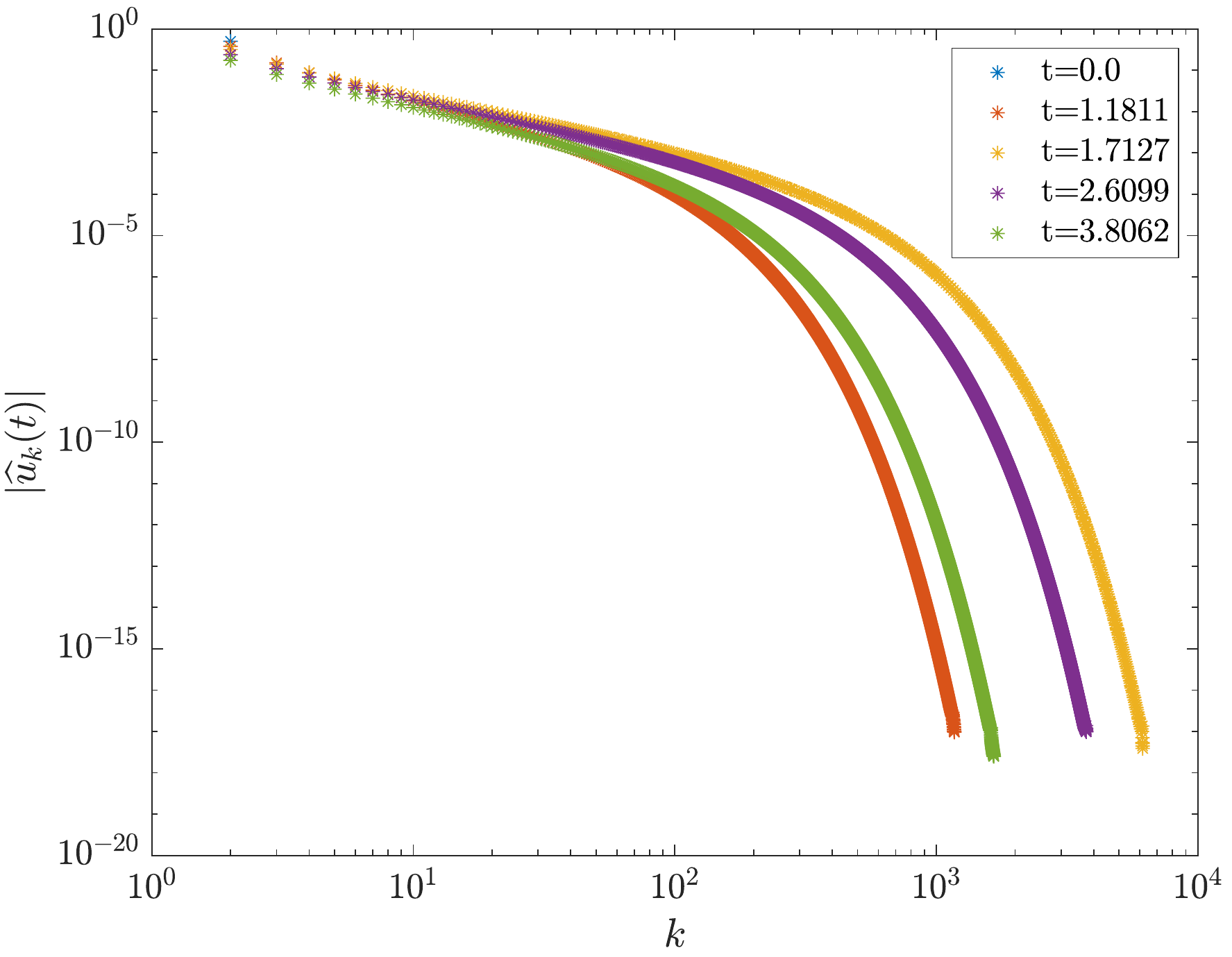}}
\subfigure[]{\includegraphics[width=0.495\textwidth]{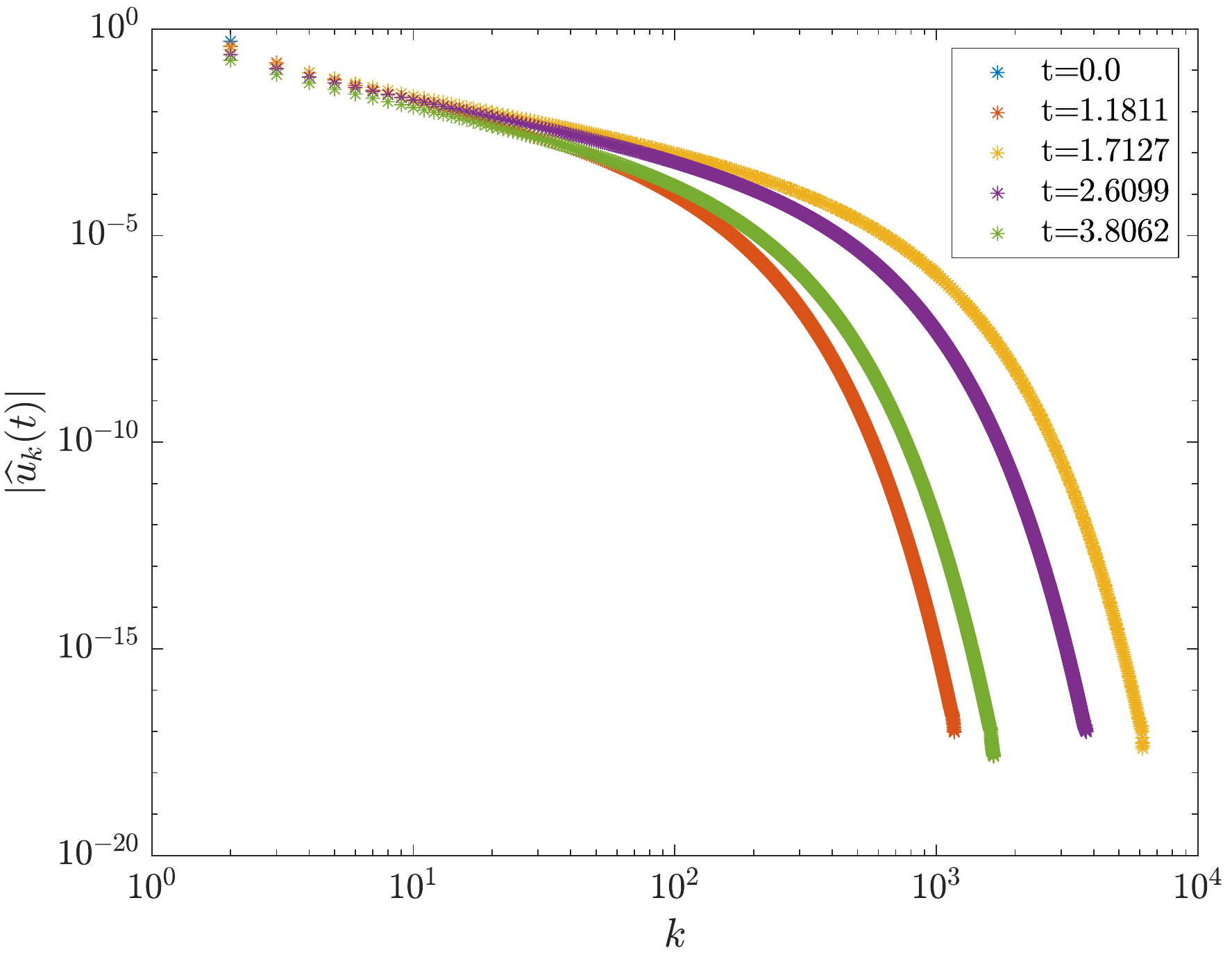}}}
\mbox{
\subfigure[]{\includegraphics[width=0.48\textwidth]{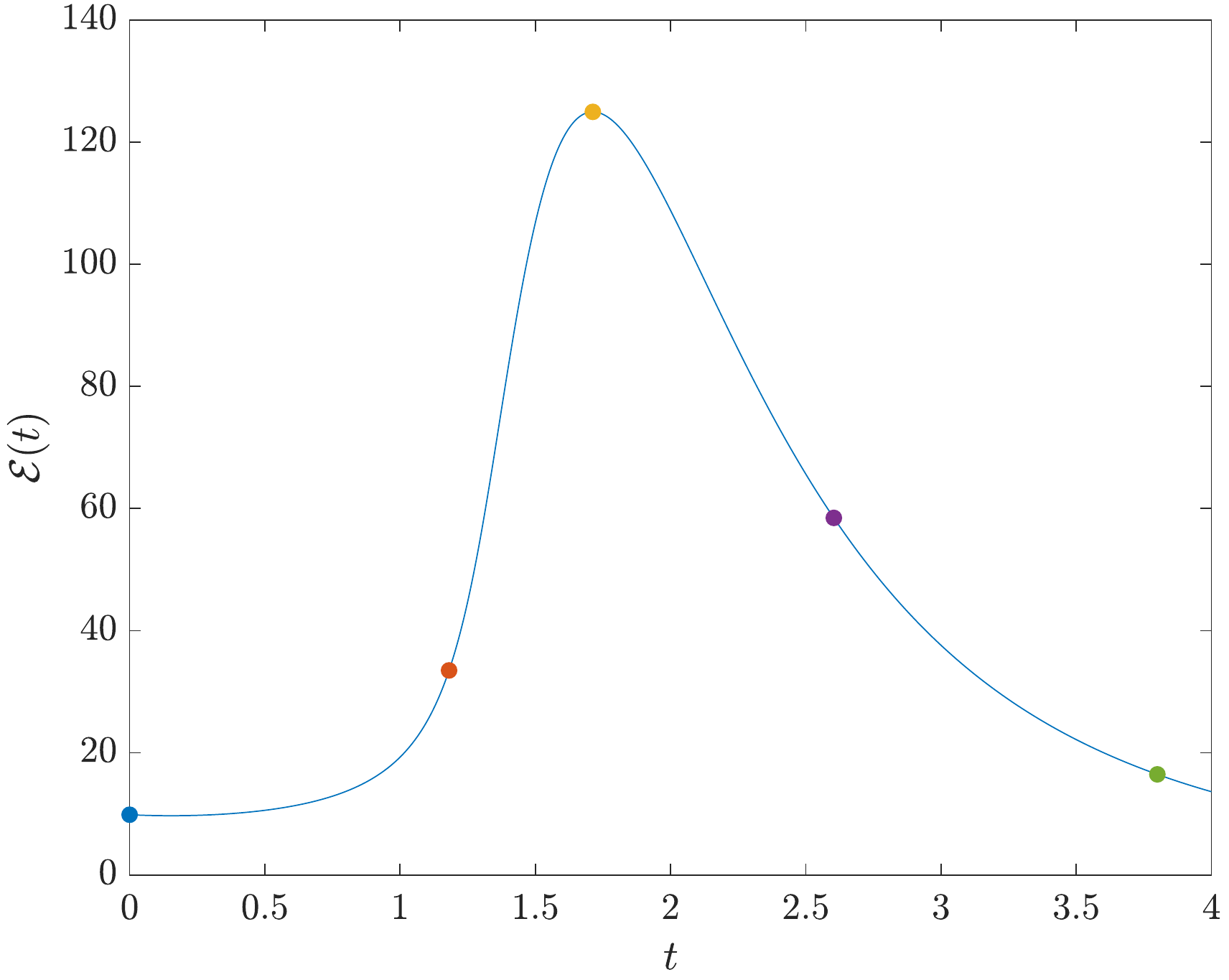}}
\subfigure[]{\includegraphics[width=0.484\textwidth]{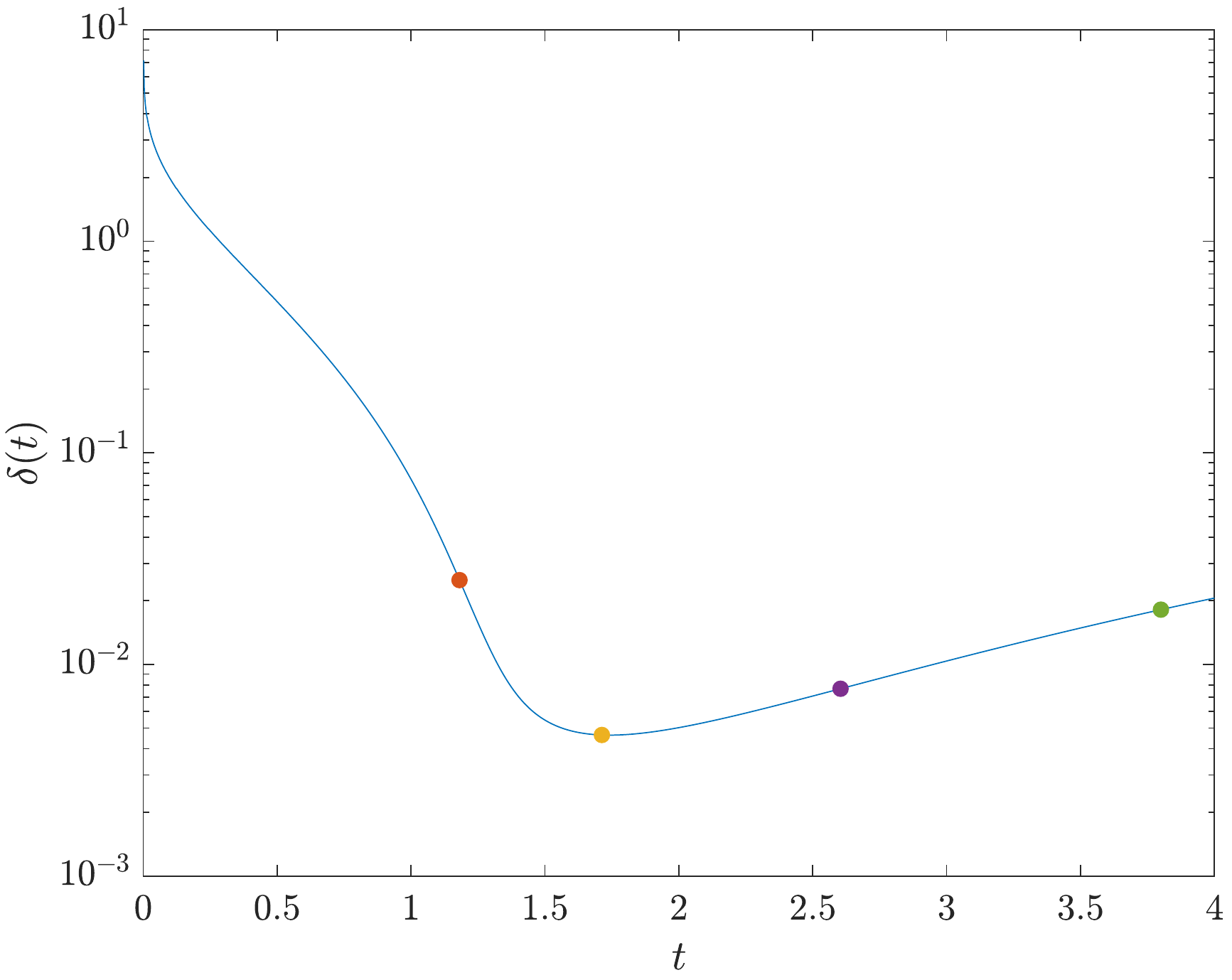}}}
\caption{Solution of system \eqref{FBE} with {$\alpha=0.6$} in (a) the
  physical space $u(t,x)$ and (b) the Fourier space
  $|\widehat{u}_{k}(t)|$ at the indicated time levels with the
  corresponding evolution of (c) the enstrophy $\mathcal{E}(t)$ and
  (d) the width of the analyticity strip $\delta(t)$. The symbols in
  panel (c) and (d) correspond to the time instances at which the
  solution is shown in panels (a) and (b).}
\label{fig:a0.6}
\end{figure}

\begin{figure}[h]
\centering
\mbox{
\subfigure[]{\includegraphics[width=0.482\textwidth]{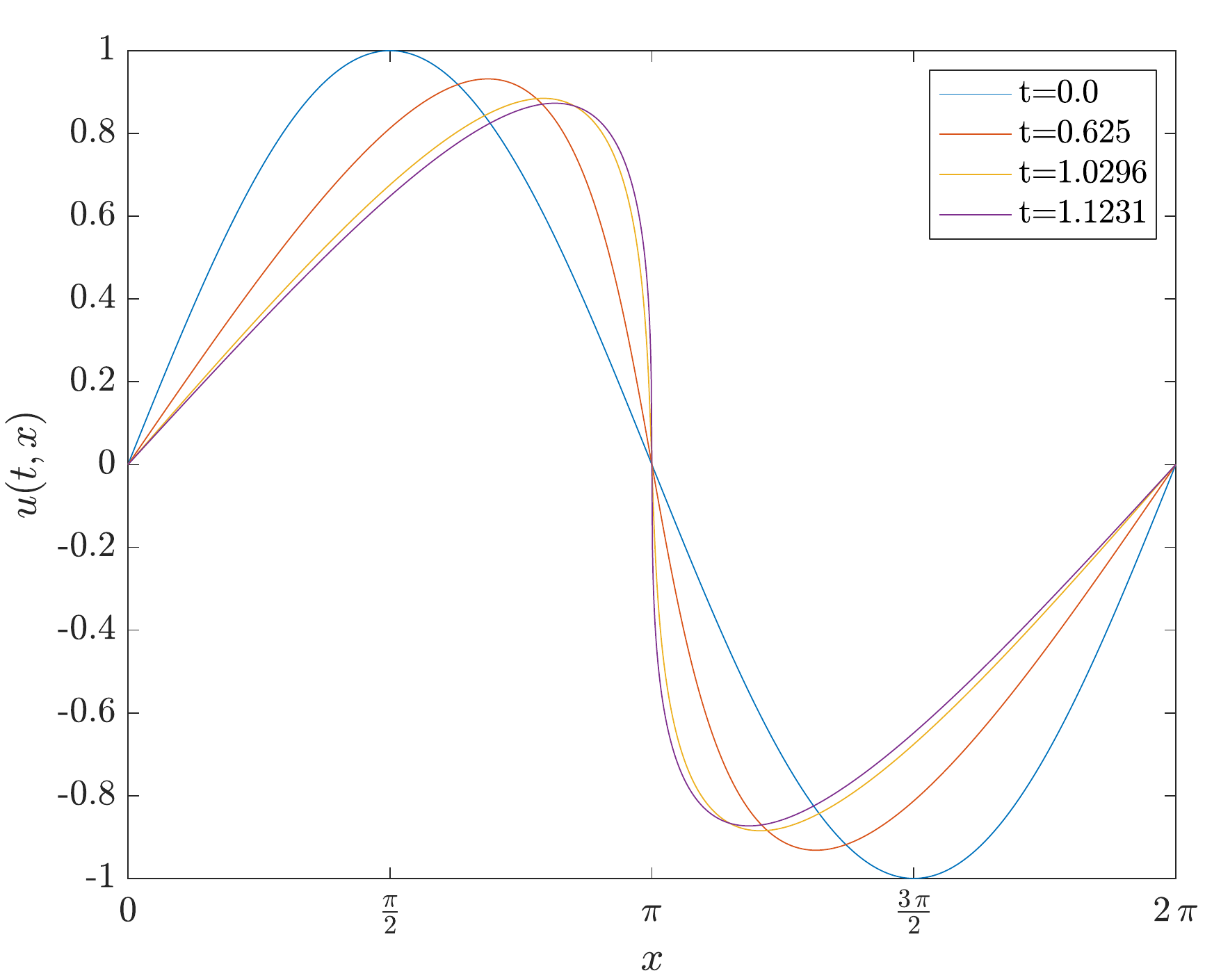}}
\subfigure[]{\includegraphics[width=0.495\textwidth]{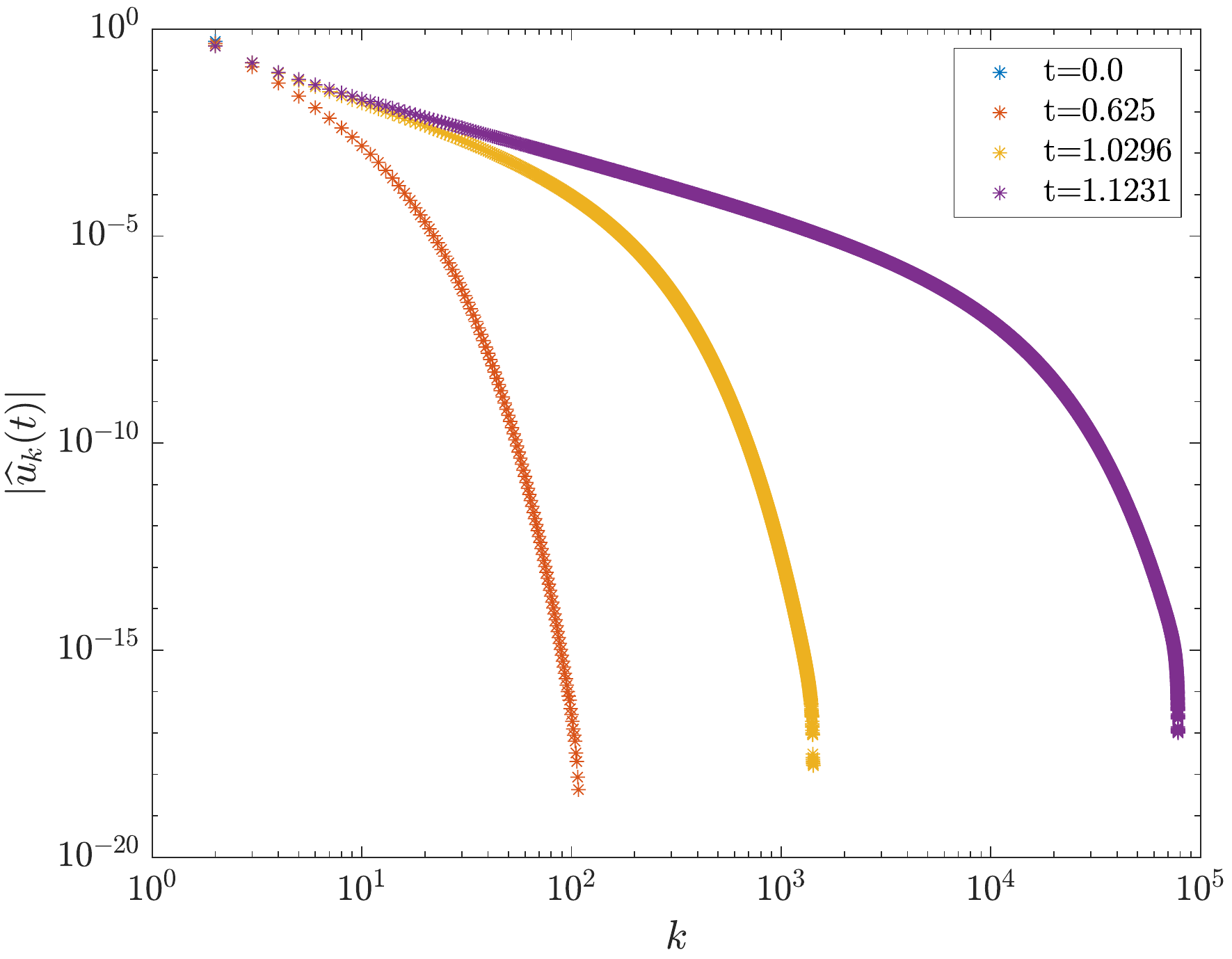}}}
\mbox{
\subfigure[]{\includegraphics[width=0.48\textwidth]{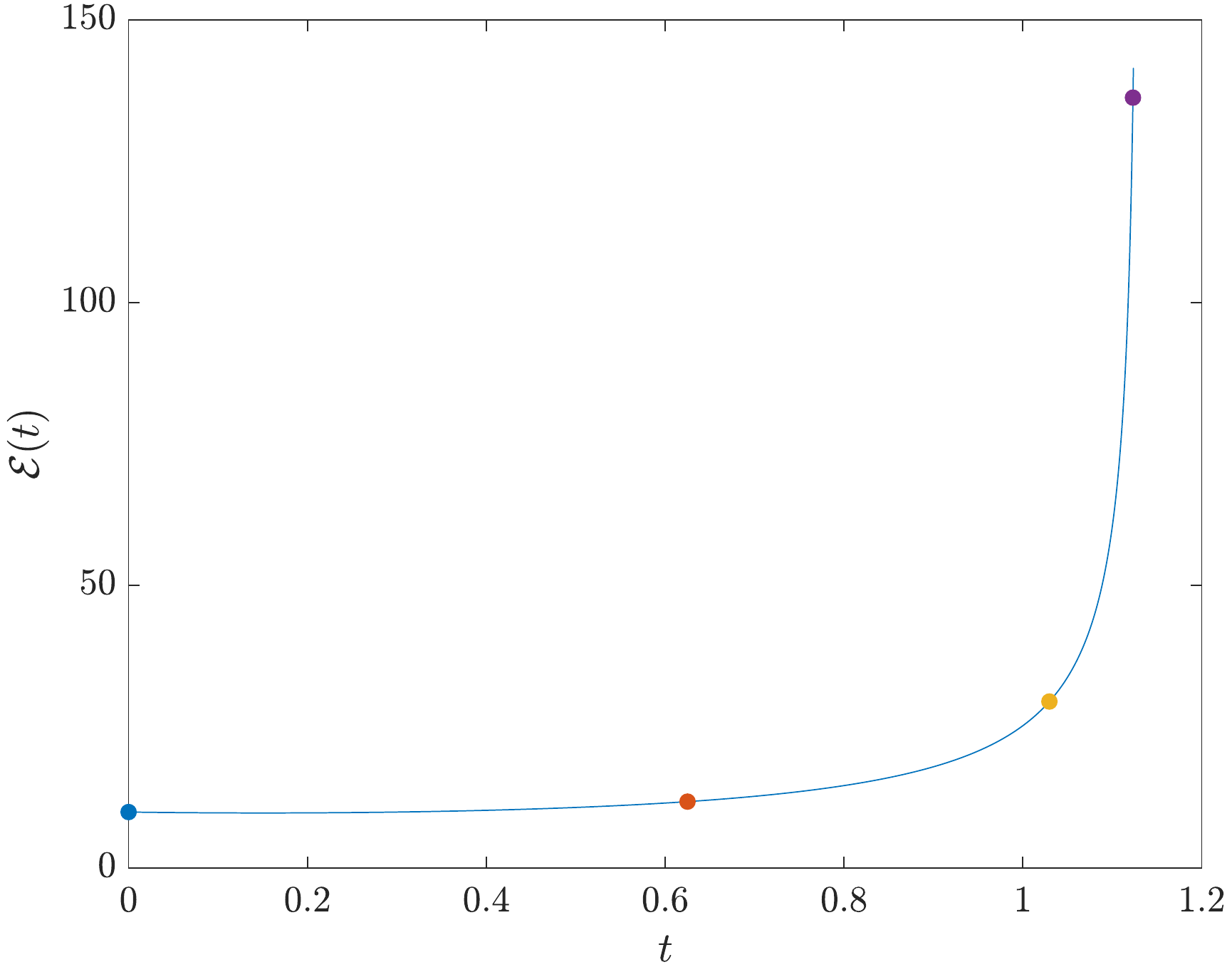}}
\subfigure[]{\includegraphics[width=0.482\textwidth]{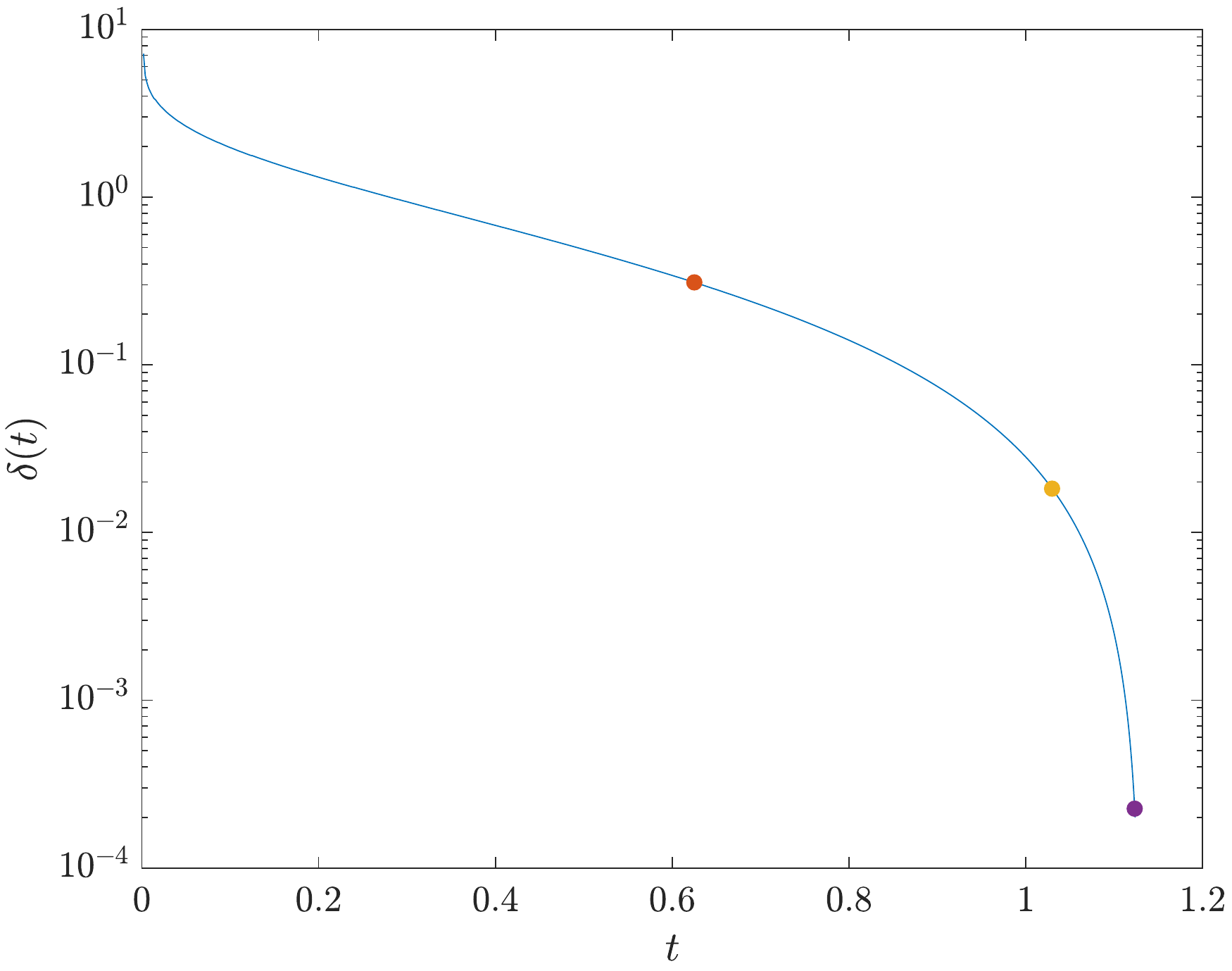}}}
\caption{Solution of system \eqref{FBE} with {$\alpha=0.4$} in (a) the
  physical space $u(t,x)$ and (b) the Fourier space
  $|\widehat{u}_{k}(t)|$ at the indicated time levels with the
  corresponding evolution of (c) the enstrophy $\mathcal{E}(t)$ and
  (d) the width of the analyticity strip $\delta(t)$. The symbols in
  panel (c) and (d) correspond to the time instances at which the
  solution is shown in panels (a) and (b).}
\label{fig:a0.4}
\end{figure}

\begin{figure}[h]
\mbox{
\subfigure[]{\includegraphics[width=0.482\textwidth]{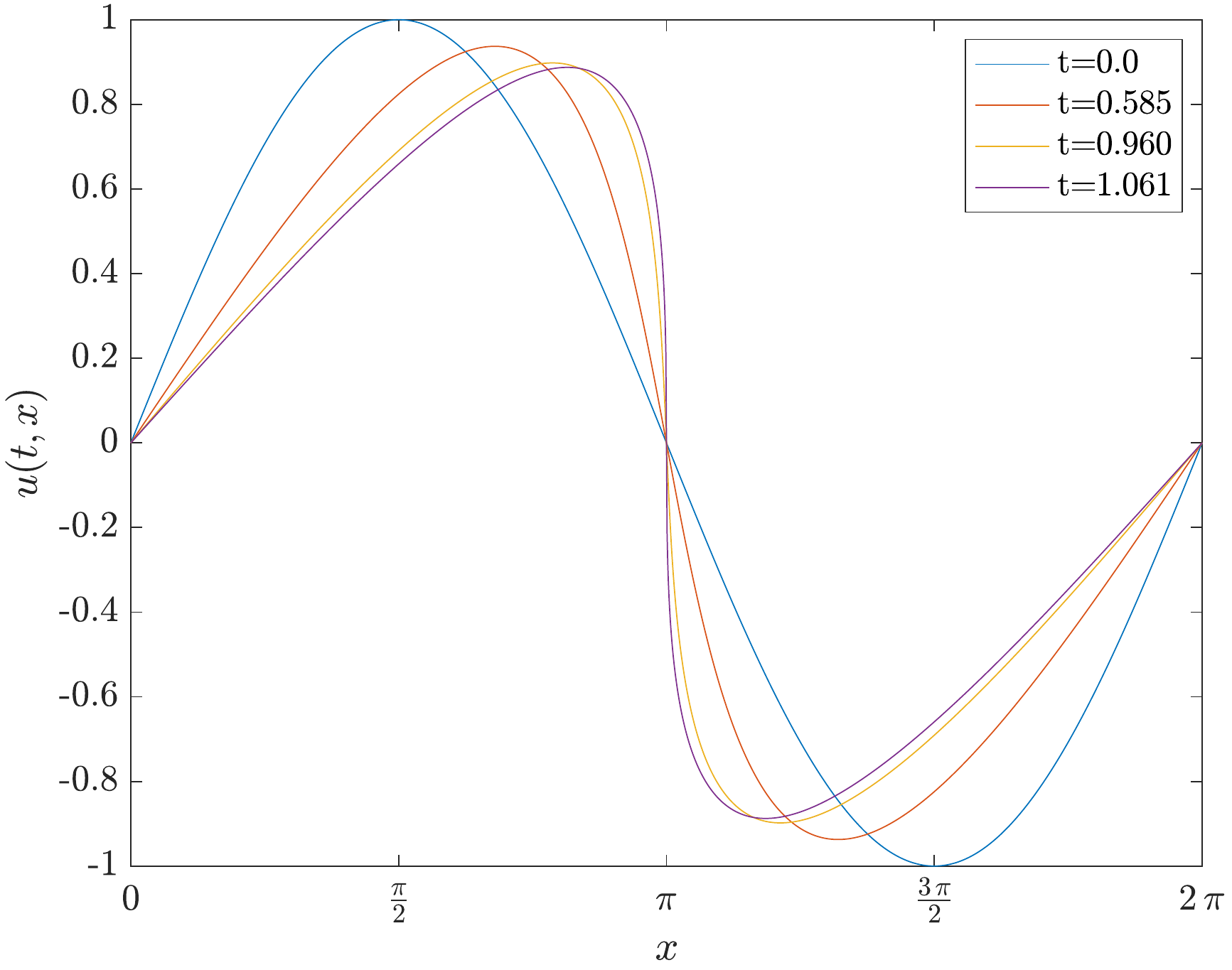}}
\subfigure[]{\includegraphics[width=0.495\textwidth]{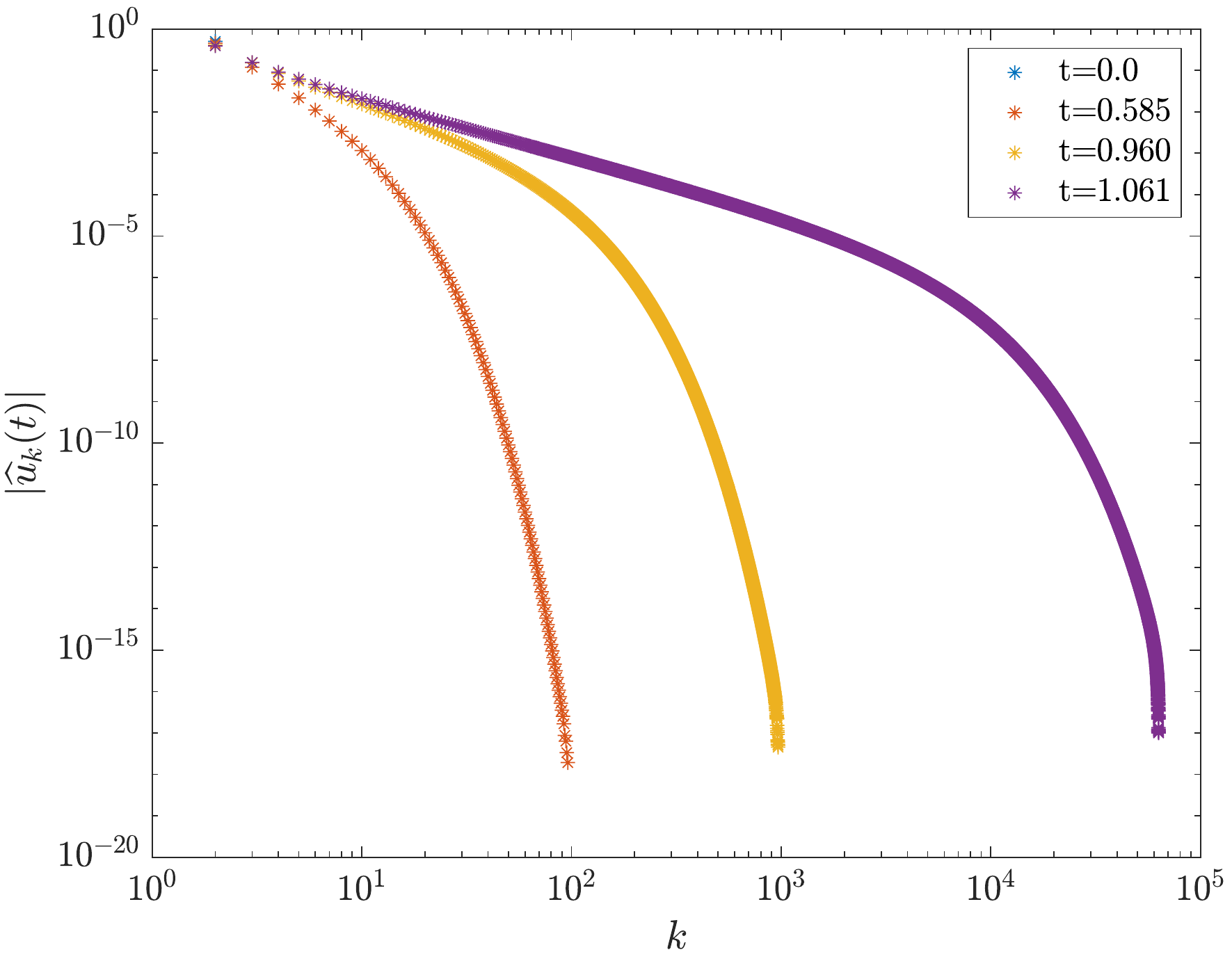}}}
\mbox{
\subfigure[]{\includegraphics[width=0.48\textwidth]{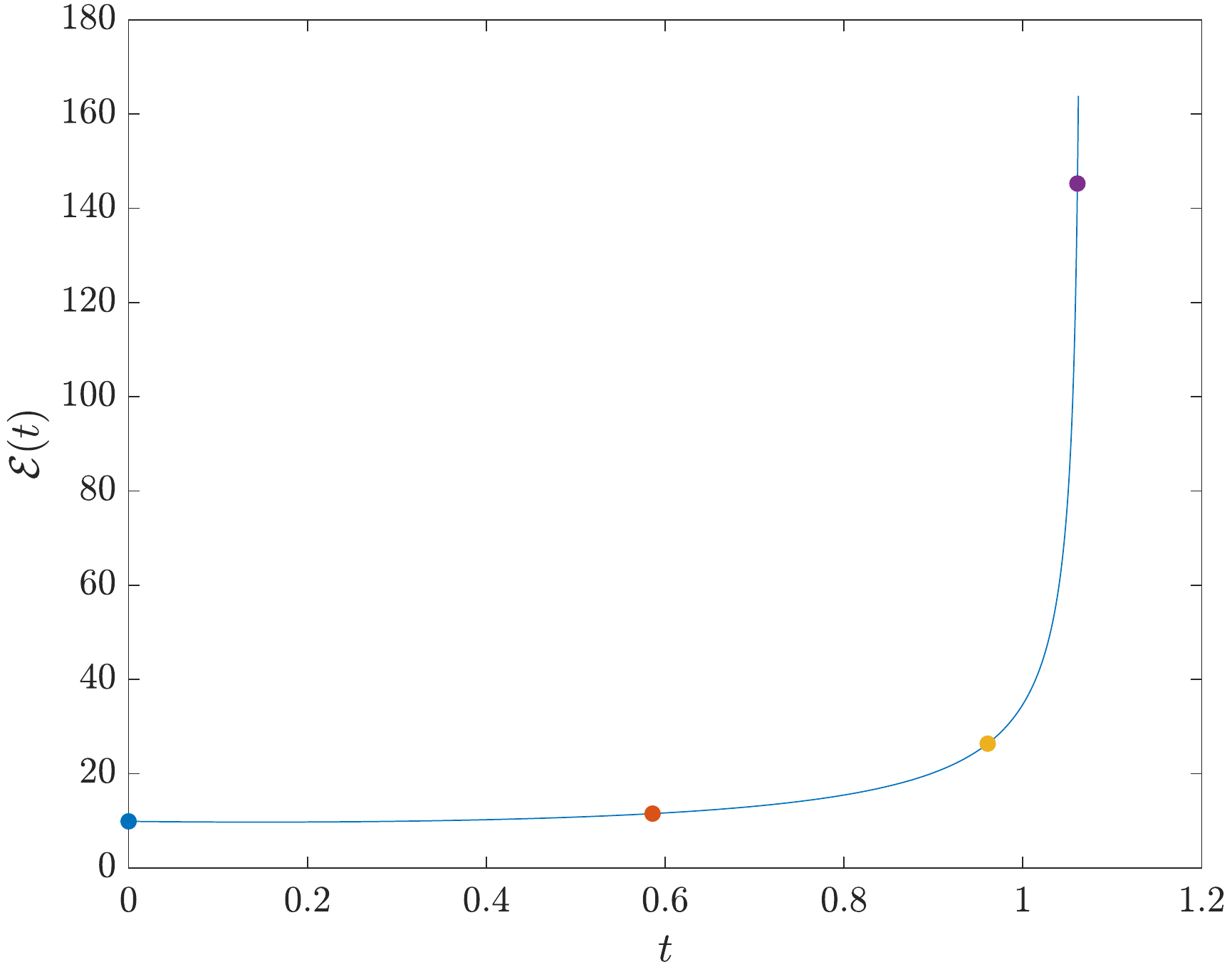}}
\subfigure[]{\includegraphics[width=0.482\textwidth]{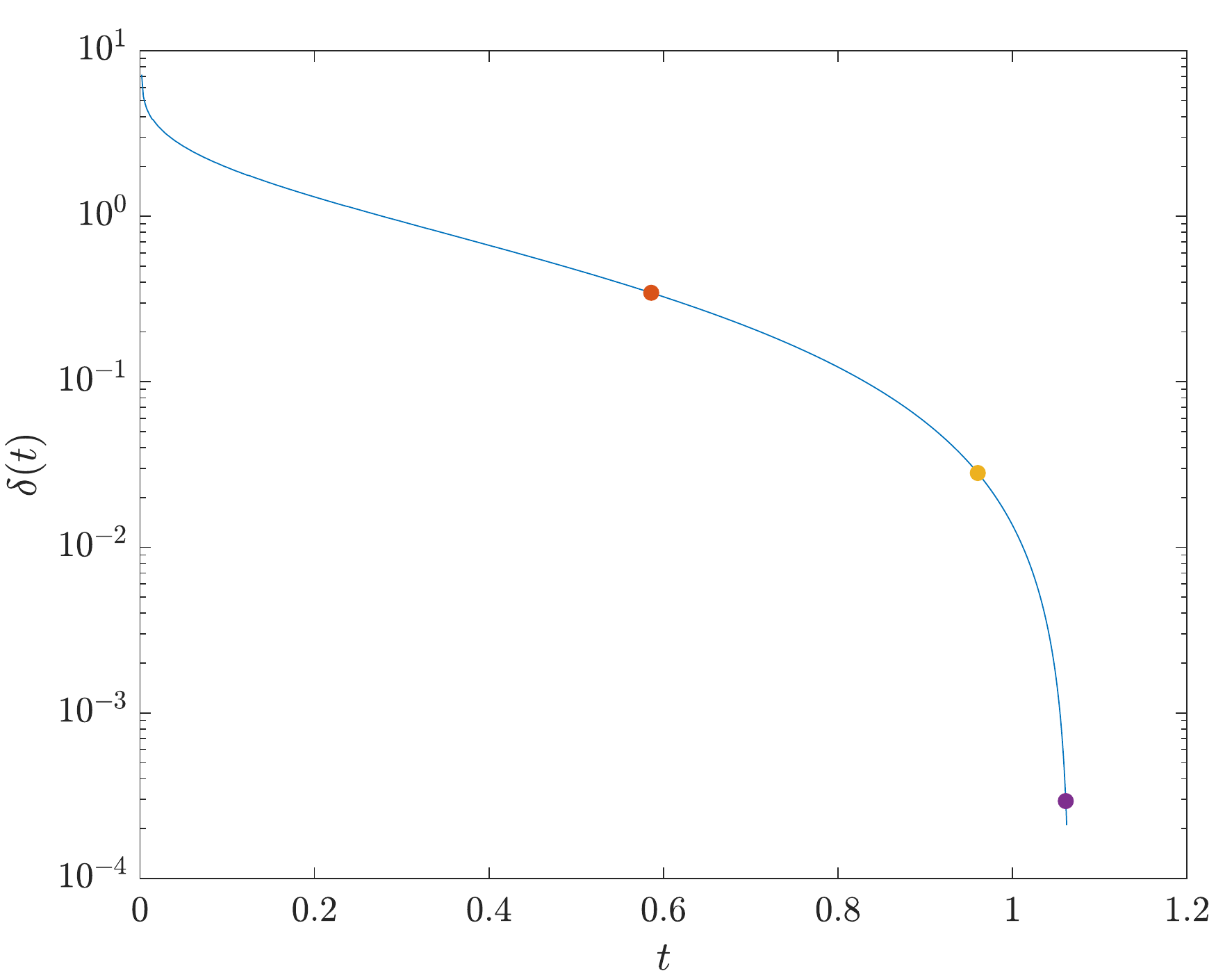}}}
\caption{Solution of system \eqref{FBE} with $\alpha=0.1$ in (a) the
  physical space $u(t,x)$ and (b) the Fourier space
  $|\widehat{u}_{k}(t)|$ at the indicated time levels with the
  corresponding evolution of (c) the enstrophy $\mathcal{E}(t)$ and
  (d) the width of the analyticity strip $\delta(t)$. The symbols in
  panel (c) and (d) correspond to the time instances at which the
  solution is shown in panels (a) and (b).}
\label{fig:a0.1}
\end{figure}

Results obtained by solving system \eqref{FBE} with ${\alpha = 0.6,
  0.4}, 0.1$ are presented in Figures \ref{fig:a0.6}, \ref{fig:a0.4}
and \ref{fig:a0.1} where we show the solutions $u(t,x)$ in the
physical space and in the spectral (Fourier) space at different time
levels as well as the time evolution of the enstrophy $\E(t)$ and of
the width of the analyticity strip $\delta(t)$. As regards the
subcritical case with {$\alpha = 0.6$}, in Figure \ref{fig:a0.6}a we
see that around the time $t \approx 1.71$, the solution (marked with a
yellow line) develops a steep front which is then smeared out by
dissipation. This behaviour is also evident in the evolution of the
Fourier spectrum $|\hu_k(t)|$ which at large times $t$ exhibits a
faster decay with the wavenumber $k$, cf.~Figure \ref{fig:a0.6}b.
Since the solution remains smooth (analytic) throughout its entire
evolution, the enstrophy $\E(t)$ and the width of the analyticity
strip $\delta(t)$ remain bounded, respectively, from above and below.
These quantities achieve their maximum and minimum values at the time
when the front in the solution is the most steep, cf.~Figures
\ref{fig:a0.6}c--d.

In the supercritical cases with {$\alpha = 0.4$} and $\alpha = 0.1$
presented in Figures \ref{fig:a0.4} and \ref{fig:a0.1} the
computations are stopped at some point before the singularity occurs
at $t = T^*$, which is caused by a rapid increase of the computational
cost. In Figures \ref{fig:a0.4}a--b we see that at the singularity
time the derivative of the solution at $x = \pi$ becomes unbounded, a
blow-up mechanism which is also known from the inviscid Burgers
problem \eqref{eq:iBE} \cite{kl04}. The singularity formation is also
manifested by the enstrophy becoming unbounded and the width of the
analyticity strip vanishing as $t \rightarrow (T^*)^-$, cf.~Figures
\ref{fig:a0.4}c--d and \ref{fig:a0.1}c--d. As regards the Fourier
spectra, in Figures \ref{fig:a0.4}b and \ref{fig:a0.1}b we observe
that they extend to higher and higher wavenumbers as the blow-up time
is approached and the region where the Fourier coefficients $\hu_k(t)$
decay exponentially fast with $k$ shrinks. This behavior marks the
transition from an exponential to algebraic decay of the spectrum
$|\hu_k(t)|$ with $k$ as the solution develops a singularity.

\begin{figure}[h]
\mbox{
\subfigure[]{\includegraphics[width=0.48\textwidth,trim={1.4cm 6.8cm 2.3cm 7.3cm},clip]{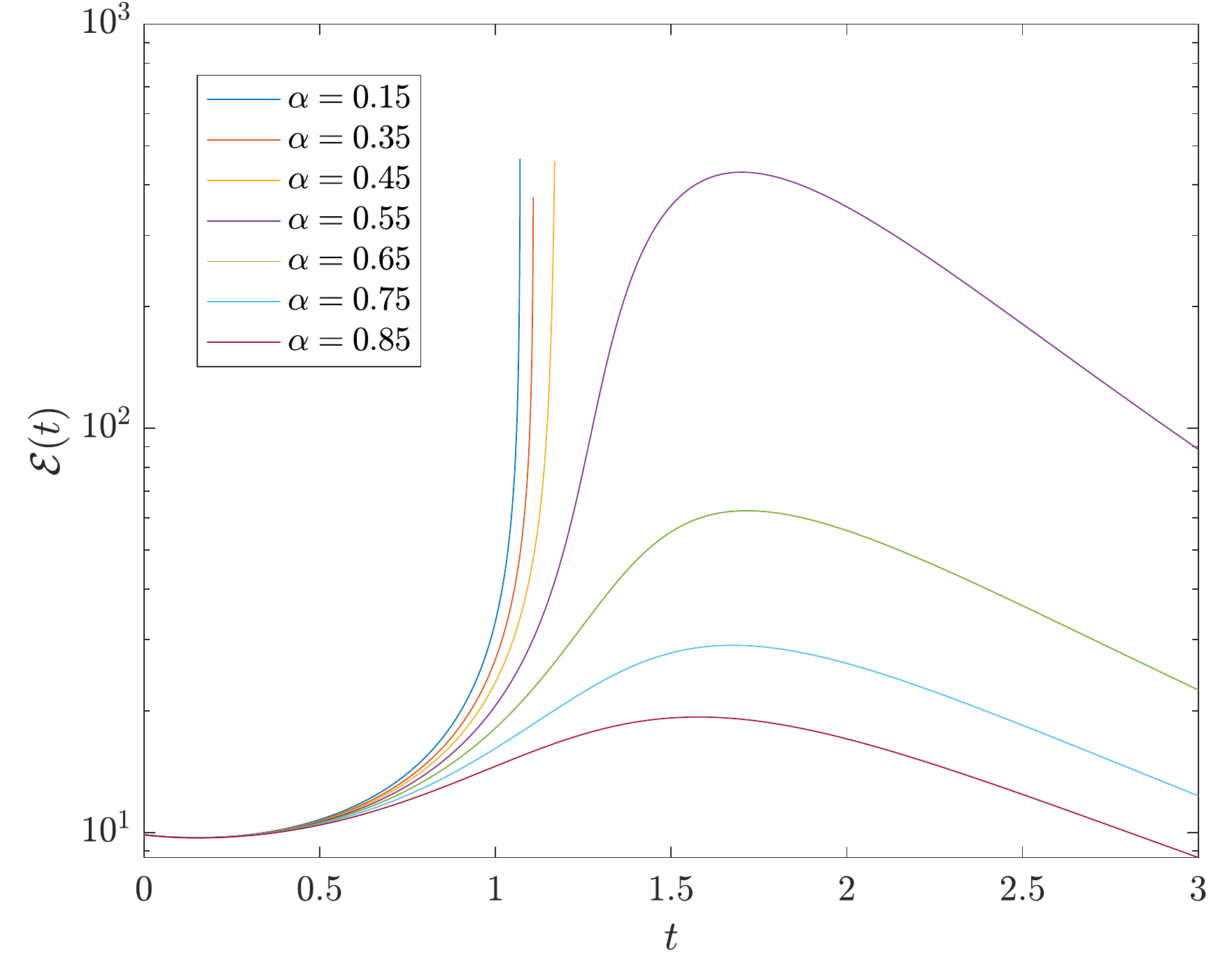}}
\subfigure[]{\includegraphics[width=0.48\textwidth,trim={1.4cm 6.8cm 2.3cm 7.3cm},clip]{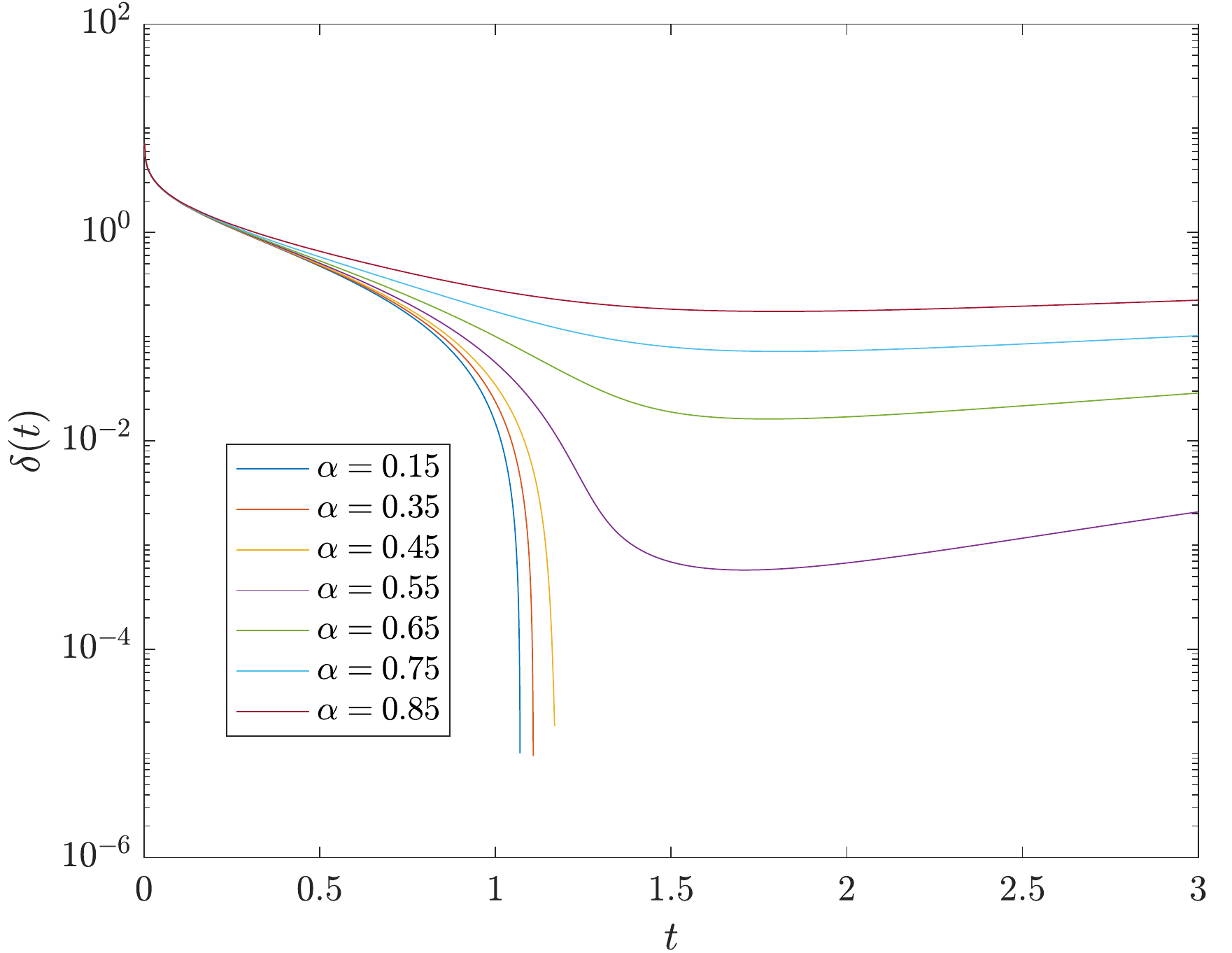}}}
\caption{Time evolution of (a) the enstrophy $\E(t)$ and (b) the width
  of the analyticity strip $\delta(t)$ in solutions of the fractional
  Burgers equation \eqref{FBE} with different indicated values of the
  fractional dissipation exponent $\alpha$.}
\label{fig:Ed}
\end{figure}

\begin{figure}[h]
\mbox{
\subfigure[]{\includegraphics[width=0.48\textwidth,trim={1.60cm 6.8cm 2cm 7.3cm},clip]{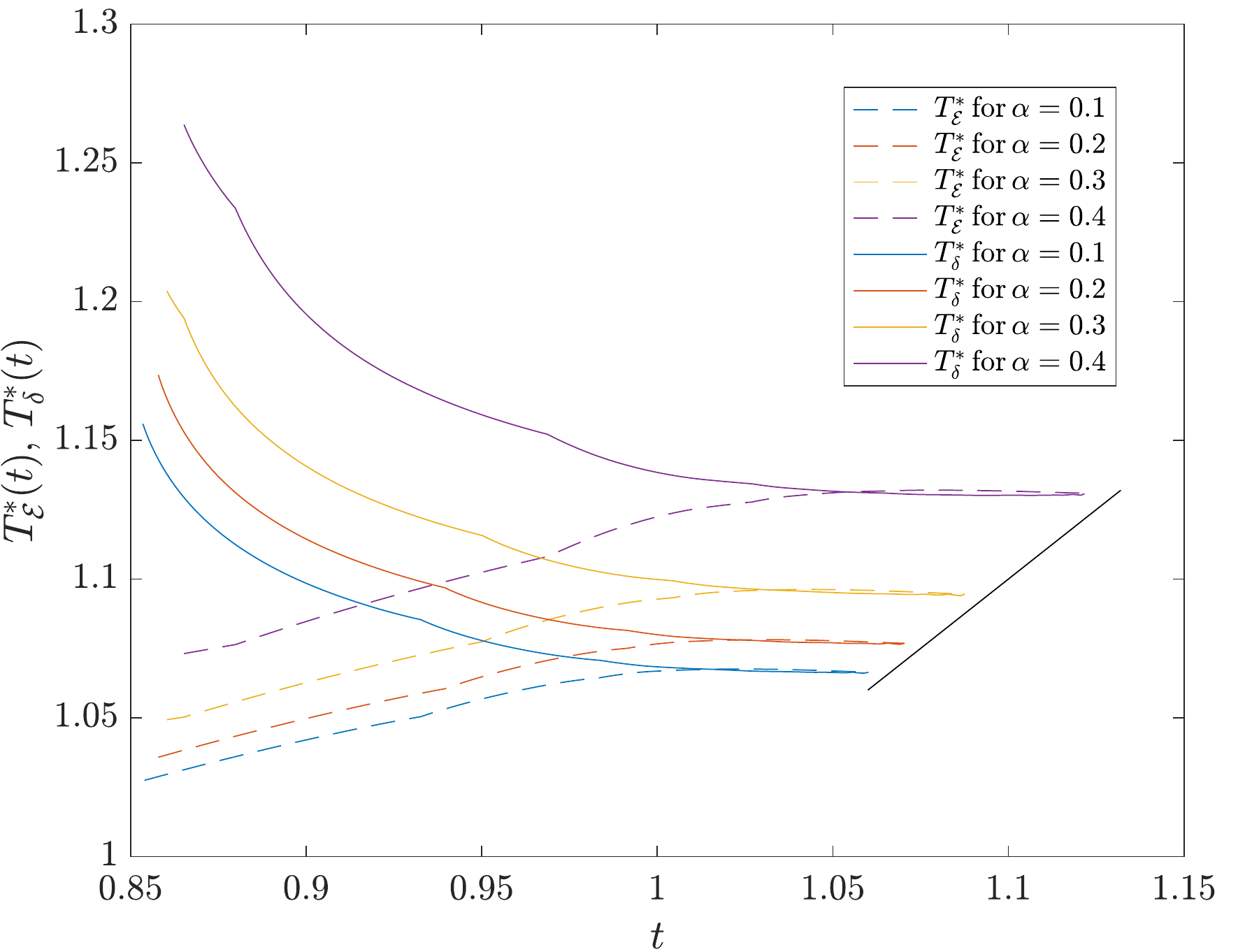}}
\subfigure[]{\includegraphics[width=0.46\textwidth]{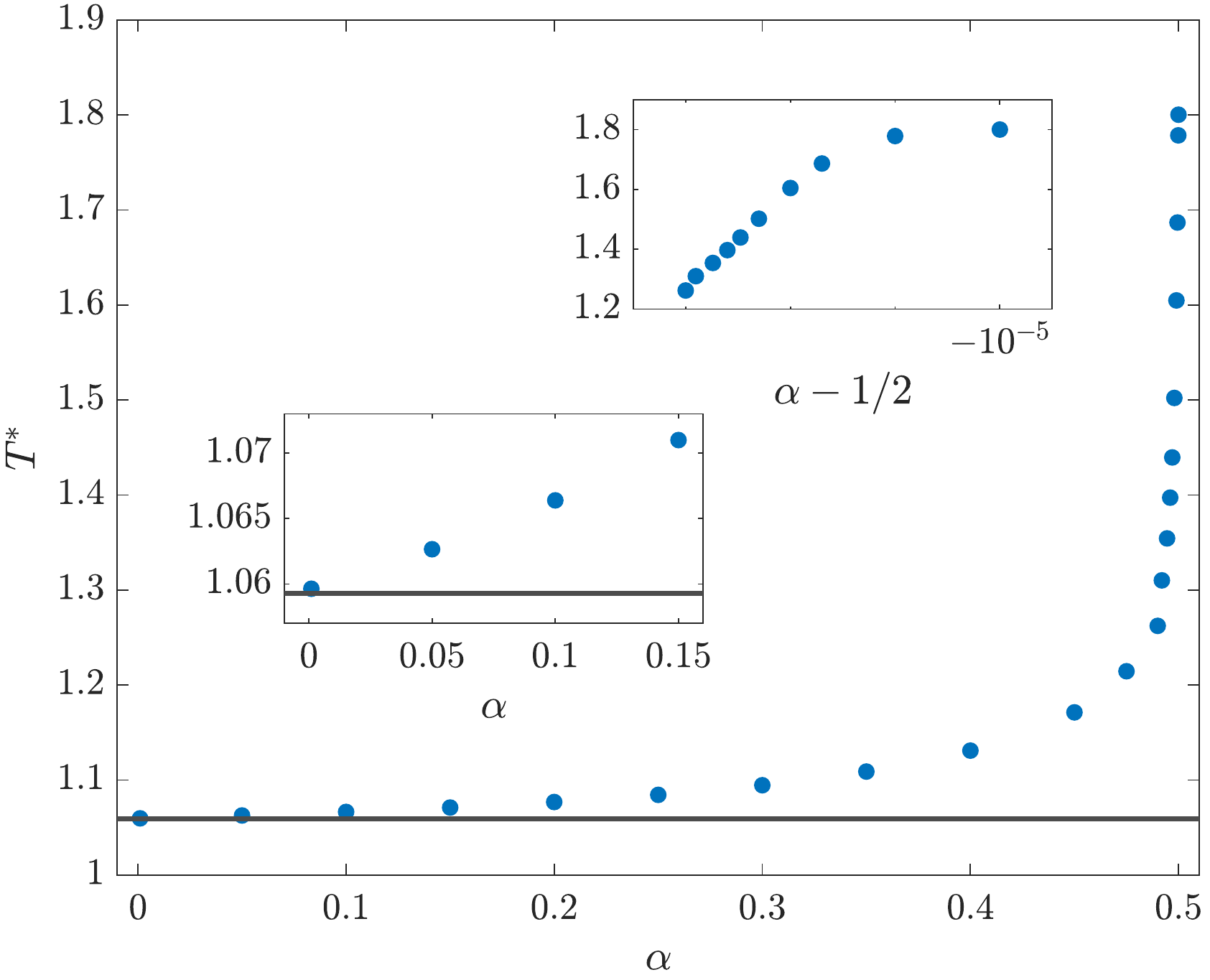}}}
\caption{(a) Estimates of the blow-up time based on the enstrophy
  $T^{*}_{\mathcal{E}}$ (dashed lines) and on the width of the
  analyticity strip $T^{*}_{\delta}$ (solid lines) as functions of
  time for different values of $\alpha$. The time $t$ is understood
  here as the position of the center of the window $I_j$,
  $j=1,\dots,K$, over which fits are performed when solving
  optimization problem \eqref{eq:minf}, cf.~Algorithm \ref{alg:Ts}.
  Curves with the same color indicate the estimates obtained for the
  same value of $\alpha$. The black slanted line corresponds to
  $T^{*}=t$ and serves as a reference for how far the estimate of
  $T^{*}$ can be continued.  (b) Dependence the estimated blow-up
  times $T^*$ on the fractional dissipation exponent $\alpha$. The top
  and bottom inset represent small neighborhoods near the endpoints
  $\alpha = 1/2$ and $\alpha = 0$, whereas the horizontal line denotes
  the blow up time $T^*_1$ in the limiting case when
  $\alpha\rightarrow 0^{+}$, cf.~\eqref{eq:TsfBE}.}
\label{fig:Ts}
\end{figure}

The time evolution of the enstrophy $\E(t)$ and of the width of the
analyticity strip $\delta(t)$ are shown for several values of the
exponent $\alpha$ in both the subcritical and supercritical regime in
Figures \ref{fig:Ed}a--b (the critical case with $\alpha = 1/2$, where
system \eqref{FBE} is known to be globally well-posed, cf.~Theorem
\ref{thm:crit}, is omitted due to its excessive computational cost).
We see that as $\alpha \rightarrow (1/2)^+$ the maximum attained
values of the enstrophy $\max_t \E(t)$ increase and the minimum
attained values of the width of the analyticity strip $\min_t
\delta(t)$ decrease. When $\alpha$ decreases to zero in the
supercritical regime, blow-up occurs earlier as indicated by the
instances of time when the enstrophy becomes unbounded and the width
of the analyticity strip vanishes. In order to analyze this aspect
quantitatively, in Figure \ref{fig:Ts}a we show the blow-up times
$T_\E^*(t)$ and $T_\delta^*(t)$ estimated using Algorithm \ref{alg:Ts}
based on, respectively, the evolution of the enstrophy and of the
width of the analyticity strip as functions of the location of the
time window $I_j$ over which the optimization problem \eqref{eq:minf}
is solved. These results are presented in Figure \ref{fig:Ts}a for
different supercritical values of $\alpha \in (0,1/2)$. We see that,
interestingly, fits performed based on $\E(t)$ at early times
underestimate the blow-up time while those performed based on
$\delta(t)$ overestimate the blow-up time. However, as the windows
$I_j$ approach the blow-up time, the two estimates $T_\E^*(t)$ and
$T_\delta^*(t)$ agree with each other very well. In fact, the
difference between these two estimates is
$\mathcal{O}(10^{-3}-10^{-4})$ depending on $\alpha$, which allows us
to conclude that the two approaches to estimating the blow-up time are
consistent.  Hereafter, we will estimate the blow-up time based on the
evolution of the enstrophy as this approach can also be used in the
stochastic case (since solutions of the stochastic problem
\eqref{SBE2} are not in general analytic functions of the space
coordinate, we have $\delta(t) = 0$, $\forall t>0$ in this case). We
will therefore set $T^{*}=T^{*}_{\mathcal{E}}(t_{K})$.

The blow-up times $T^*$ determined as discussed above are shown as
functions of $\alpha \in (0,1/2)$ in Figure \ref{fig:Ts}b. We remark
that in order to explore behaviors corresponding to the limiting
values of this parameter, its smallest and largest values are,
respectively, $\alpha = 10^{-3}$ and $\alpha = 1/2 - 10^{-5}$. As is
evident from Figure \ref{fig:Ts}b, the blow-up time $T^*$ is an
increasing function of $\alpha$. In the limit $\alpha \rightarrow 0$
it tends to the value given by expression \eqref{eq:TsfBE} for the
limiting system \eqref{eq:fBE}. The data in Figure \ref{fig:Ts}b also
suggests that the blow-up time $T^*$ remains bounded in the opposite
limit when $\alpha \rightarrow (1/2)^-$.

\begin{figure}[h!]
\centering
\includegraphics[scale=0.49]{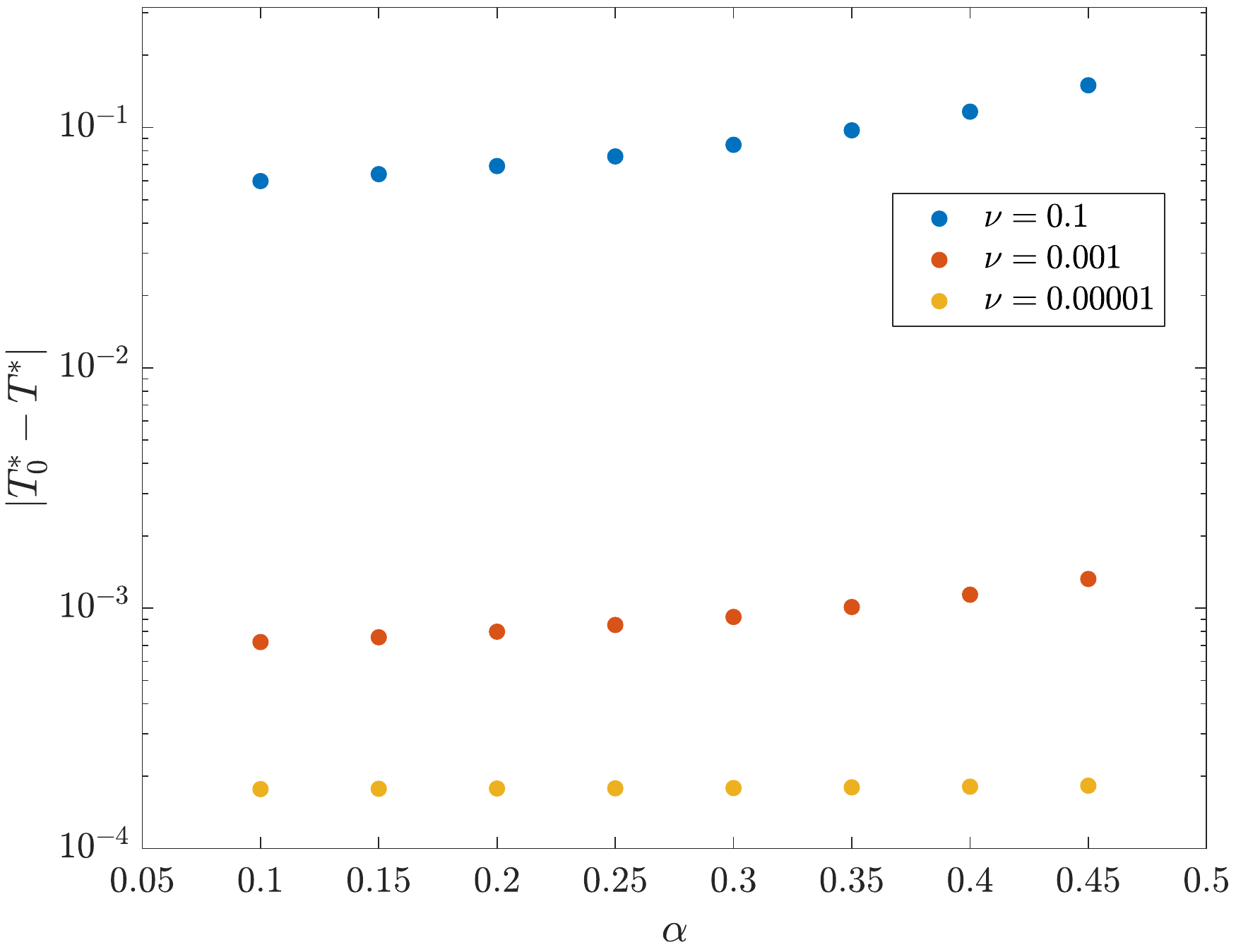}
\caption{The difference between the blow-up times $T^*$ characterizing
  solutions of the fractional Burgers system \eqref{FBE} with given
  $\alpha$ and $\nu$ and the blow-up time $T_0^* = 1$ in the solution
  of the inviscid problem \eqref{eq:iBE}, which is independent of
  $\alpha$, cf.~\eqref{eq:TsiBE}.}
\label{fig:Tnu0}
\end{figure}

Finally, we consider the behavior of the blow-up time $T^*$ in the
inviscid limit $\nu \rightarrow 0$ with fixed $\alpha \in (0,1/2)$.
The results shown in Figure \ref{fig:Tnu0} demonstrate that as $\nu$
vanishes the blow-up time $T^*$ approaches the blow-up time for the
inviscid problem \eqref{eq:iBE}, cf.~expression \eqref{eq:TsiBE},
uniformly in $\alpha$. In this figure we also observe that the
relative range of variation of $T^*$ is reduced as $\nu$ decreases.
Some additional details concerning the results presented in this
section are available in \cite{Ramirez2020}. In the next section we
discuss how the behavior of fractional Burgers flows is affected by
the stochastic forcing introduced in Section \ref{sec:BurgersStoch}.

\section{Results --- Stochastic Case}
\label{sec:results_stoch}

In this section we present results obtained for the stochastic
fractional Burgers system \eqref{SBE2}, first in the supercritical
regime (with $\alpha = 0.4$) and then in the subcritical regime (with
$\alpha = 0.6$). We use the same initial condition and same value of
the viscosity coefficient as in the deterministic case, i.e., $g(x) =
\sin(x)$ and $\nu = 0.11$. The problem is solved using the approach
described in Section \ref{sec:numerstoch} with the spatial resolution
given by $N = 2^{17} = 131,072$ and $N = 2^{15}=32,768$ in the
supercritical and subcritical regimes, respectively. The number of
Monte-Carlo samples $M$ used to discretize the probability space will
be specified individually in the different cases below.

\subsection{Supercritical Regime}
\label{sec:results_super}

\begin{figure}[h]
\mbox{
\subfigure[]{\includegraphics[width=0.48\textwidth]{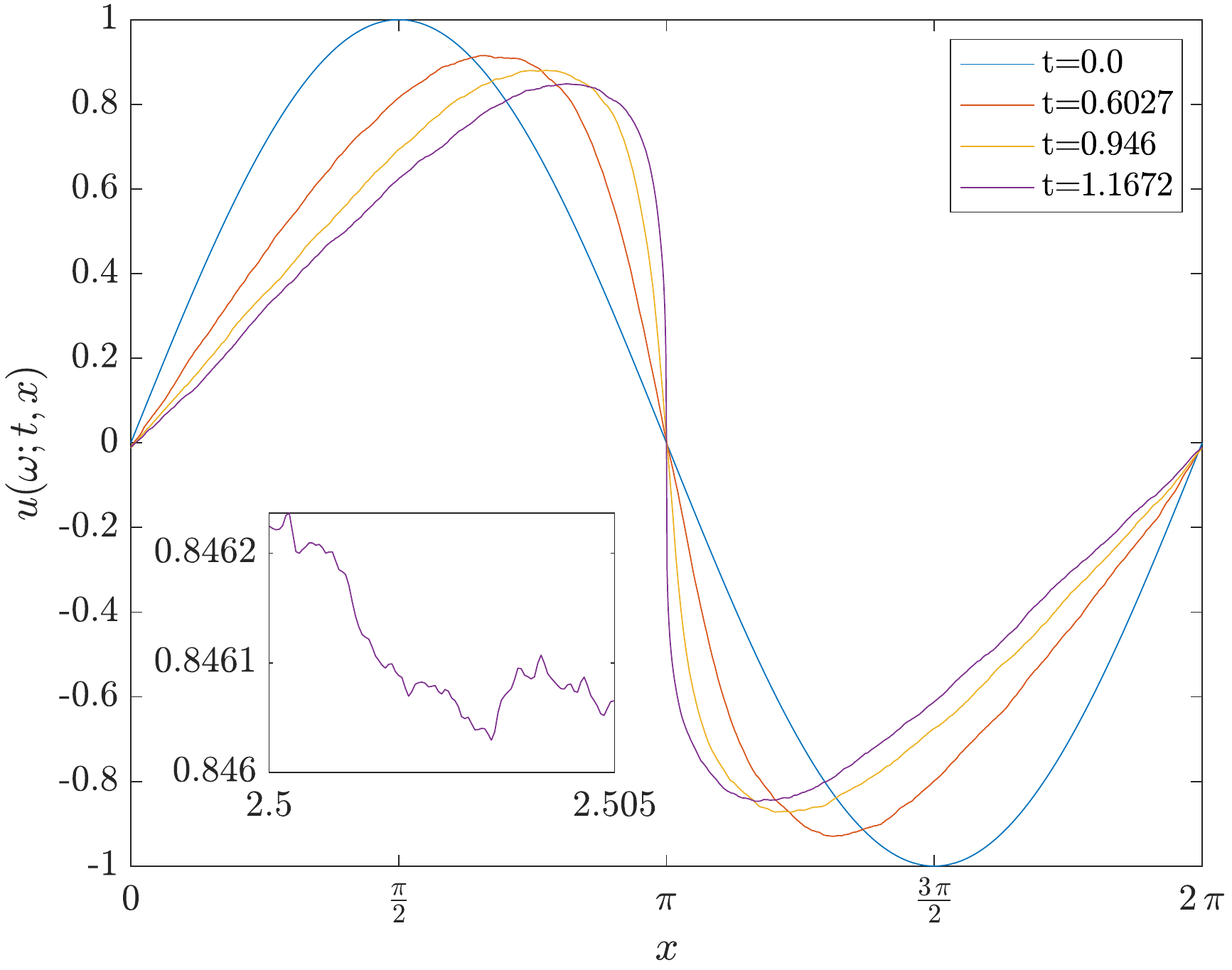}}
\subfigure[]{\includegraphics[width=0.48\textwidth]{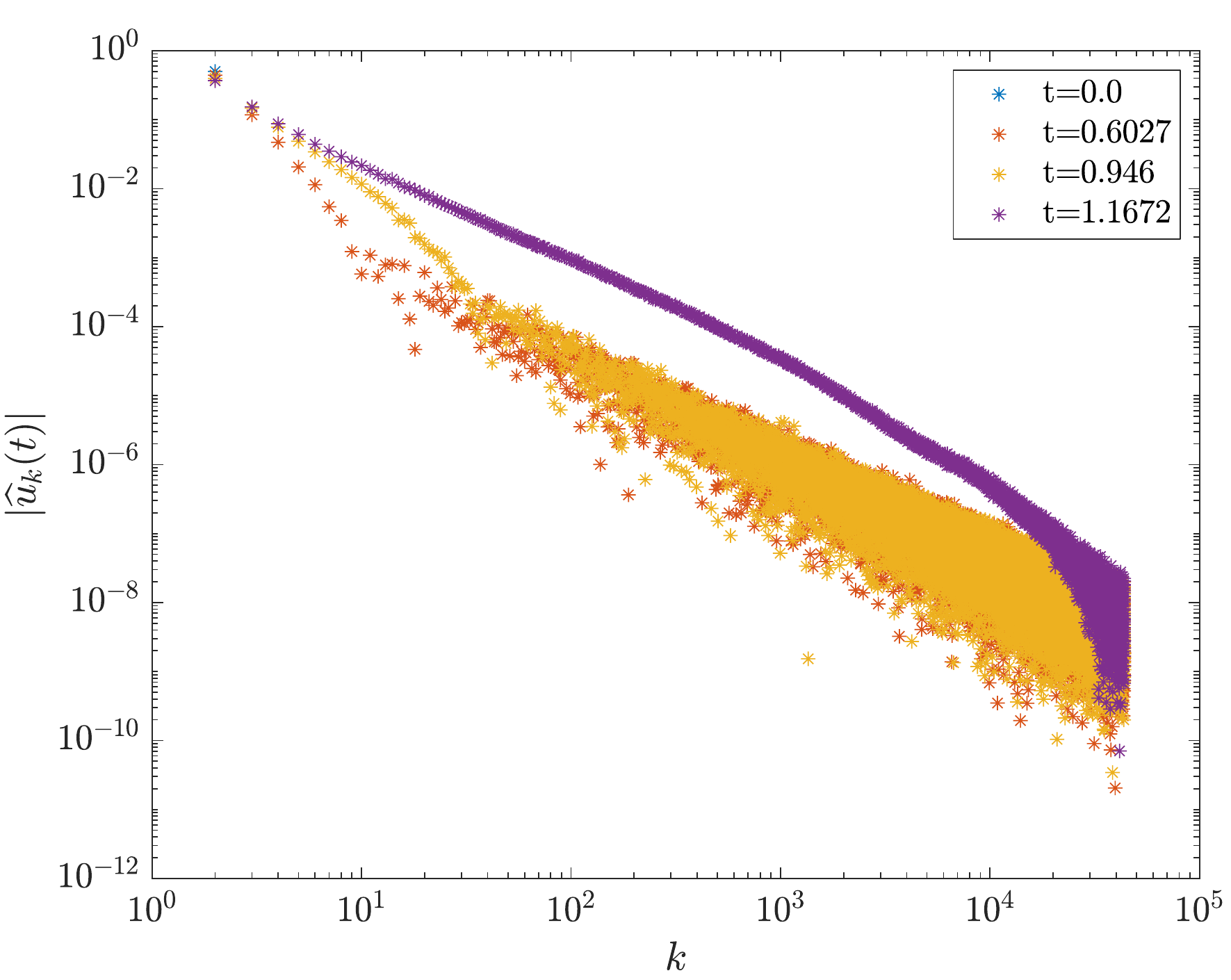}}}
\mbox{
\subfigure[]{\includegraphics[width=0.48\textwidth]{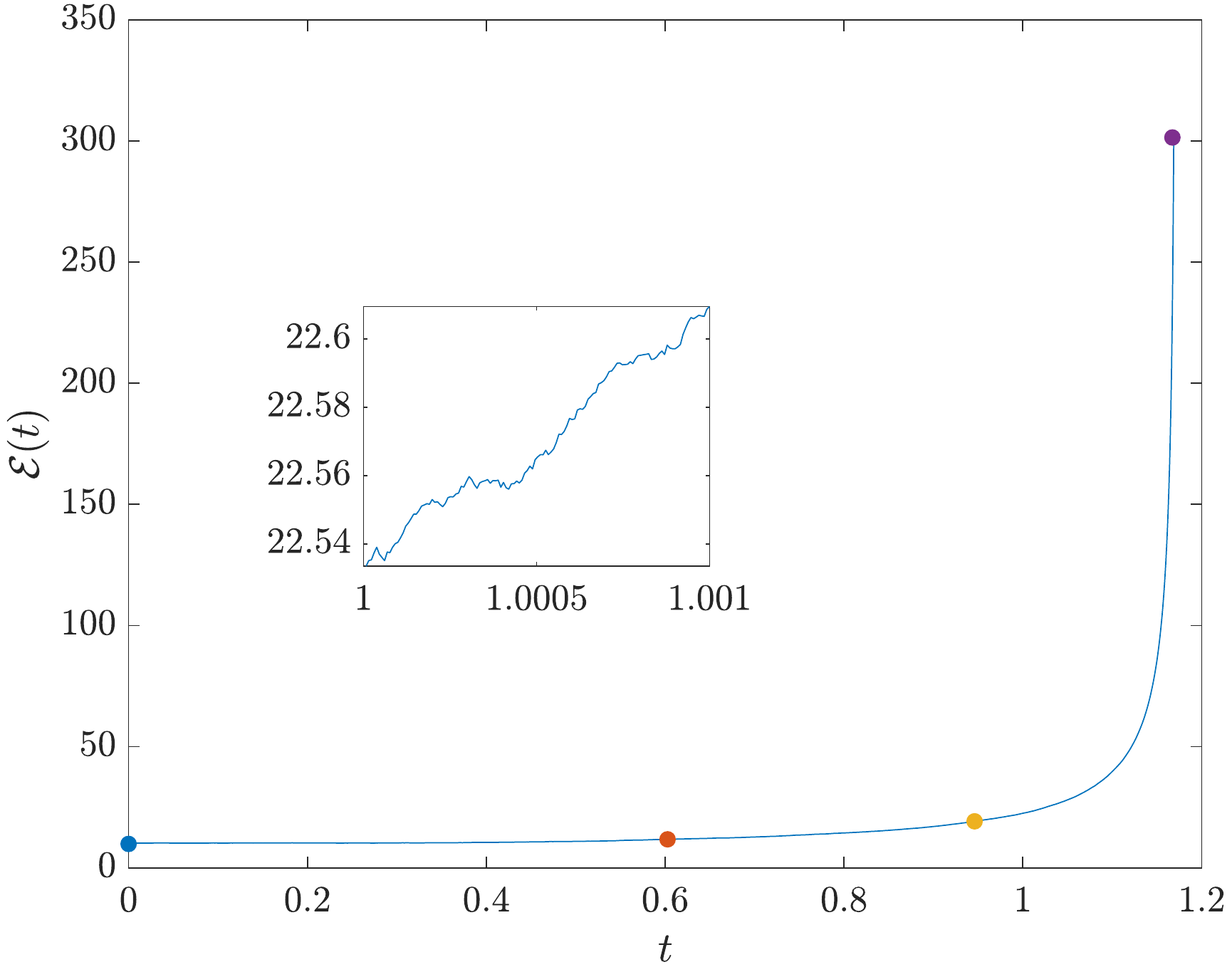}}
\subfigure[]{\includegraphics[width=0.48\textwidth]{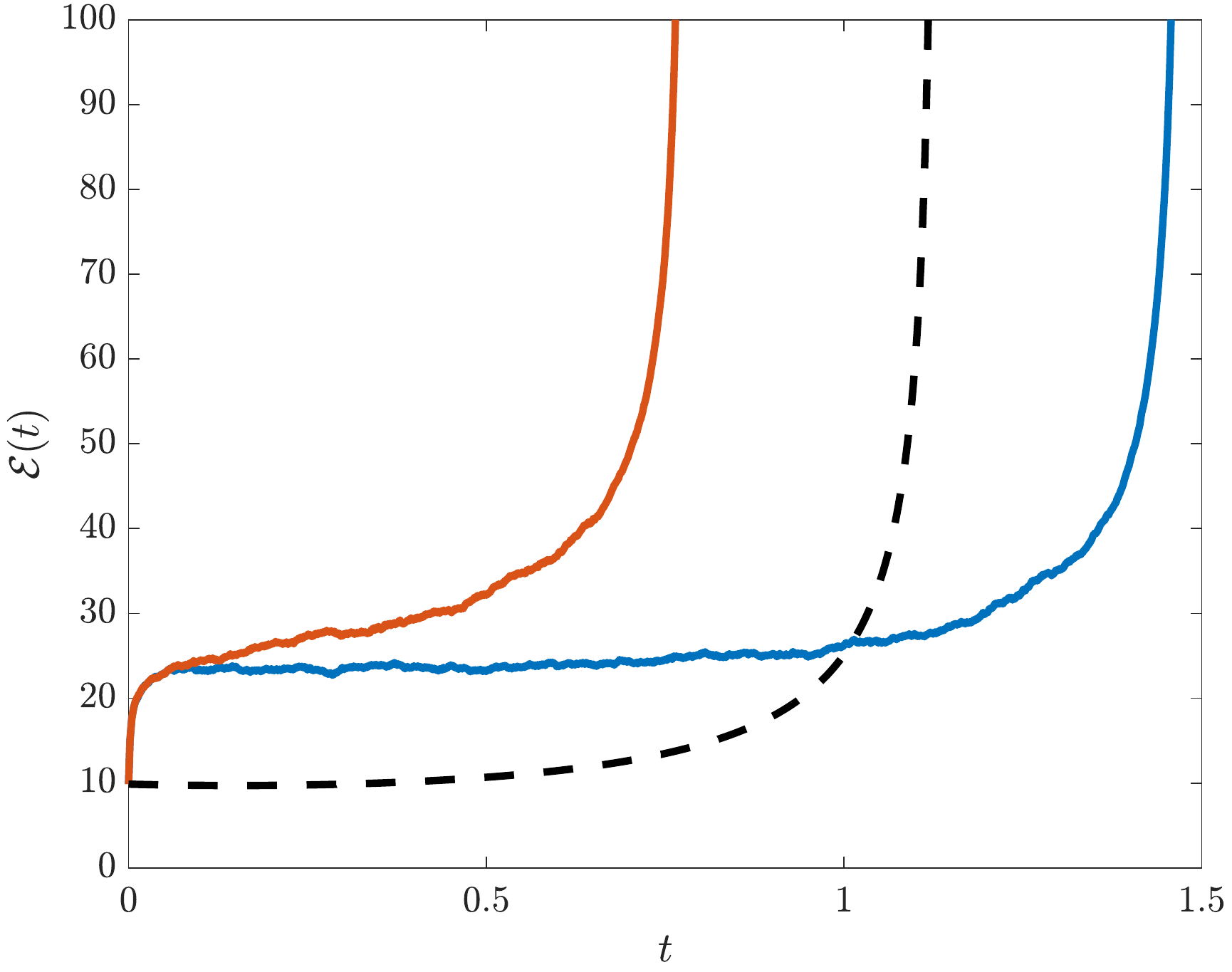}}}
\caption{A single stochastic realization of the solution of system
  \eqref{SBE2} with {$\alpha=0.4$} and $\rho = 10^{-2}$ in (a) the
  physical space $u(t,x)$ and (b) the Fourier space
  $|\widehat{u}_{k}(t)|$ at the indicated time levels with the
  corresponding evolution of the enstrophy $\mathcal{E}(t)$ in panel
  (c). The symbols in panel (c) correspond to the time instances at
  which the solution is shown in panels (a) and (b). (d) Time
  evolution of the enstrophy in the stochastic realizations of the
  solution of system \eqref{SBE2} with {$\alpha=0.4$} and $\rho =
  5\cdot 10^{-2}$ with blow-up occurring at the shortest and latest
  recorded times $T^*$ (solid lines), and in the solution of the
  deterministic problem \eqref{FBE} with the same value of $\alpha$
  (dashed line).}
\label{fig:sup1}
\end{figure}

As regards the supercritical regime, the main question we want to
address is the effect of the stochastic excitation on the blow-up time
$T^*$ and to this end we have solved system \eqref{SBE2} with $\alpha
= 0.4$ and with different noise amplitudes $\rho = 10^{-6}, 10^{-4},
10^{-2}, 2\cdot 10^{-2}, 5\cdot 10^{-2}$, in each case using $M =
2\cdot 10^4$ Monte Carlo samples. The time evolution of a
representative stochastic realization of the solution with an
intermediate noise level $\rho = 10^{-2}$ is shown in Figures
\ref{fig:sup1}a and \ref{fig:sup1}b in the physical and Fourier space,
respectively.  Comparing these results with Figures \ref{fig:a0.4}a,b
corresponding to the deterministic case, we observe that the
large-scale features of the solution remain unaffected and the effect
of noise is to produce small-scale oscillations. As is evident from
Figure \ref{fig:sup1}b, the Fourier spectrum of the solution no longer
decays exponentially fast for large wavenumbers. The corresponding
time evolution of the enstrophy shown in Figures \ref{fig:sup1}c
indicates that its growth is now non-monotonic.  Different stochastic
realizations with a given value of $\rho$ exhibit singularity
formation at times $T^*$ which can be both shorter and longer than the
blow-up time in the deterministic case. However, in none of these
cases did we find any evidence for noise to regularize solutions, so
that blow-up would not occur in finite time at all. To illustrate this
point, in Figures \ref{fig:sup1}d we show the enstrophy evolutions in
the stochastic realizations corresponding to the singularity occurring
at the earliest and latest recorded times $T^*$ when the noise
amplitude is $\rho = 5\cdot 10^{-2}$. It is clear from this figure
that noise may nonetheless significantly accelerate or delay
singularity formation.

\begin{figure}[h]
\centering
\mbox{
\subfigure[]{\includegraphics[width=0.48\textwidth]{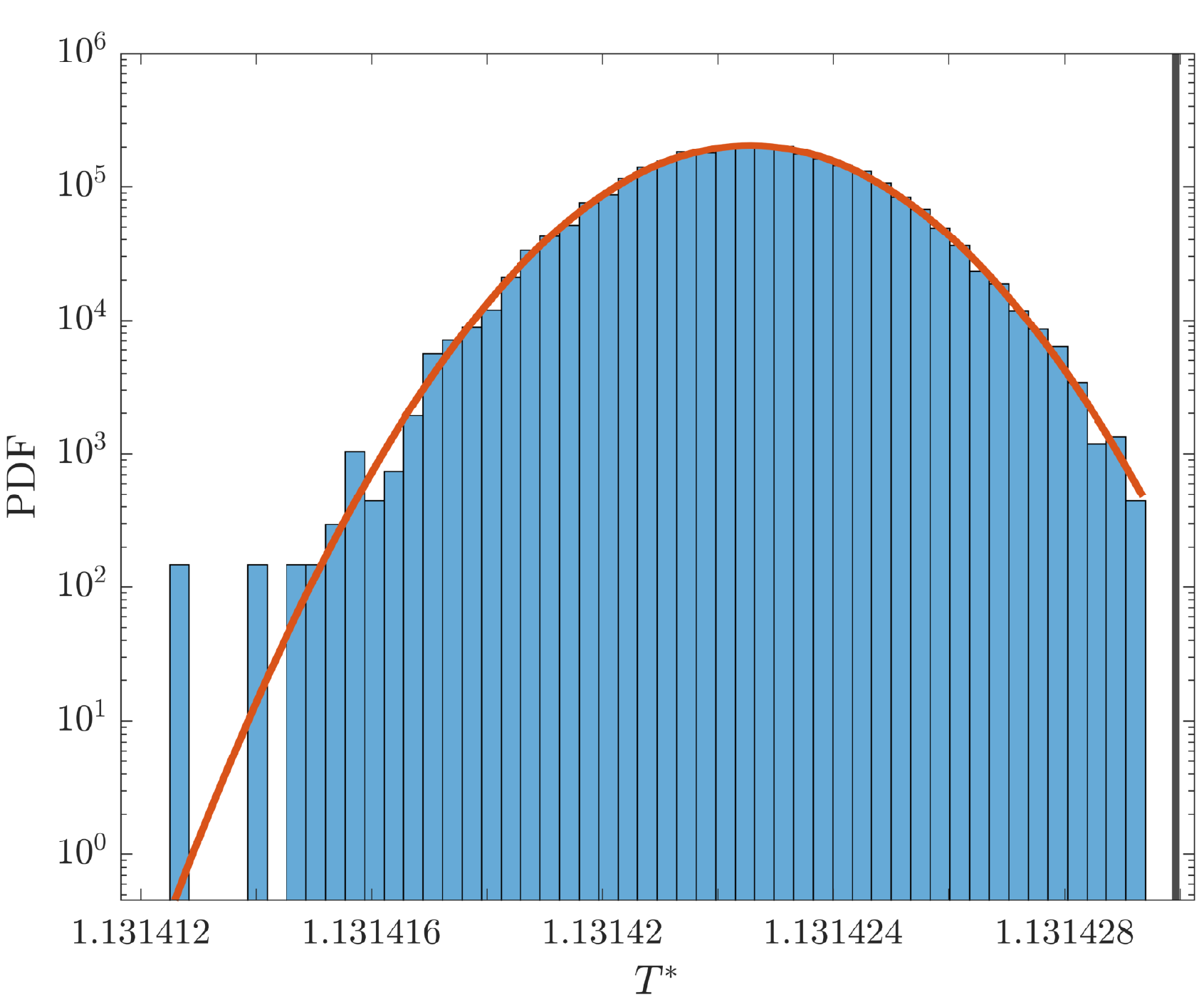}}
\subfigure[]{\includegraphics[width=0.48\textwidth]{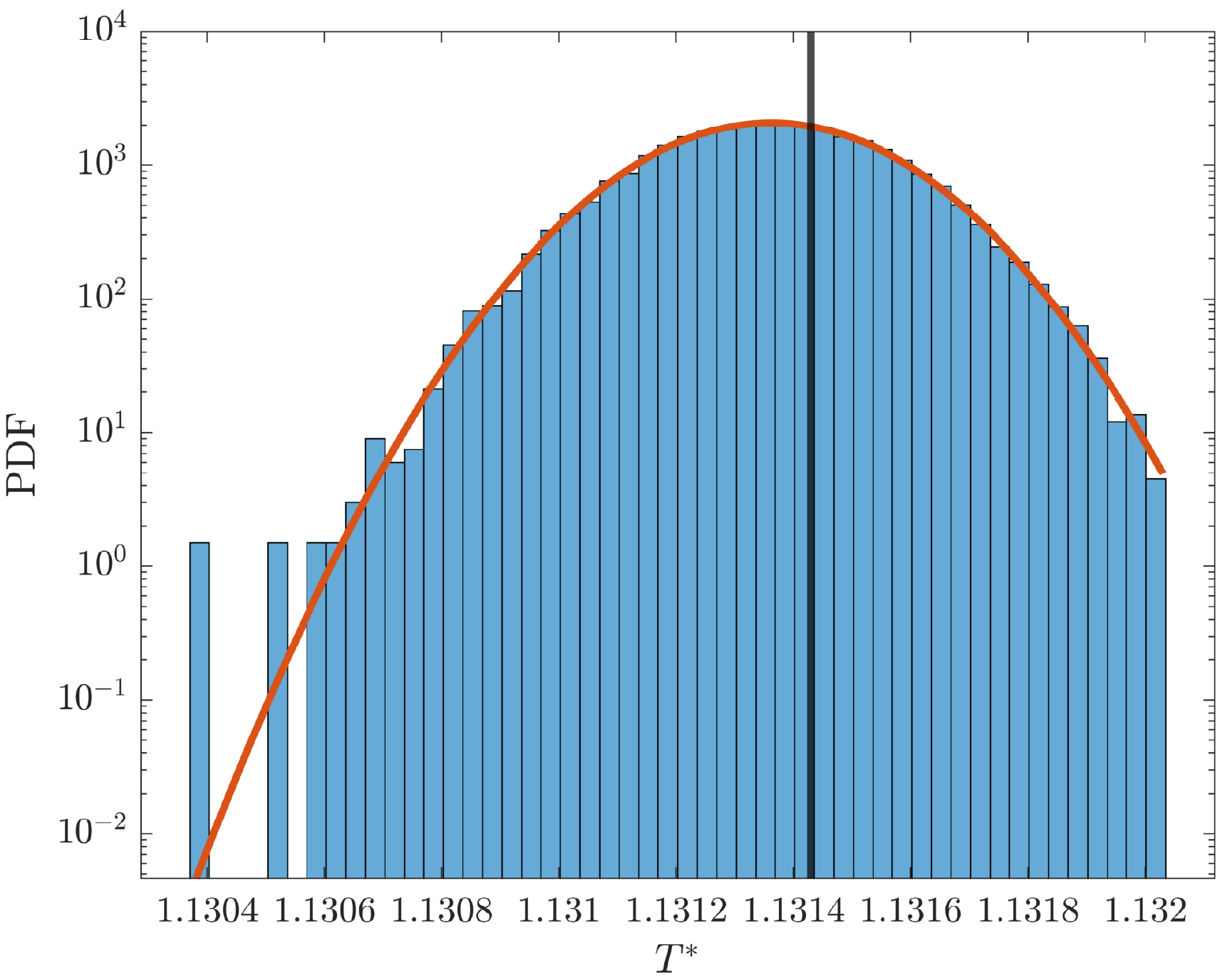}}}
\mbox{
\subfigure[]{\includegraphics[width=0.48\textwidth]{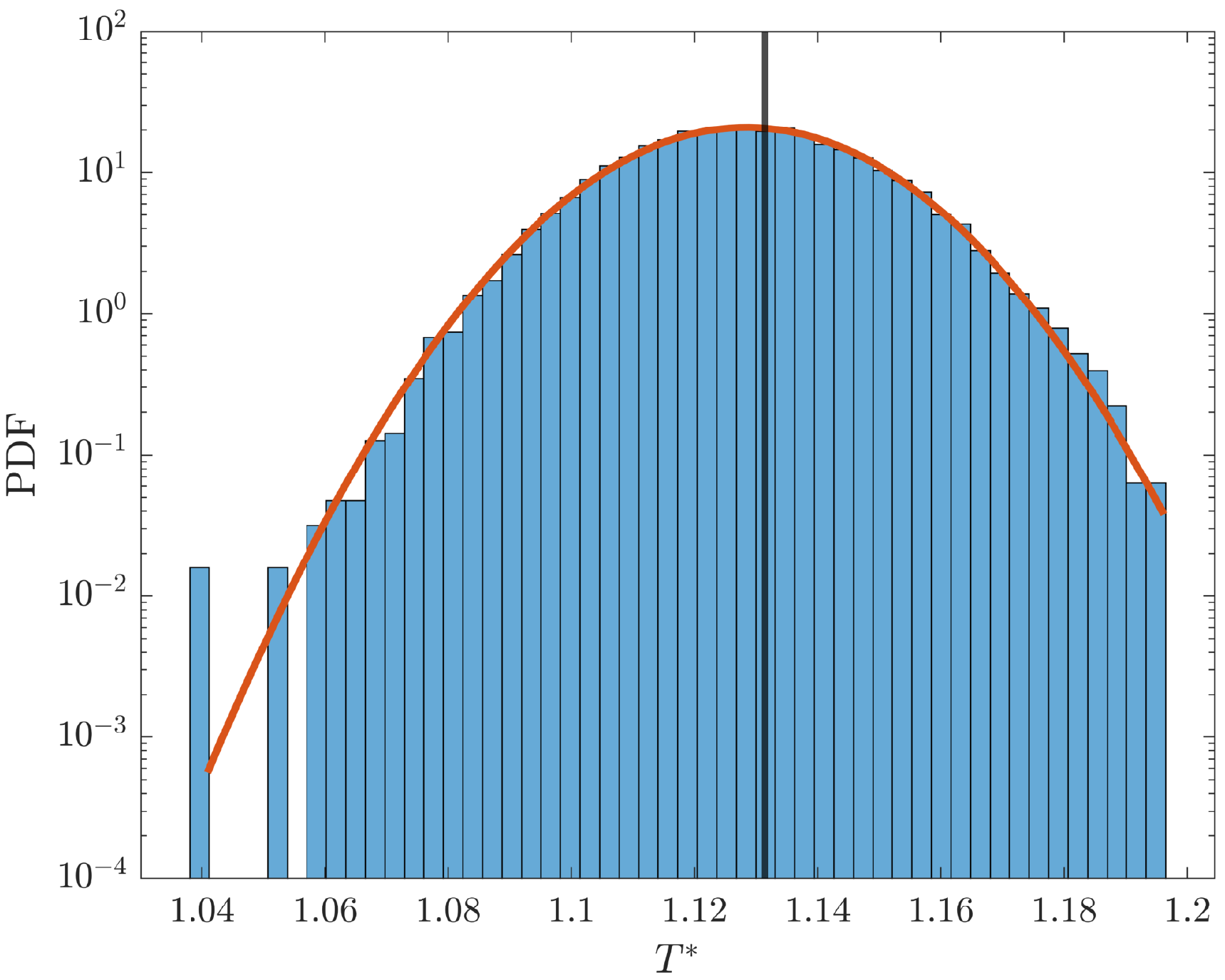}}
\subfigure[]{\includegraphics[width=0.48\textwidth]{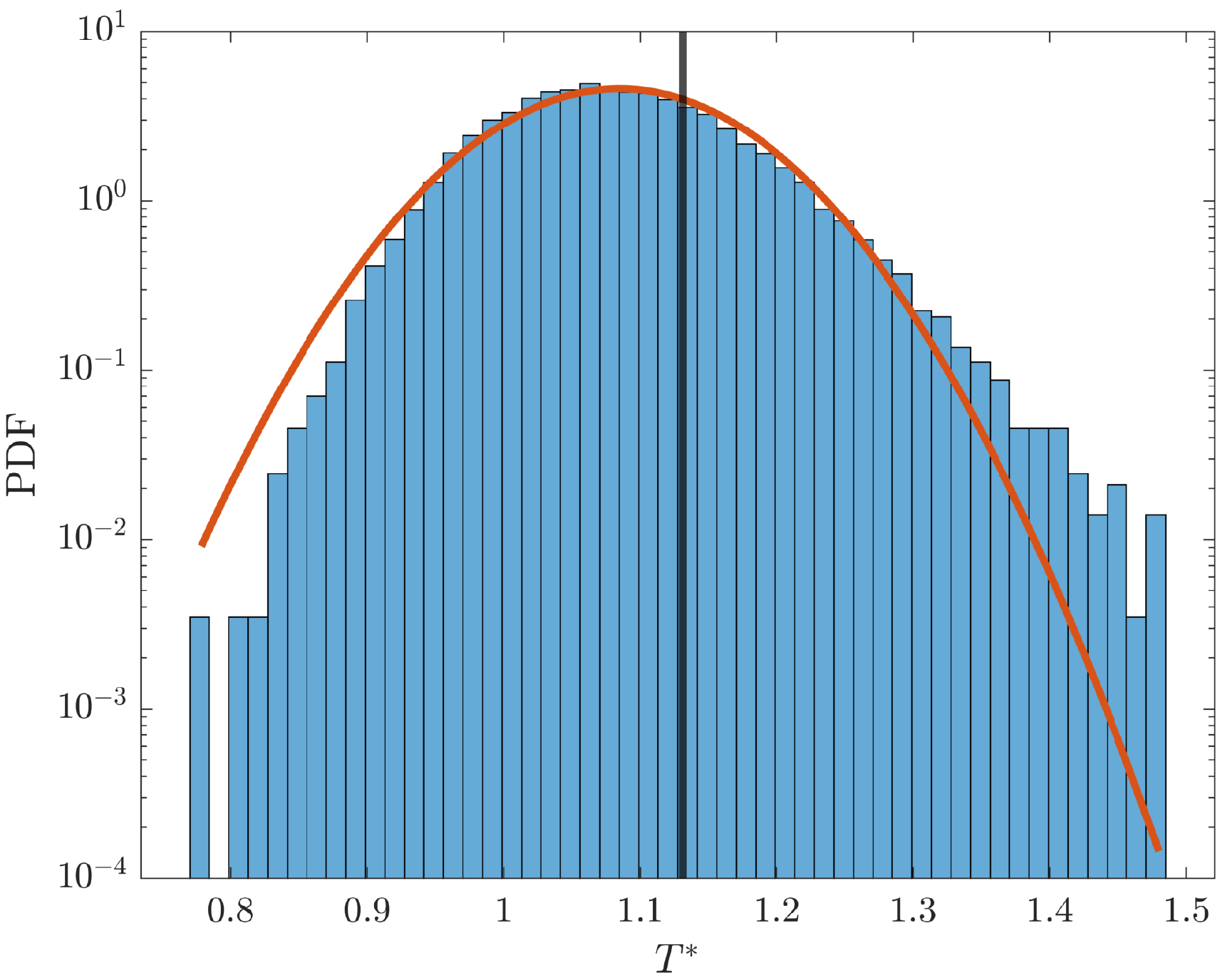}}}
\caption{PDFs of the blow-up time $T^*$ in solutions of the stochastic
  problem \eqref{SBE2} with different noise amplitudes (a)
  $\rho=10^{-6}$ , (b) $\rho=10^{-4}$, (c) $\rho=10^{-2}$ and
  (d) $\rho=5\cdot10^{-2}$. Red curves represent the Gaussian
  distributions with the same means and standard deviations, whereas
  black vertical lines denote the blow-up times in the deterministic
  case.}
\label{fig:supPDF}
\end{figure}

The blow-up time $T^*$, estimated here using Algorithm \ref{alg:Ts},
is a random variable and we now go on to analyze its statistical
properties. Probability density functions (PDFs) of $T^*$ are shown in
Figures \ref{fig:supPDF}a--d for four different values of the noise
amplitude $\rho$. We see that while for small noise magnitudes the
PDFs of the blow-up time have distributions close to the Gaussian
distribution, for increasing $\rho$ they develop ``heavy'' tails and
become skewed towards longer times. Interestingly, this increasing
non-Gaussianity of the distributions is accompanied by the mean value of
$T^*$ shifting towards shorter times while its standard deviation
increases. This means that as the noise amplitude $\rho$ increases,
blow-up on average occurs earlier than in the deterministic case,
however, realizations in which blow-up is significantly delayed also
become more likely.

\begin{figure}[h]
\centering
\mbox{
\subfigure[]{\includegraphics[width=0.48\textwidth]{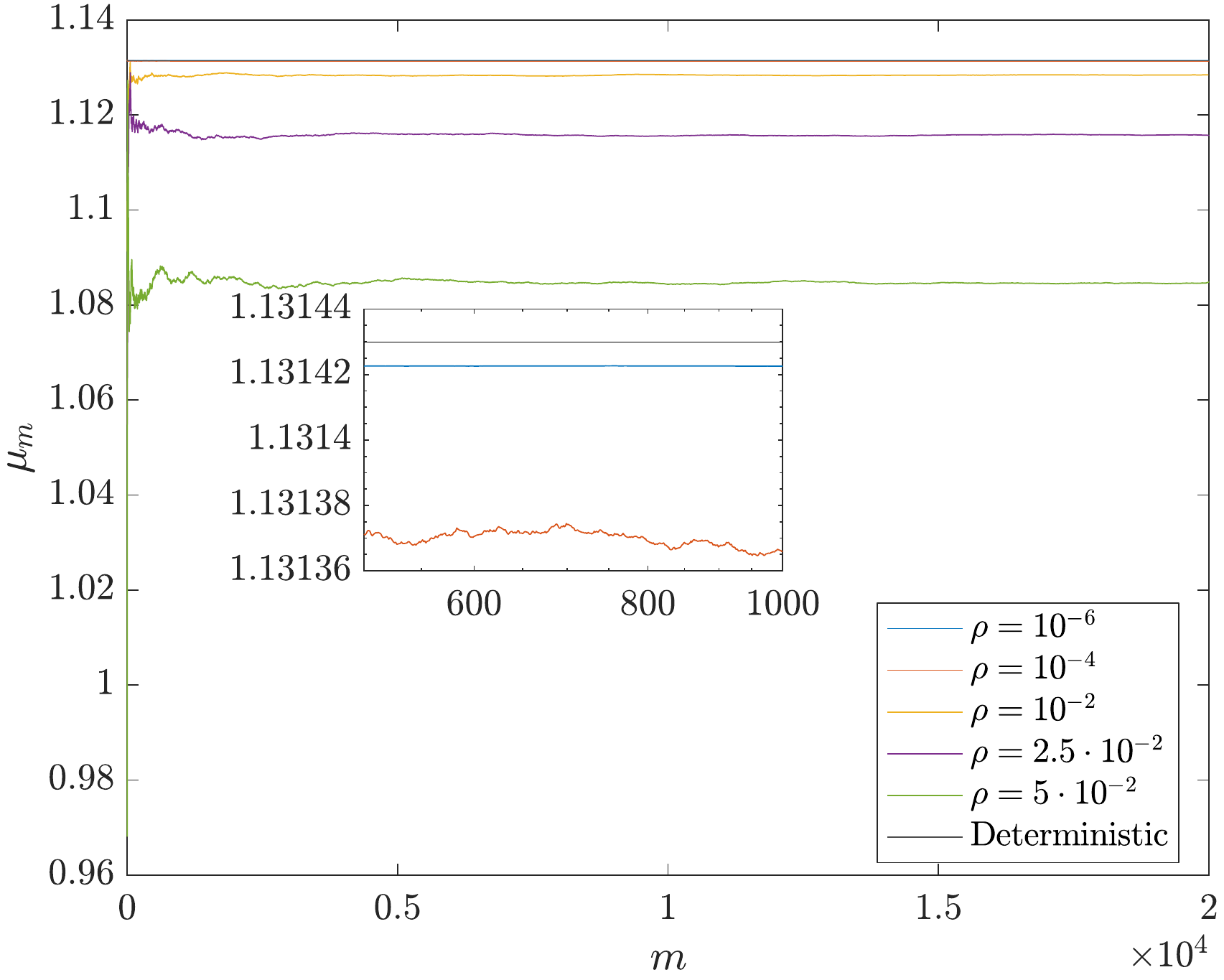}}
\subfigure[]{\includegraphics[width=0.48\textwidth]{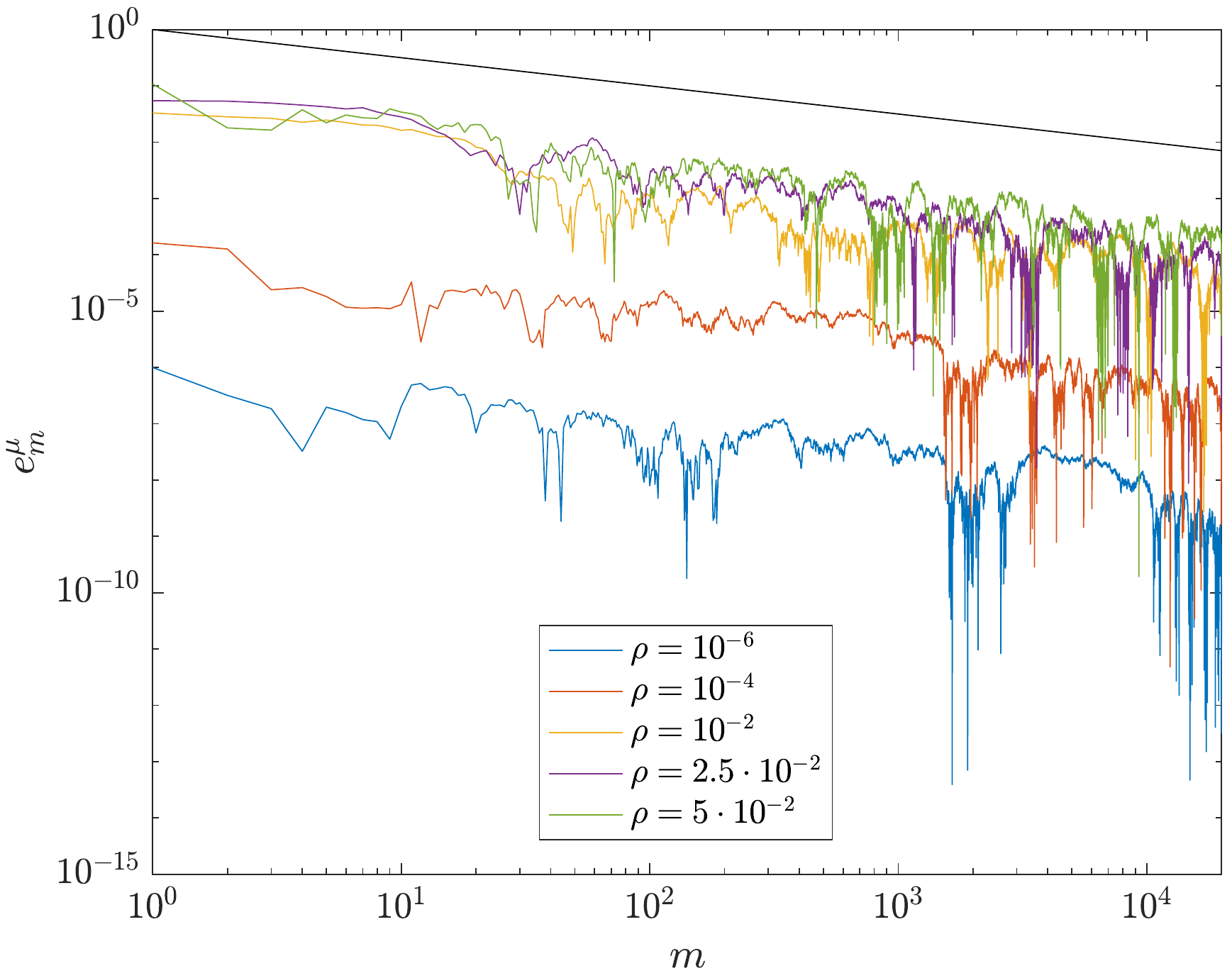}}}
\mbox{
\subfigure[]{\includegraphics[width=0.48\textwidth]{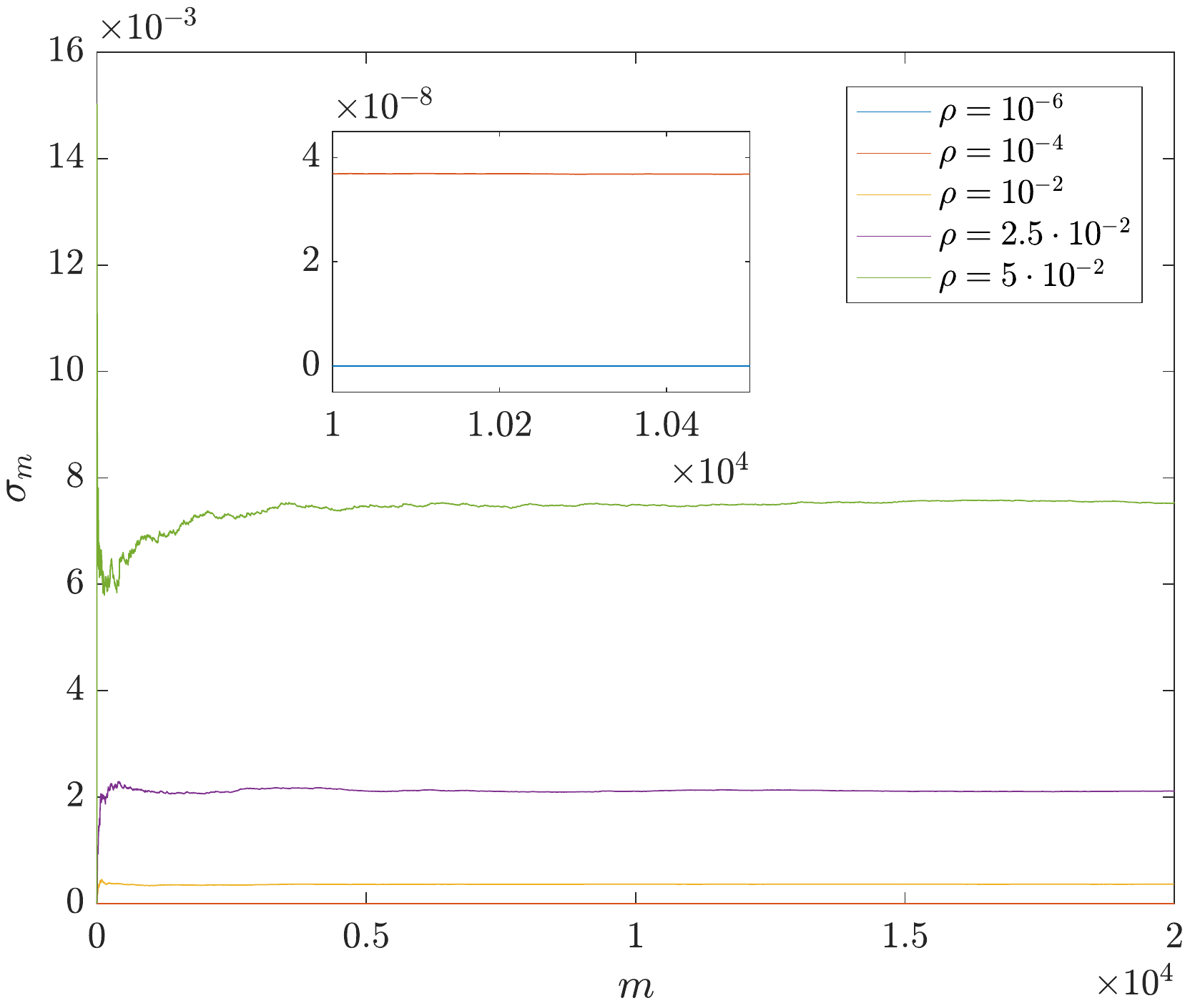}}
\subfigure[]{\includegraphics[width=0.48\textwidth]{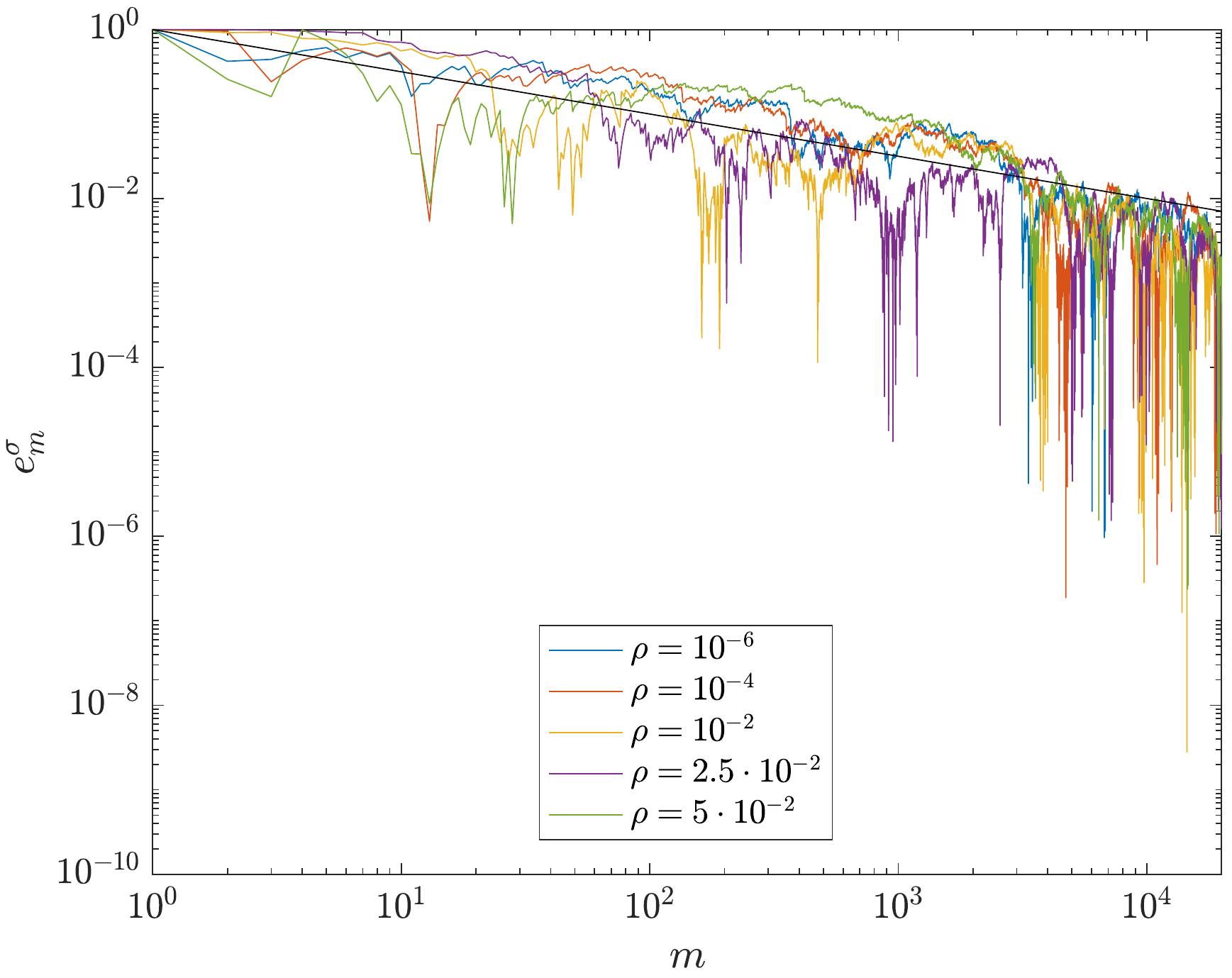}}}
\caption{Approximations (a) $\mu_m$ of the mean and (c) $\sigma_m$ of
  the standard deviation together with their relative errors (b)
  $e^{\mu}_m$ and (d) $e^{\sigma}_m$ as functions of the number of
  samples $m$ used in the computation.  The straight black line in
  panels (b) and (d) represents the expression $m^{-1/2}$.}
\label{fig:m1}
\end{figure}

\begin{figure}[h]
\centering
\mbox{
\subfigure[]{\includegraphics[width=0.48\textwidth]{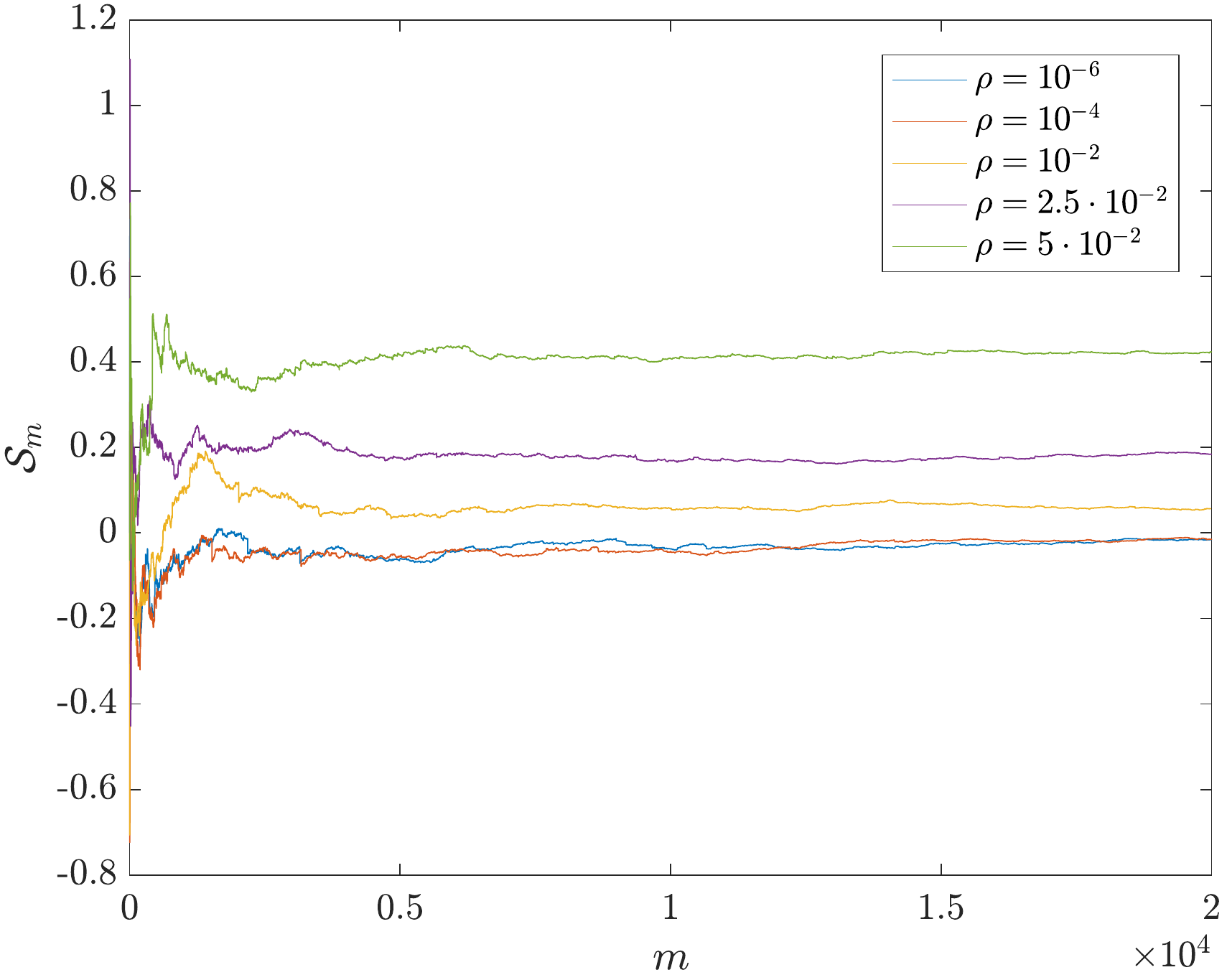}}
\subfigure[]{\includegraphics[width=0.48\textwidth]{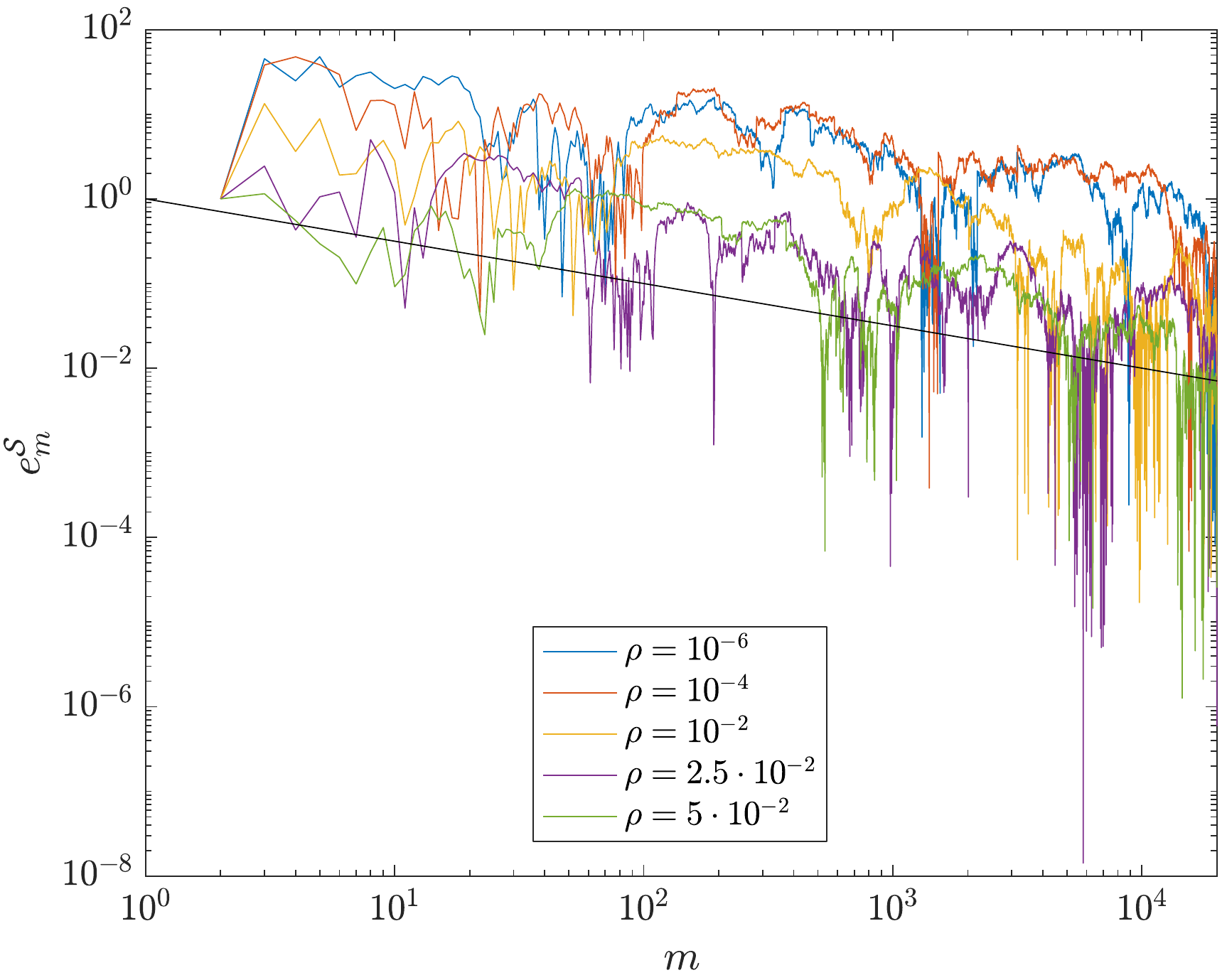}}}
\mbox{
\subfigure[]{\includegraphics[width=0.48\textwidth]{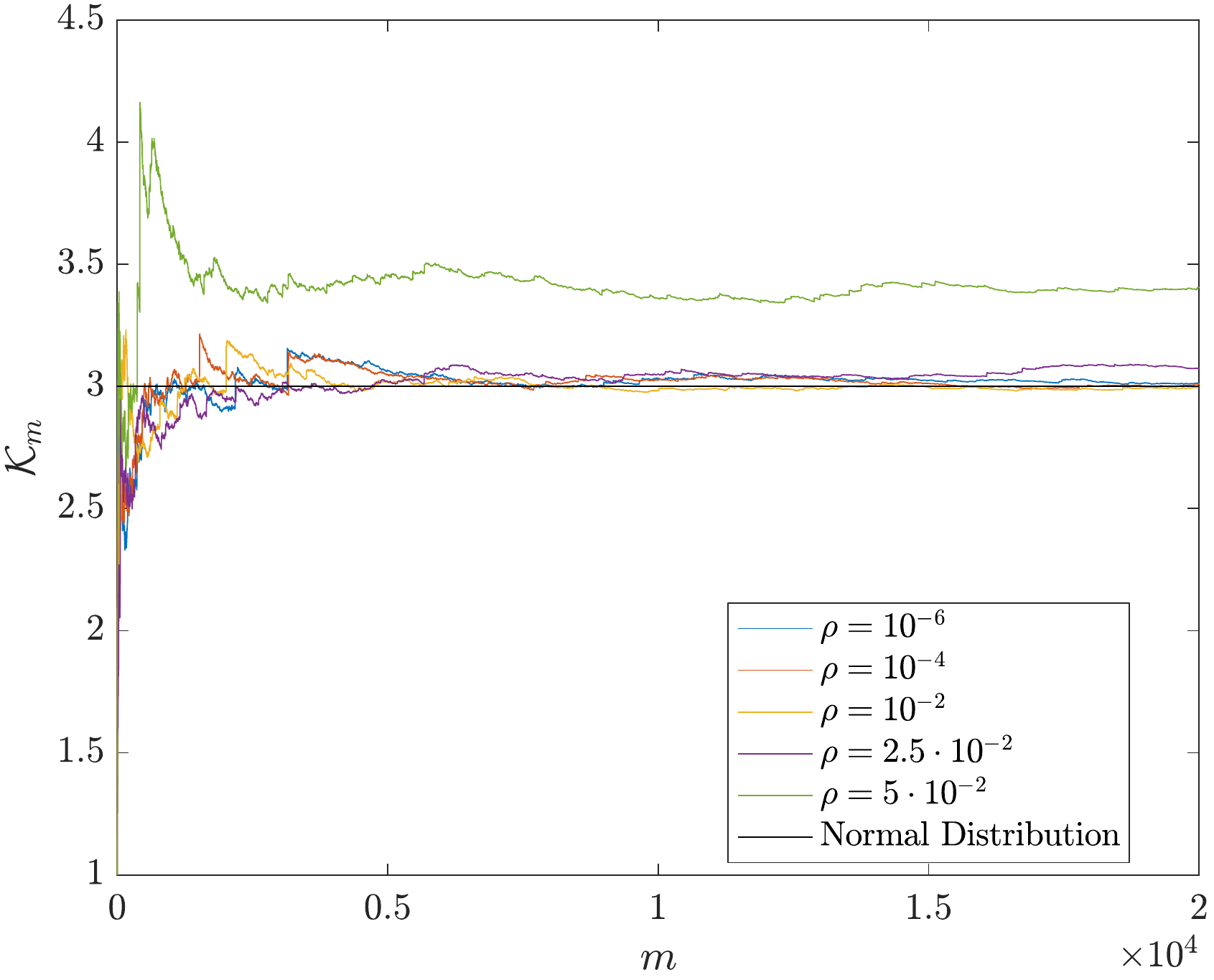}}
\subfigure[]{\includegraphics[width=0.48\textwidth]{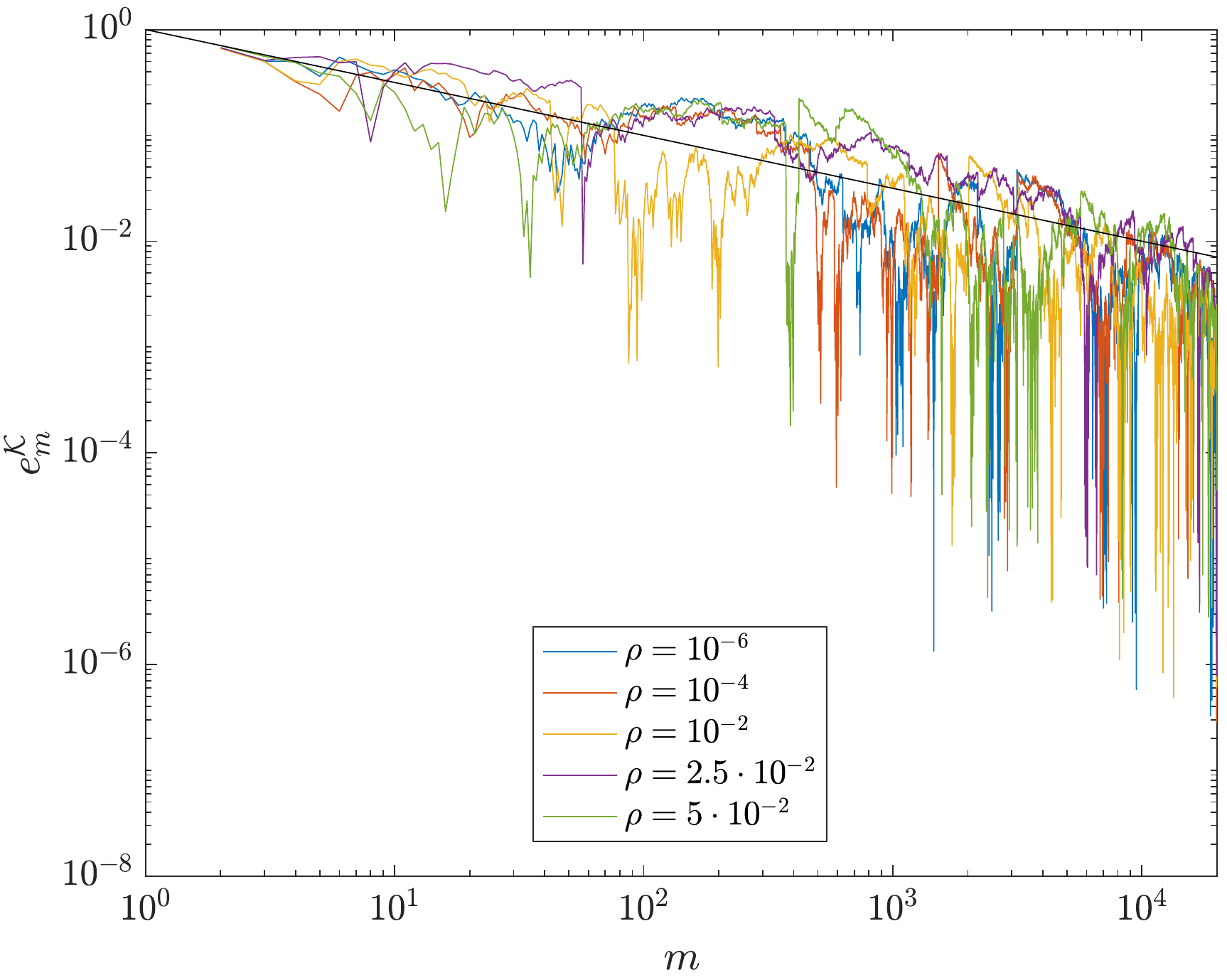}}}
\caption{Approximations (a) $\S_m$ of the skewness and (c) $\K_m$ of
  the kurtosis together with their relative errors (b) $e^{\S}_m$ and
  (d) $e^{\K}_m$ as functions of the number of samples $m$ used in the
  computation.  The straight black line in panels (b) and (d)
  represents the expression $m^{-1/2}$.}
\label{fig:m2}
\end{figure}

In order to quantify the observations made above, we compute the first
four statistical moments, namely the mean $\mu := \EE[T^*]$, standard
deviation $\sigma := \EE[(T^* - \mu)^2]$, skewness $\S :=
\EE\left[\left((T^*-\mu)/\sigma\right)^3\right]$ and
kurtosis
$\K:=\EE\left[\left((X-\mu)/\sigma\right)^4\right]$, where
$\EE[\cdot]$ is the expected value, of the blow-up times obtained for
different noise amplitudes $\rho$. These moments are approximated as
follows
\begin{subequations}
\label{eq:moments}
\begin{alignat}{2}
\mu & \approx & \mu_m & := \frac{1}{m} \sum_{j=1}^m T_j^*, \label{eq:mu} \\
\sigma & \approx & \sigma_m & := \frac{1}{m} \sum_{j=1}^m (T_j^* - \mu_m)^2, \label{eq:sig}  \\
\S & \approx & \S_m & := \frac{1}{m} \sum_{j=1}^m \left(\frac{T_j^* - \mu_m}{\sigma_m}\right)^3, \label{eq:S}  \\
\K & \approx & \K_m & := \frac{1}{m} \sum_{j=1}^m \left(\frac{T_j^* - \mu_m}{\sigma_m}\right)^4, \label{eq:K} 
\end{alignat}
\end{subequations}
where $m = 1,\dots,M$ is the number of samples used to construct the
approximations and $\{ T_j^* \}_{j=1}^M$ are the blow-up times in the
realizations from the Monte-Carlo ensemble. It is known that moments
of increasing order are harder to approximate accurately via relations
\eqref{eq:moments} for random variables characterized by PDFs with
heavy tails, such as the PDF shown in Figure \ref{fig:supPDF}d.  In
order to estimate how many samples are required to accurately estimate
the moments for different values of the noise amplitude $\rho$, in
Figures \ref{fig:m1}a,c and \ref{fig:m2}a,c we show the estimates
$\mu_m$, $\sigma_m$, $\S_m$ and $\K_m$ for $m =
1,\dots,M$. These results indicate how the estimated moments depend on
the number of samples $m$ used in the calculation. They are
complemented by plots of the relative errors defined as
\begin{equation}
e^X_m := \left|\frac{X_m - X_M}{X_M}\right|, \qquad m=1,\dots,M, \quad X = \mu, \sigma, \S, \K,
\label{eq:merr}
\end{equation}
where the estimates $X_M$ computed based on all available $M$ samples
are used as the ``true'' values, shown in Figures \ref{fig:m1}b,d and
\ref{fig:m2}b,d. We see in these last plots that the relative errors
decrease in proportion to $m^{-1/2}$, as expected for a Monte-Carlo
approximation. Overall, in Figures \ref{fig:m1}a--d and
\ref{fig:m2}a--d we observe that approximation errors are larger for
moments of higher order and they also increase with the noise
magnitude $\rho$ as the PDFs become increasingly non-Gaussian,
cf.~Figure \ref{fig:supPDF}. However, we can conclude that the total
number of Monte-Carlo samples $M = 2\cdot10^4$ is sufficient to ensure
required accuracy in all cases. On the other hand, the trends evident
in Figures \ref{fig:m1}a--d and \ref{fig:m2}a--d indicate that a
significantly larger number of Monte-Carlo samples would be required
to obtain converged statistics for a higher value of $\rho$, which
was not attempted due to a prohibitive computational cost.

\begin{figure}[h]
\centering
\mbox{
\subfigure[]{\includegraphics[width=0.48\textwidth]{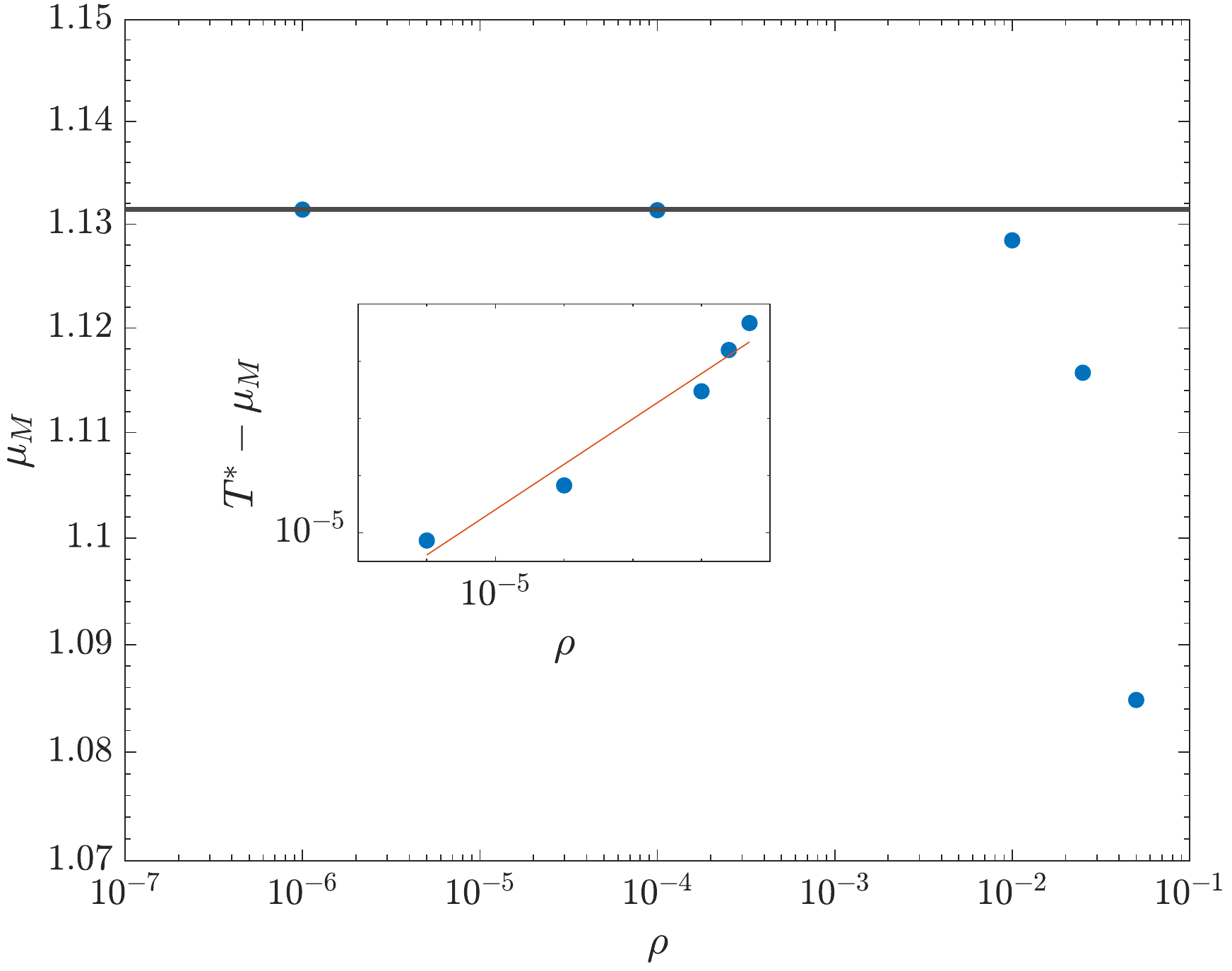}}
\subfigure[]{\includegraphics[width=0.495\textwidth]{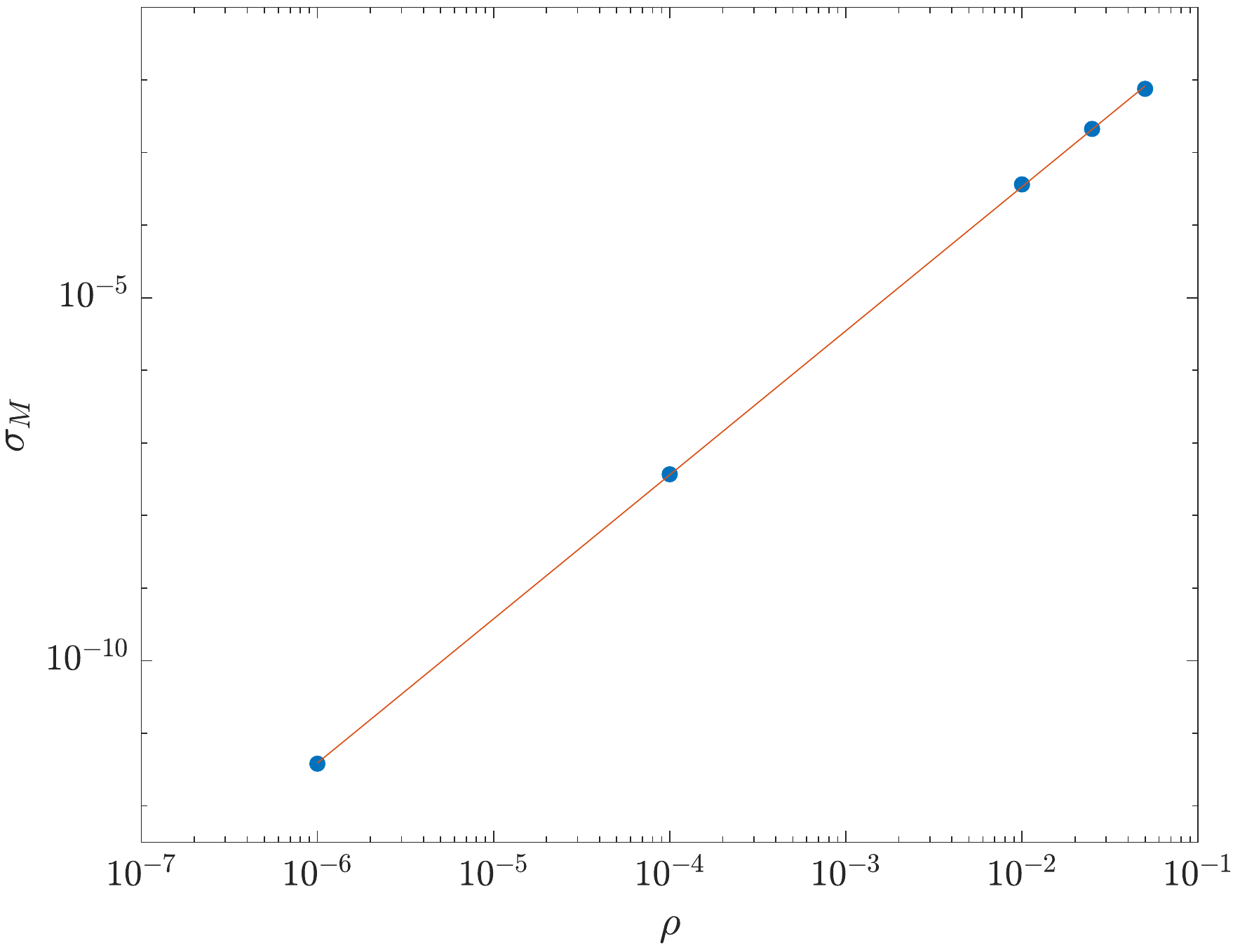}}}
\mbox{
\subfigure[]{\includegraphics[width=0.48\textwidth]{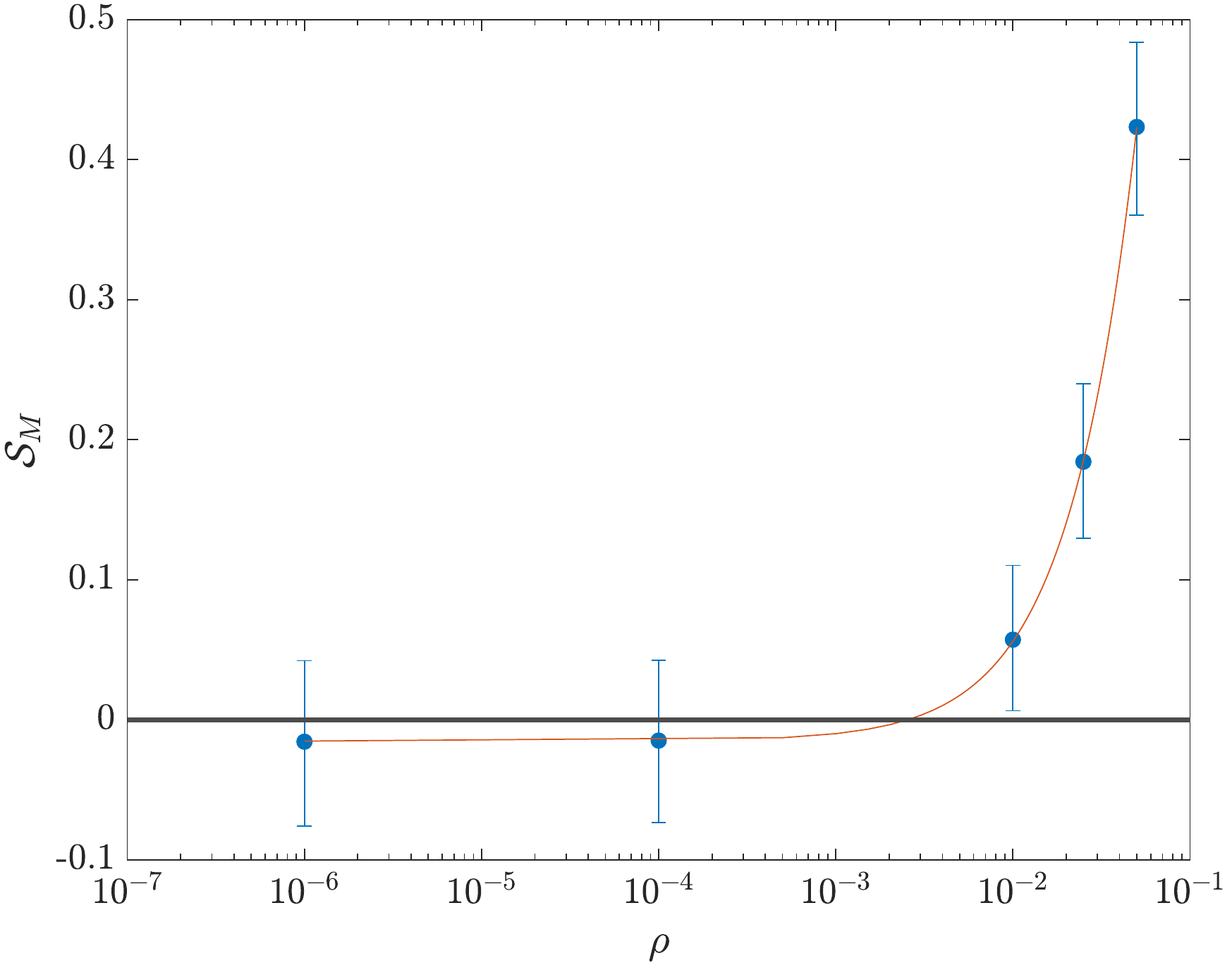}}
\subfigure[]{\includegraphics[width=0.48\textwidth]{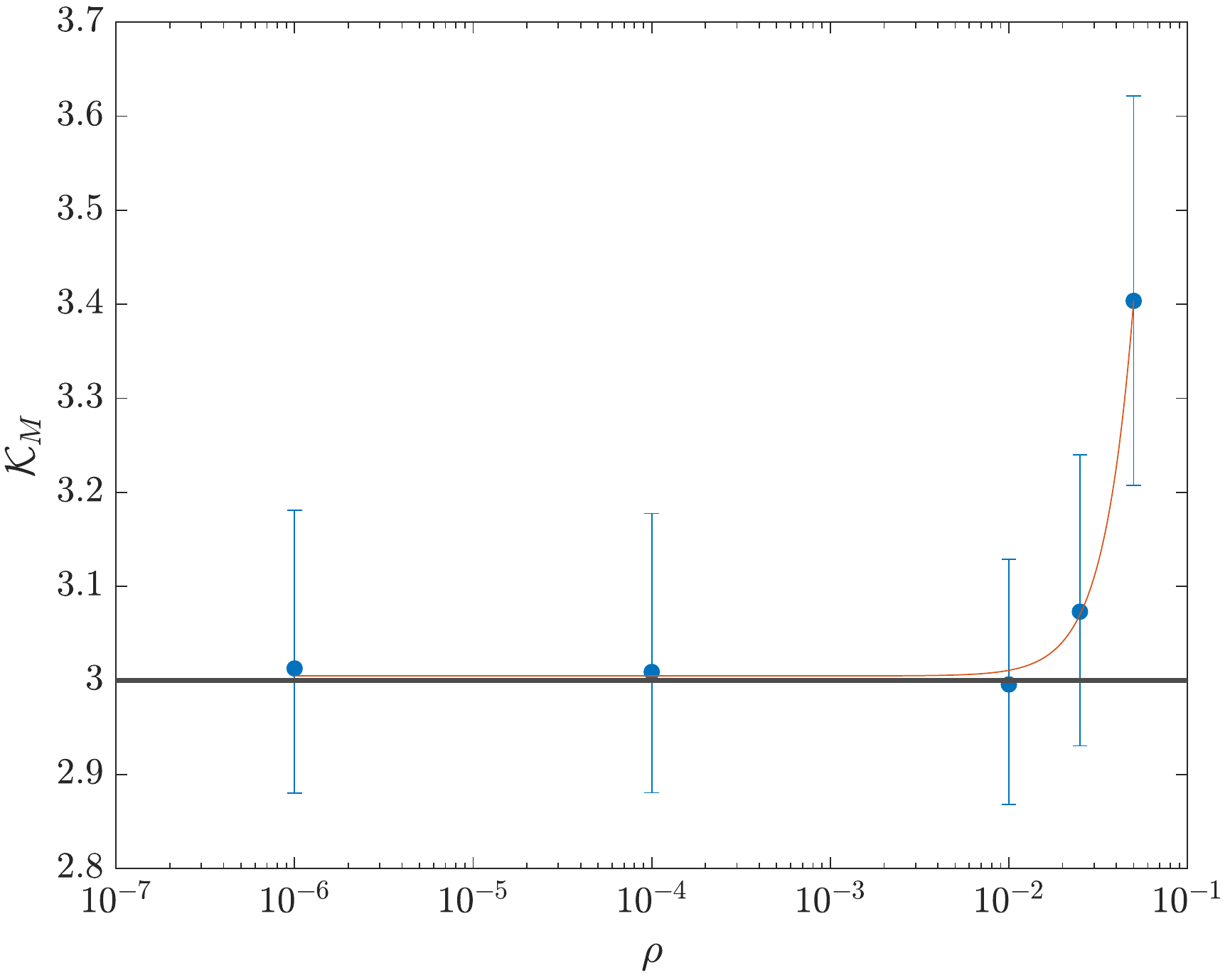}}}
\caption{Estimated (a) mean blow time $T^*$, (b) its standard
  deviation, (d) skewness and (e) kurtosis as functions of the noise
  amplitude $\rho$. The horizontal line in (a) corresponds to the
  blow-up time in the deterministic case, whereas the horizontal lines
  in (c) and (d) represent the values characterizing the Gaussian
  distribution, The red curves correspond to the fits obtained
    using the power-law relation \eqref{eq:fit} with the parameters
    reported in Table \ref{tab:mf}.}
\label{fig:msig}
\end{figure}

Finally, the mean blow-up time, its standard deviation, skewness and
kurtosis, which were estimated using relations
\eqref{eq:mu}--\eqref{eq:K}, are shown as functions of $\rho$ in
Figures \ref{fig:msig}a--d. The plots also involve error bars computed
using the bootstrapping method \cite{bt96,ad97} assuming 95\%
confidence which reflect the approximation errors due to finite
numbers of Monte-Carlo samples, cf.~Figures \ref{fig:m1} and
\ref{fig:m2}.  Figures \ref{fig:msig}a and \ref{fig:msig}b show that
over the considered range of $\rho$ the mean blow-up time decreases
with noise amplitude whereas its standard deviation grows.  Figures
\ref{fig:msig}c and \ref{fig:msig}d show that skewness and
kurtosis also grow as $\rho$ becomes large providing evidence for an
increasing non-Gaussianity of the distribution of blow-up times. In
order to quantify these trends, Figures \ref{fig:msig}a--d include
fits constructed using the power-law relation
\begin{equation}
f(\rho) := a\,\rho^b +c, \qquad a,b,c \in \RR,
\label{eq:fit}
\end{equation}
with the obtained parameters reported in Table \ref{tab:mf} (we set
$c=0$ in the fits to the mean blow-up time and its variance). One can
see that while the dependence of the variance, skewness and kurtosis
on $\rho$ is well-represented by relation \eqref{eq:fit}, this is not
the case for the mean blow-up time. However, for each statistical
moment the optimal power-law relation involves a different exponent
$b$.

\begin{table}[t]
\centering
\begin{tabular}{|c|c|c|c|}
\hline
Statistical Moment&a&b&c \\ \hline
$\mu_M$ & 0.233 & 0.792&0\\ \hline
$\rho_M$ & 3.136&1.984&0\\\hline
$\mathcal{S}_M$ & 12.888&1.128&-0.015\\ \hline
$\mathcal{K}_M$ & 1008.425&2.615&3\\\hline
\end{tabular}
\caption{Parameters of the fits to the data shown in Figures 
\ref{fig:msig}a--d using the power-law  relation \eqref{eq:fit}.}
\label{tab:mf}
\end{table}

\FloatBarrier

\subsection{Subcritical Regime}
\label{sec:results_sub}

We now go on to analyze the results obtained for the subcritical case
with $\alpha = 0.6$ and with different noise amplitudes $\rho =
10^{-6}, 10^{-4}, 10^{-2}, 5\cdot 10^{-2}, 10^{-1}$. In each case we
use $M = 1.2\cdot 10^4$ Monte-Carlo samples which was found to be a
sufficient number based on an analysis similar to that reported in
Figures \ref{fig:m1}--\ref{fig:m2}.  Since no evidence was found for
noise-induced blow-up in the subcritical case, the key question is
about the effect of noise of on the maximum attained enstrophy $\Em :=
\max_{t \ge 0} \E(t)$ and on the time when this maximum occurs $\Tm :=
\argmax_{t \ge 0} \E(t)$.

\begin{figure}[h]
\centering
\mbox{
\subfigure[]{\includegraphics[width=0.48\textwidth]{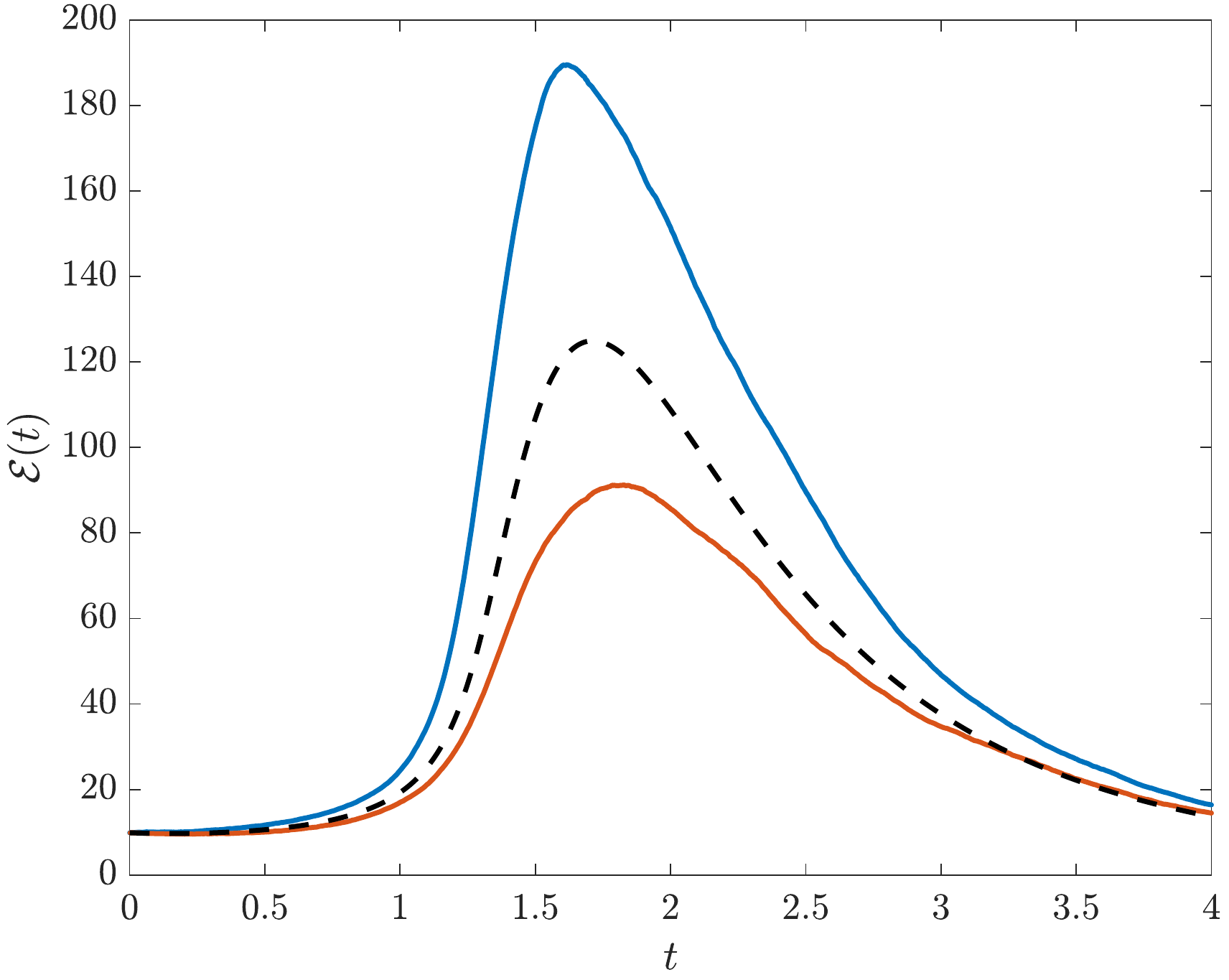}}
\subfigure[]{\includegraphics[width=0.48\textwidth]{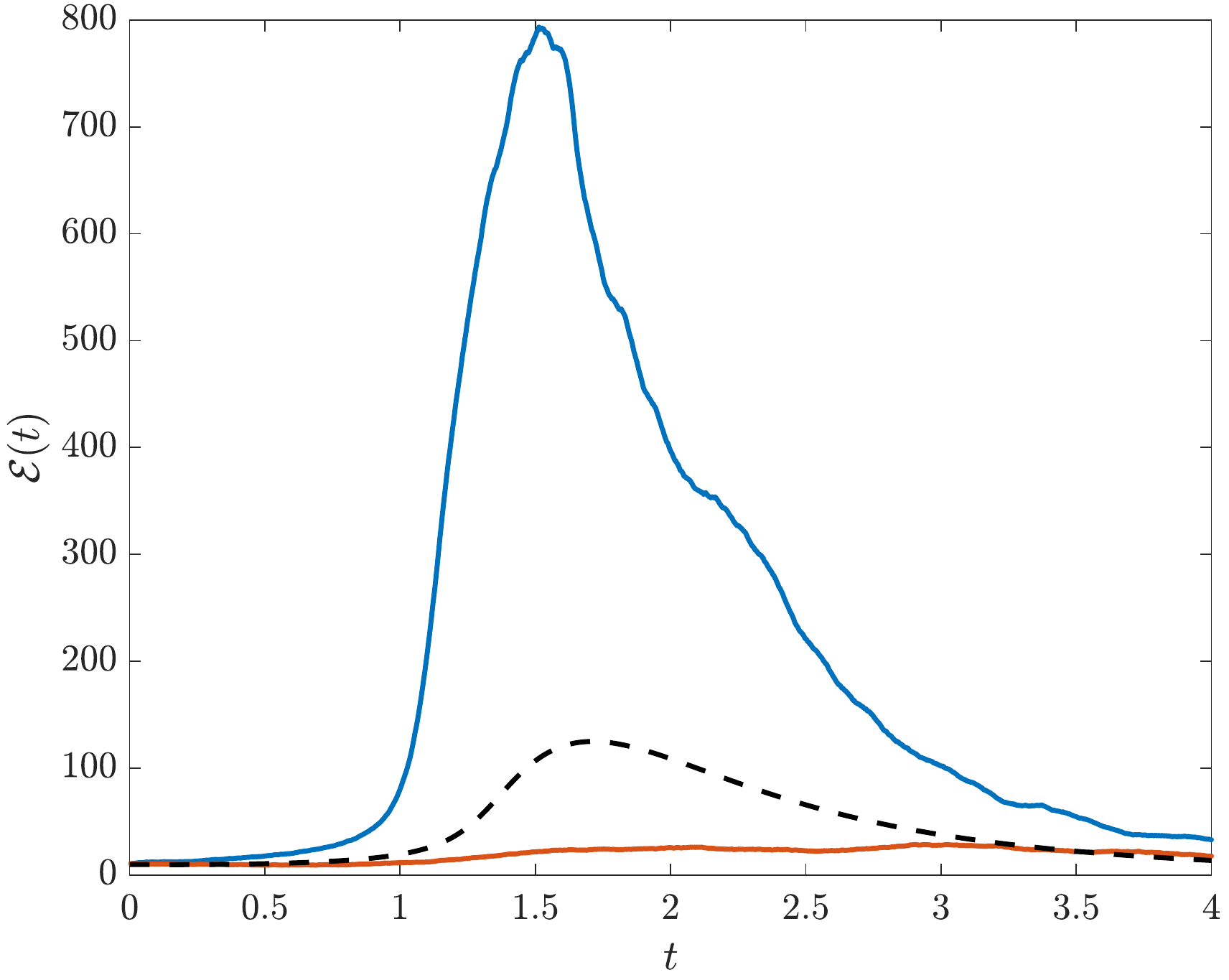}}}
\subfigure[]{\includegraphics[width=0.48\textwidth]{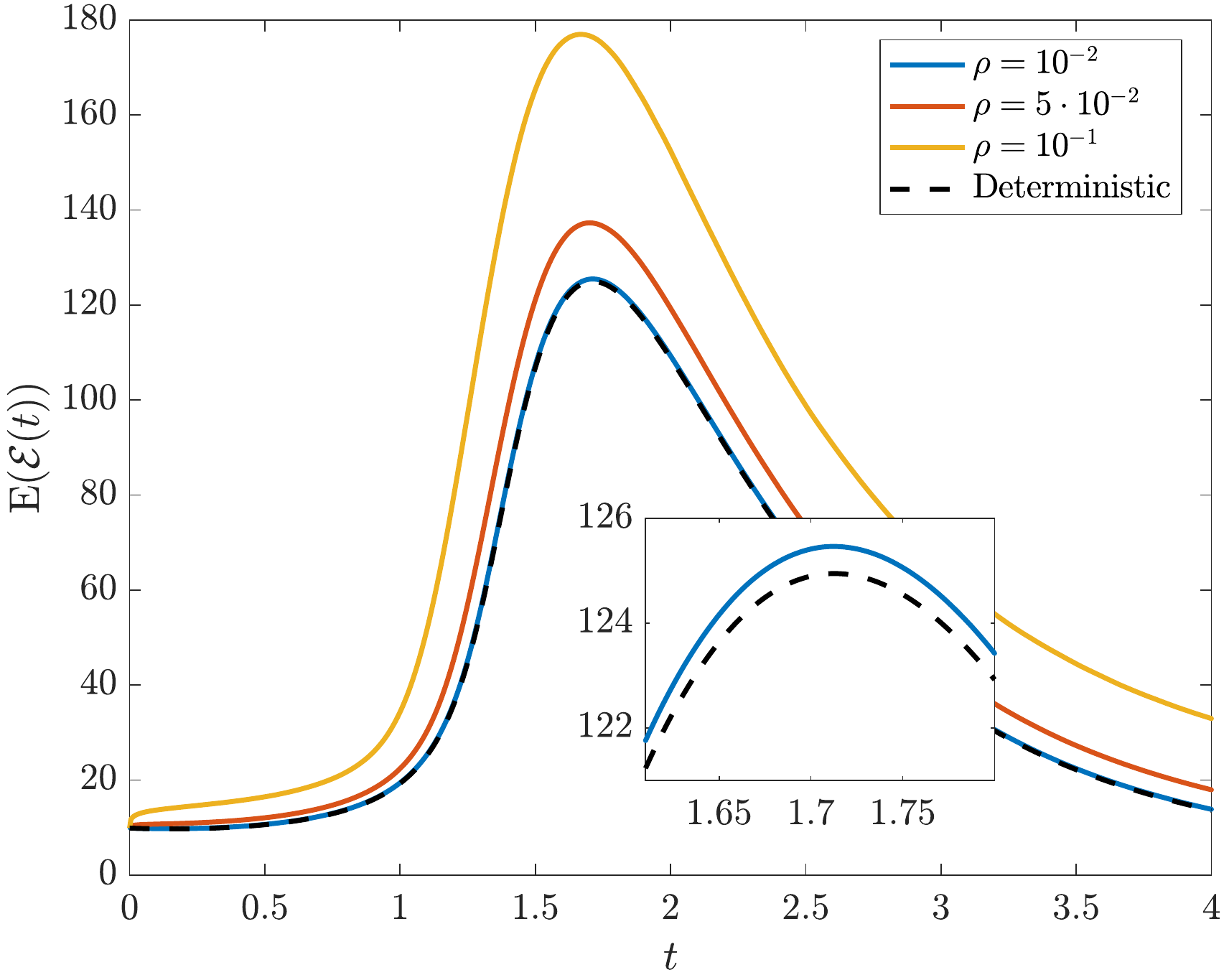}}
\caption{Solid lines in (a) and (b) represent the time evolutions of
  the enstrophy $\E(t)$ in the stochastic realizations characterized
  by the largest and smallest values of $\Em$ for (a) $\rho =
  10^{-2}$ and (b) $\rho = 5\cdot 10^{-2}$. The dashed line in (a)
  and (b) corresponds to the deterministic case. Panel (c) shows the
  time evolution of the expected value of the enstrophy
  $\EE[\E(u(t))]$ of the stochastic solution for the three largest
  values of $\rho = 10^{-2}, 5\cdot 10^{-2}, 10^{-1}$.}
\label{fig:subEt}
\end{figure}

We begin by comparing the time evolution of the enstrophy in the
stochastic realizations with the largest and smallest values of $\Em$
obtained for $\rho = 10^{-2}$ and $\rho = 5\cdot 10^{-2}$ in Figures
\ref{fig:subEt}a and \ref{fig:subEt}b, respectively, where we also
show the evolution of the enstrophy in the deterministic case. We
observe that when the maximum enstrophy achieved in a stochastic
realization is higher than in the deterministic case, the maximum
tends to occur at an earlier time, and vice versa. By comparing
Figures \ref{fig:subEt}a and \ref{fig:subEt}b we also note that the
spread between the largest and smallest values of $\Em$ increases with
the noise magnitude $\rho$. In Figure \ref{fig:subEt}c we show the
dependence of the expected value of the enstrophy $\EE[\E(u(t))]$ of
the stochastic solution on time $t$. We see that as $\rho$ increases
then so do the expected values of the enstrophy $\EE[\E(u(t))]$ at any
fixed time $t \in [0,T]$ and their maxima occur at earlier times.

We remark that one might be tempted to also consider the enstrophy of
the expected values of the solution $\E(\EE[u(t)])$. However, these
two quantities (and their estimates of the type \eqref{eq:mu}) are
related via Jensen's inequality \cite{lps14,PocasProtas2018}
\begin{equation}
\EE[\E(u(t))] \ge \E(\EE[u(t)]), \qquad \forall t,
\label{eq:Jensen}
\end{equation}
and therefore we will exclusively focus here on the former quantity.

\begin{figure}[h]
\centering
\mbox{
\subfigure[]{\includegraphics[width=0.51\textwidth]{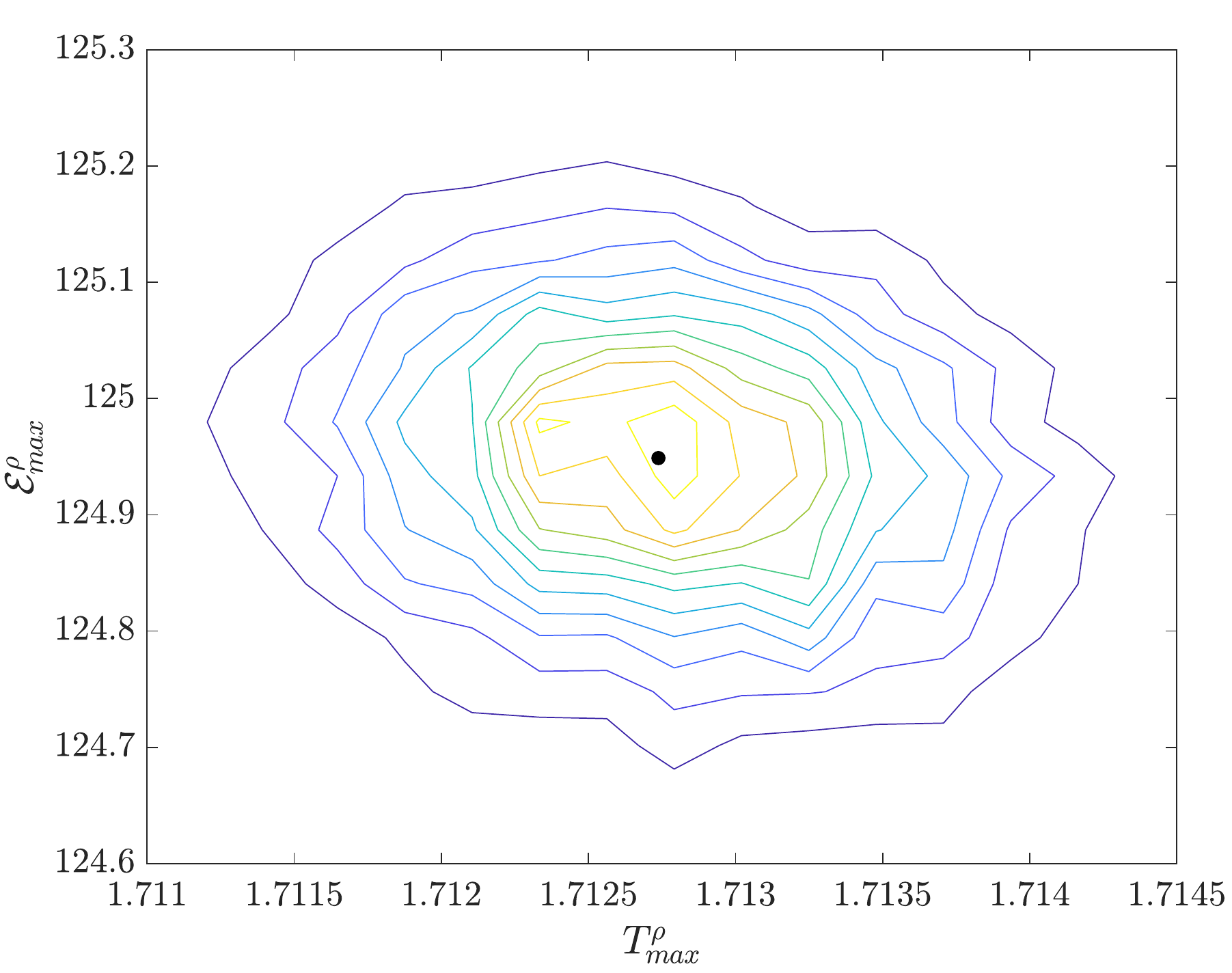}}
\subfigure[]{\includegraphics[width=0.48\textwidth]{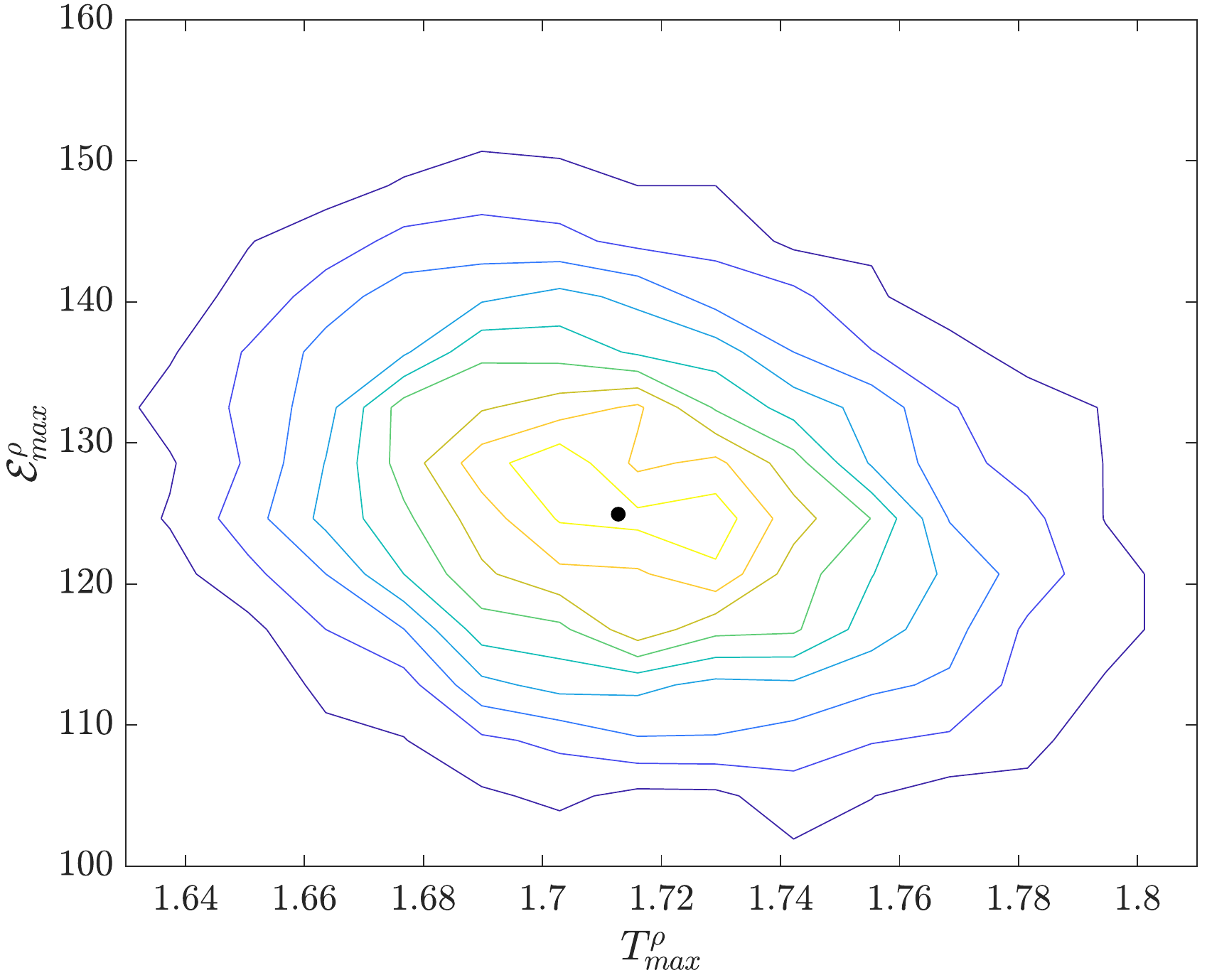}}}
\mbox{
\subfigure[]{\includegraphics[width=0.48\textwidth]{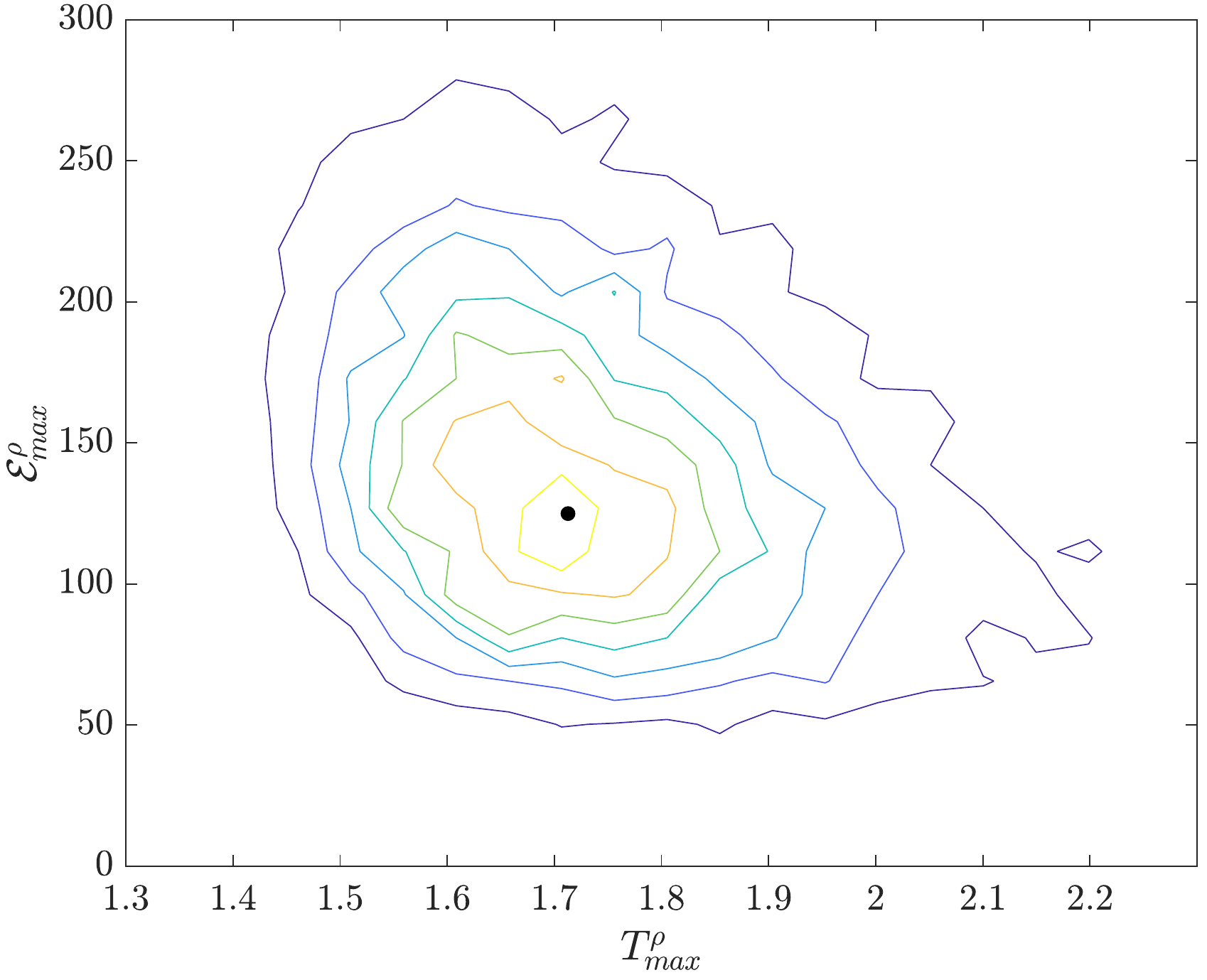}}
\subfigure[]{\includegraphics[width=0.49\textwidth]{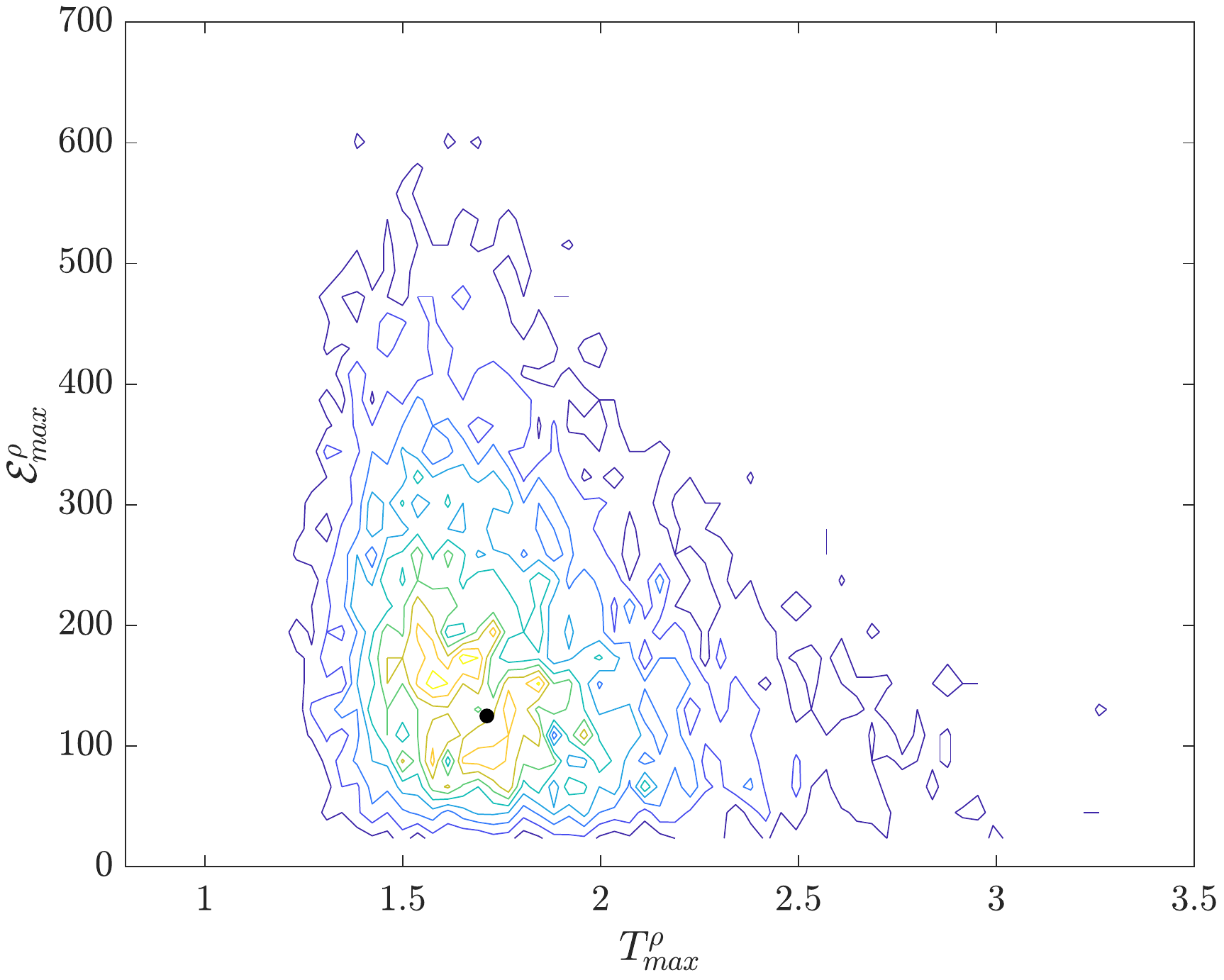}}}
\caption{Equispaced level sets of the JPDFs of $\Tm$ and $\Em$ for the
  noise amplitudes (a) $\rho=10^{-4}$, (b) $\rho=10^{-2}$, (c)
  $\rho=5\times10^{-2}$ and (d) $\rho=10^{-1}$. Yellow level
    sets correspond to large values of the JPDF, whereas black symbols
    represent the values of $T_{\text{max}}$ and $\E_{\text{max}}$ in
  the deterministic case.}
\label{fig:JPDF}
\end{figure}

In order to elucidate the correlation between the values of $\Em$ and
$\Tm$ in different stochastic realizations, joint probability density
functions (JPDFs) of these quantities are shown in Figures
\ref{fig:JPDF}a--d for different values of $\rho$. We see that as
the noise amplitude grows, larger values of $\Em$ become increasingly
correlated with times $\Tm$ shorter than in the deterministic case. By
the same token, reduced values of $\Em$ become increasingly correlated
with times $\Tm$ longer than in the deterministic case. Interestingly,
for larger values of $\rho$ there is a well-defined minimum time
$\Tm$, which depends on $\rho$, such that maximum values of
enstrophy are unlikely to be obtained before that time. We also add
that we compute our stochastic realizations over a fixed time window
$[0,T]$ with $T = 4$. Therefore, as $\rho$ grows and large values of
$\Tm \lessapprox 4$ become more likely, it becomes increasingly
difficult to capture enstrophy maxima as they begin to occur near the
end of the considered time window.

\begin{figure}[h]
\centering
\mbox{
\subfigure[]{\includegraphics[width=0.48\textwidth]{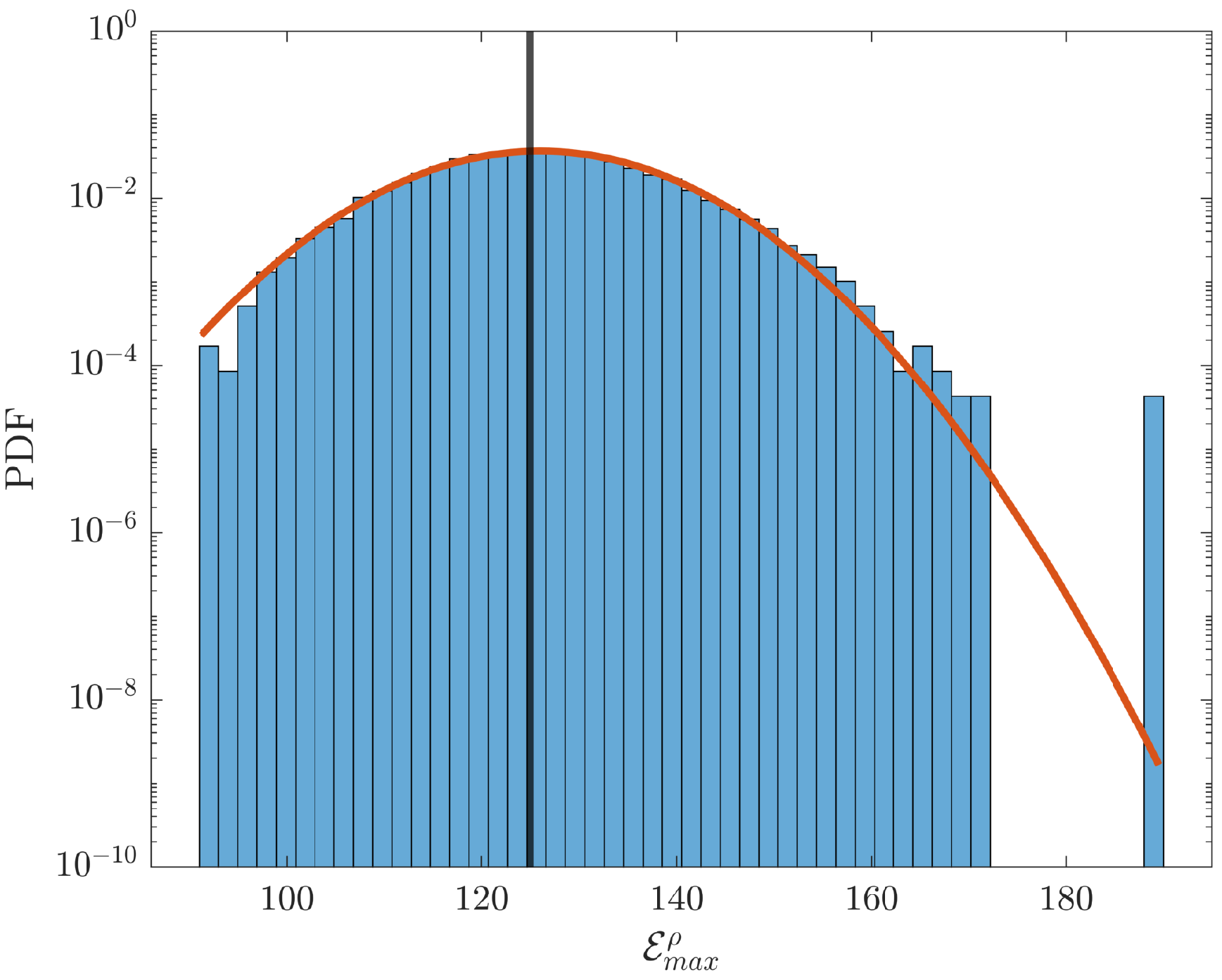}}
\subfigure[]{\includegraphics[width=0.48\textwidth]{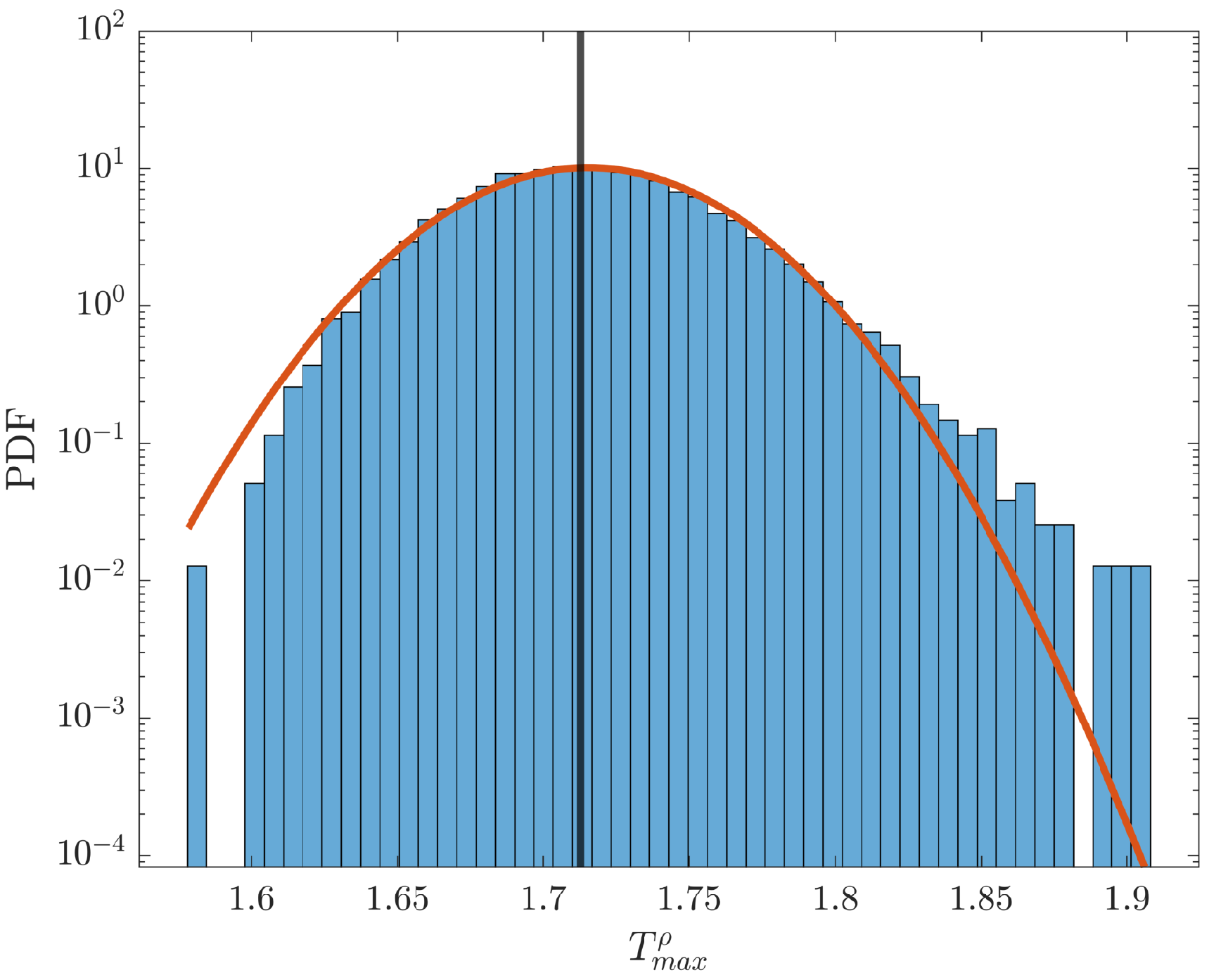}}}
\mbox{
\subfigure[]{\includegraphics[width=0.48\textwidth]{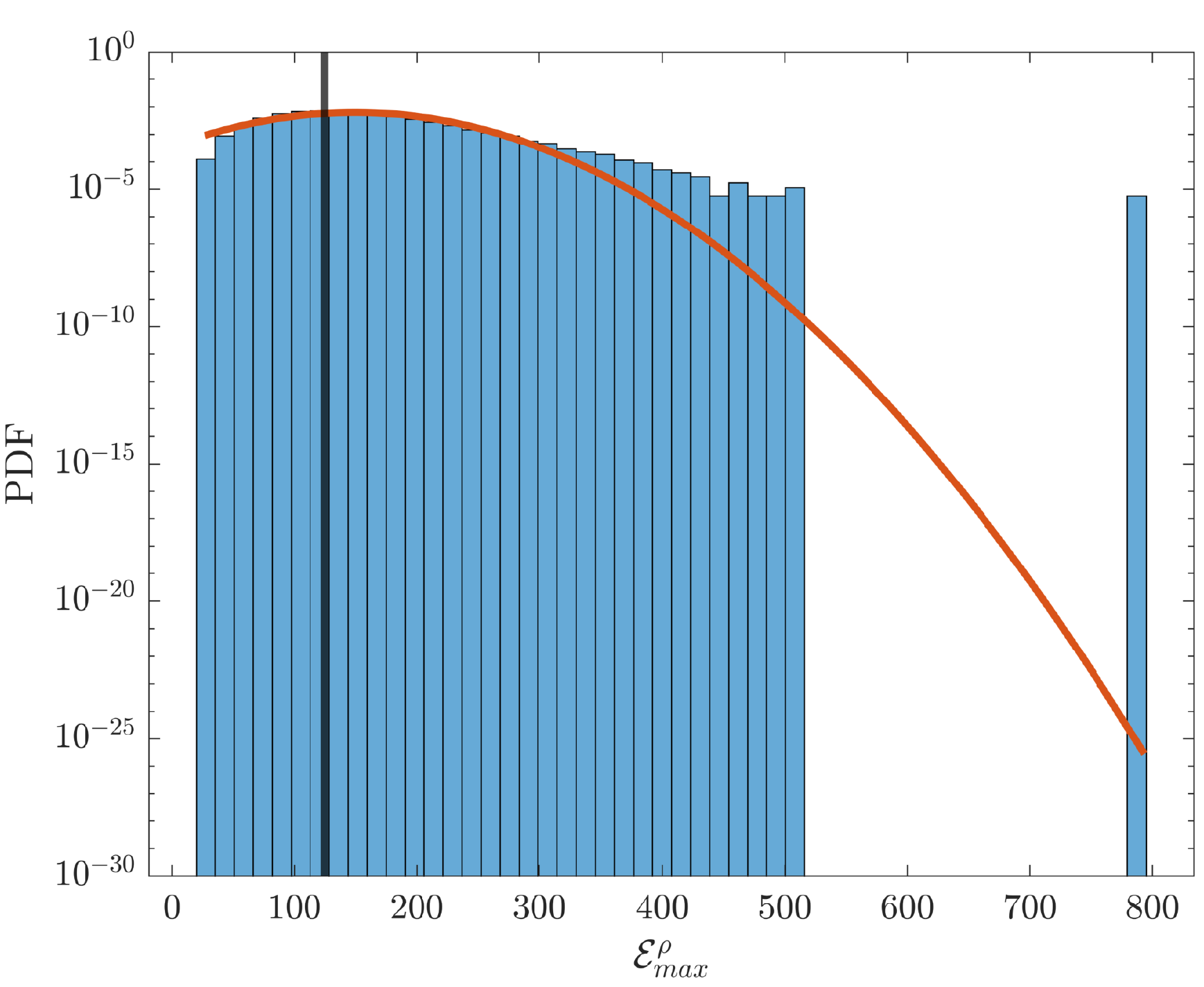}}
\subfigure[]{\includegraphics[width=0.48\textwidth]{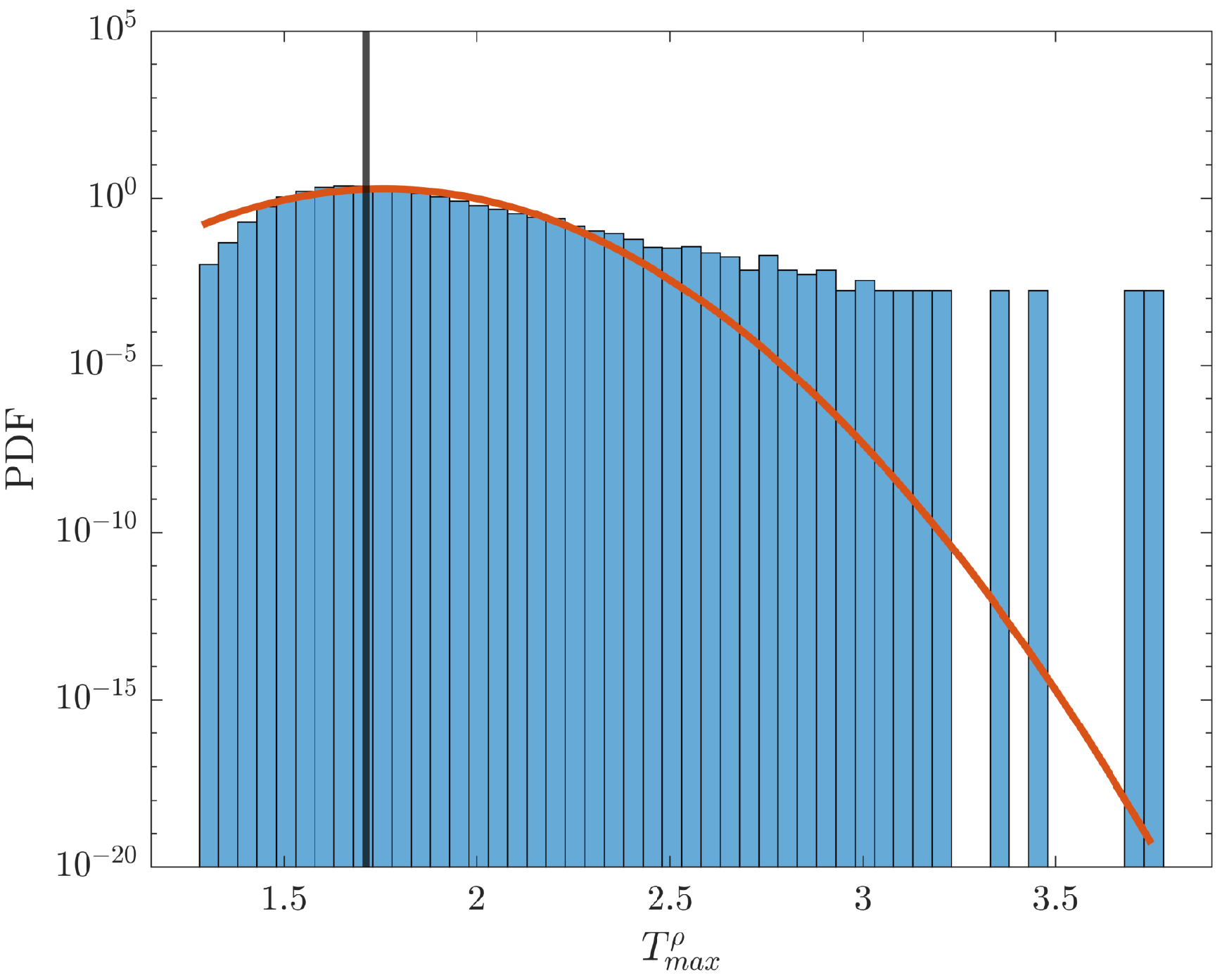}}}
\caption{PDFs of (a,c) the maximum attained enstrophy $\Em$ and (b,d)
  the time $\Tm$ when the maximum occurs in solutions of the
  stochastic problem \eqref{SBE2} with different noise amplitudes
  (a,b) $\rho=10^{-2}$ and (c,d) $\rho=5\cdot 10^{-2}$.  Red
  curves represent the Gaussian distributions with the same means and
  standard deviations, whereas black vertical lines denote the values
  of $T_{\text{max}}$  in the deterministic case.}
\label{fig:subPDF}
\end{figure}

\begin{figure}[h]
\centering
\mbox{
\subfigure[]{\includegraphics[width=0.48\textwidth]{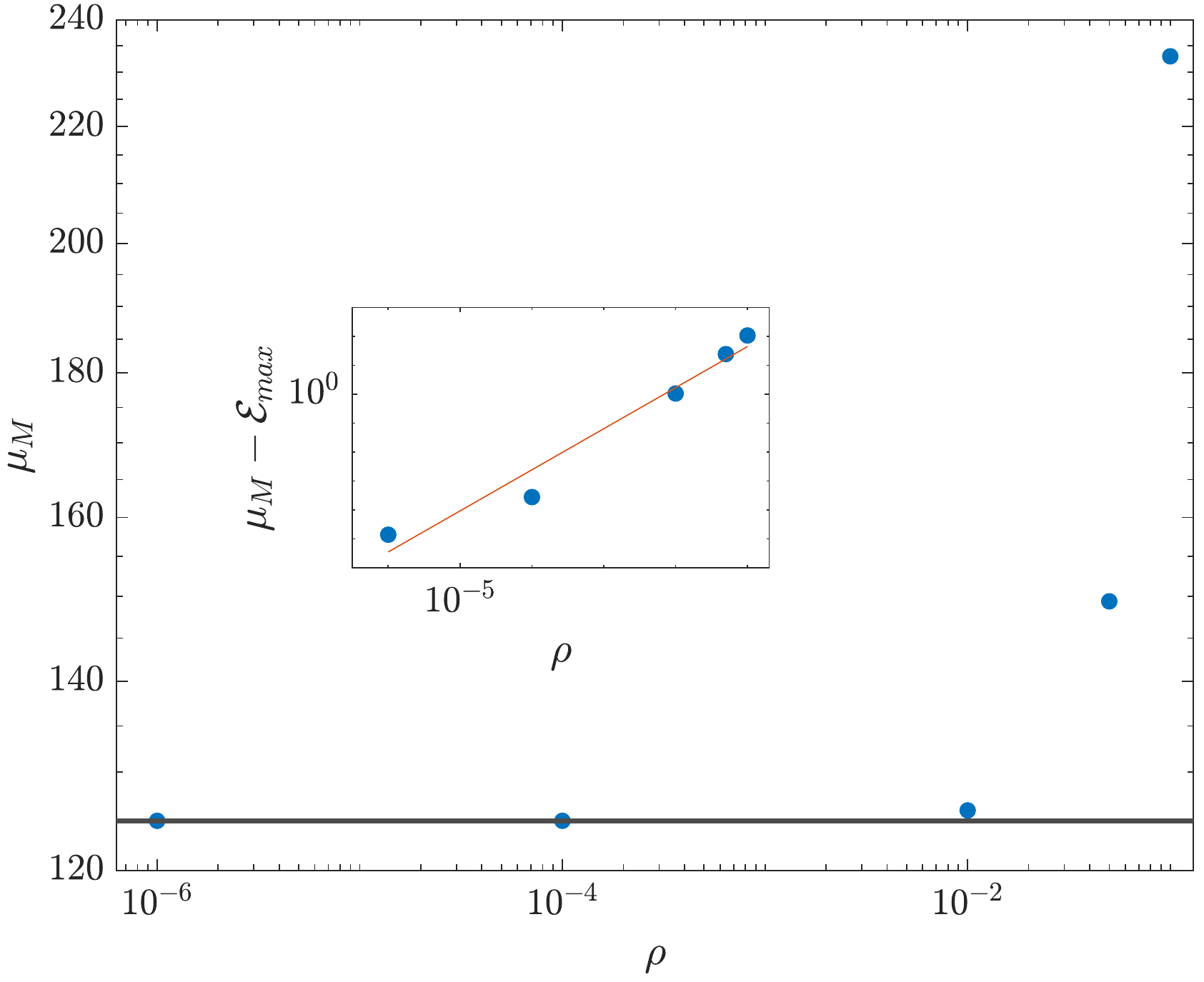}}
\subfigure[]{\includegraphics[width=0.48\textwidth]{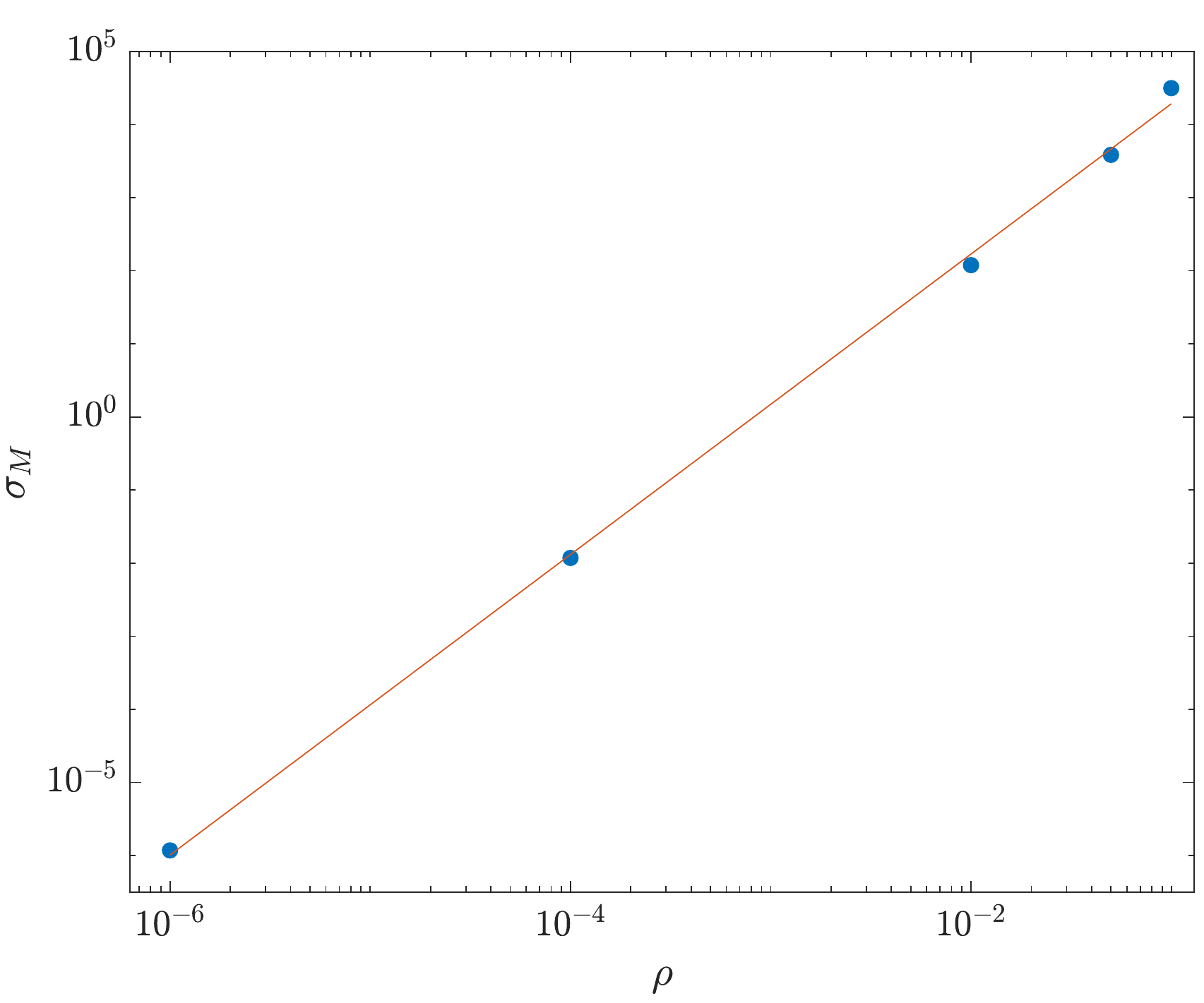}}}
\mbox{
\subfigure[]{\includegraphics[width=0.48\textwidth]{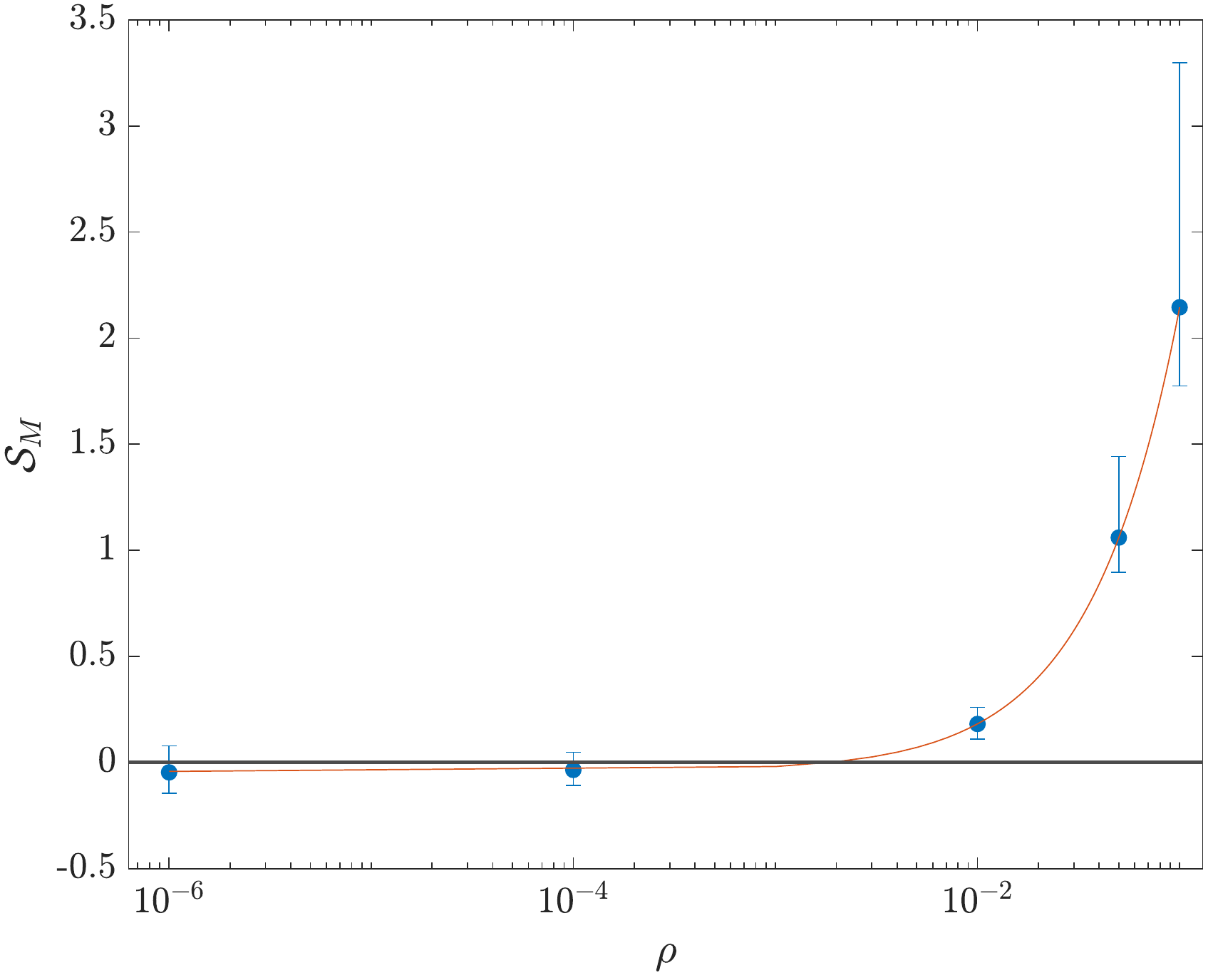}}
\subfigure[]{\includegraphics[width=0.48\textwidth]{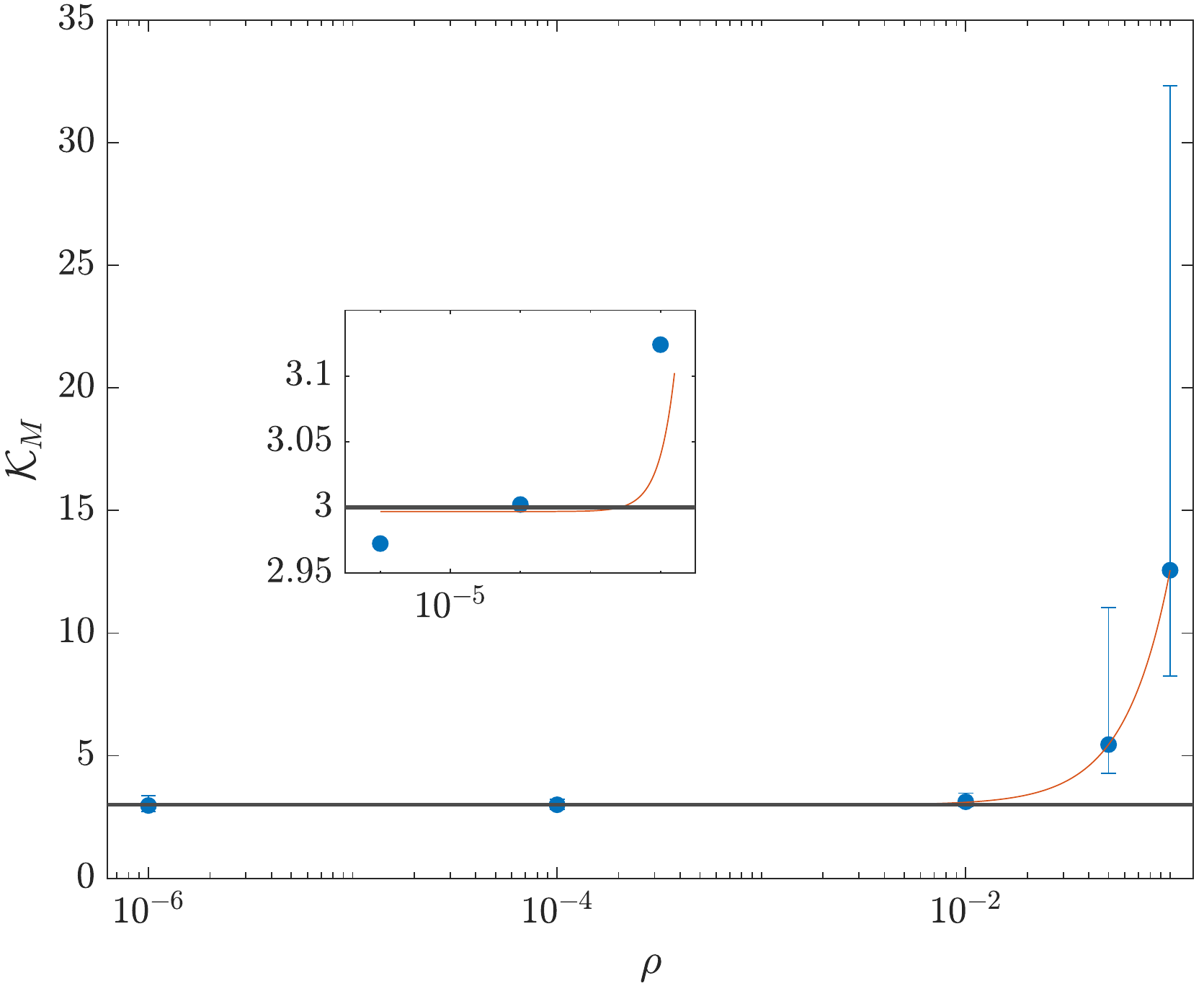}}}
\caption{Estimated (a) mean maximum enstrophy $\Em$, (b) its standard
  deviation, (c) skewness and (d) kurtosis as functions of the noise
  amplitude $\rho$. The horizontal line in (a) corresponds to the
  maximum enstrophy in the deterministic case, whereas the horizontal
  lines in (c) and (d) represent the values characterizing the
  Gaussian distribution. The red curves correspond to the fits
    obtained using the power-law relation \eqref{eq:fit} with the
    parameters reported in Table \ref{tab:FitEmsig}.}
\label{fig:Emsig}
\end{figure}

\begin{figure}[h]
\centering
\mbox{
\subfigure[]{\includegraphics[width=0.48\textwidth]{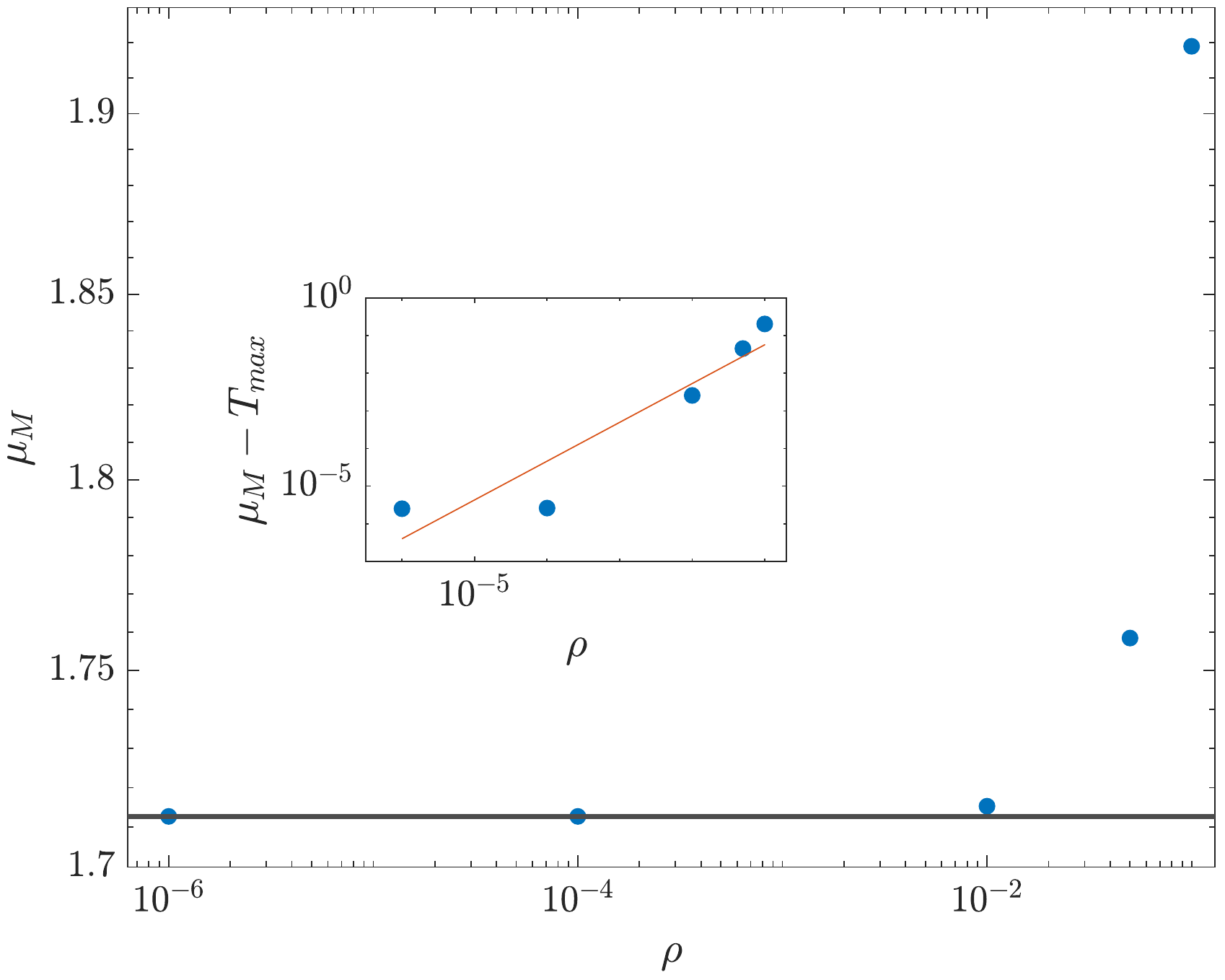}}
\subfigure[]{\includegraphics[width=0.48\textwidth]{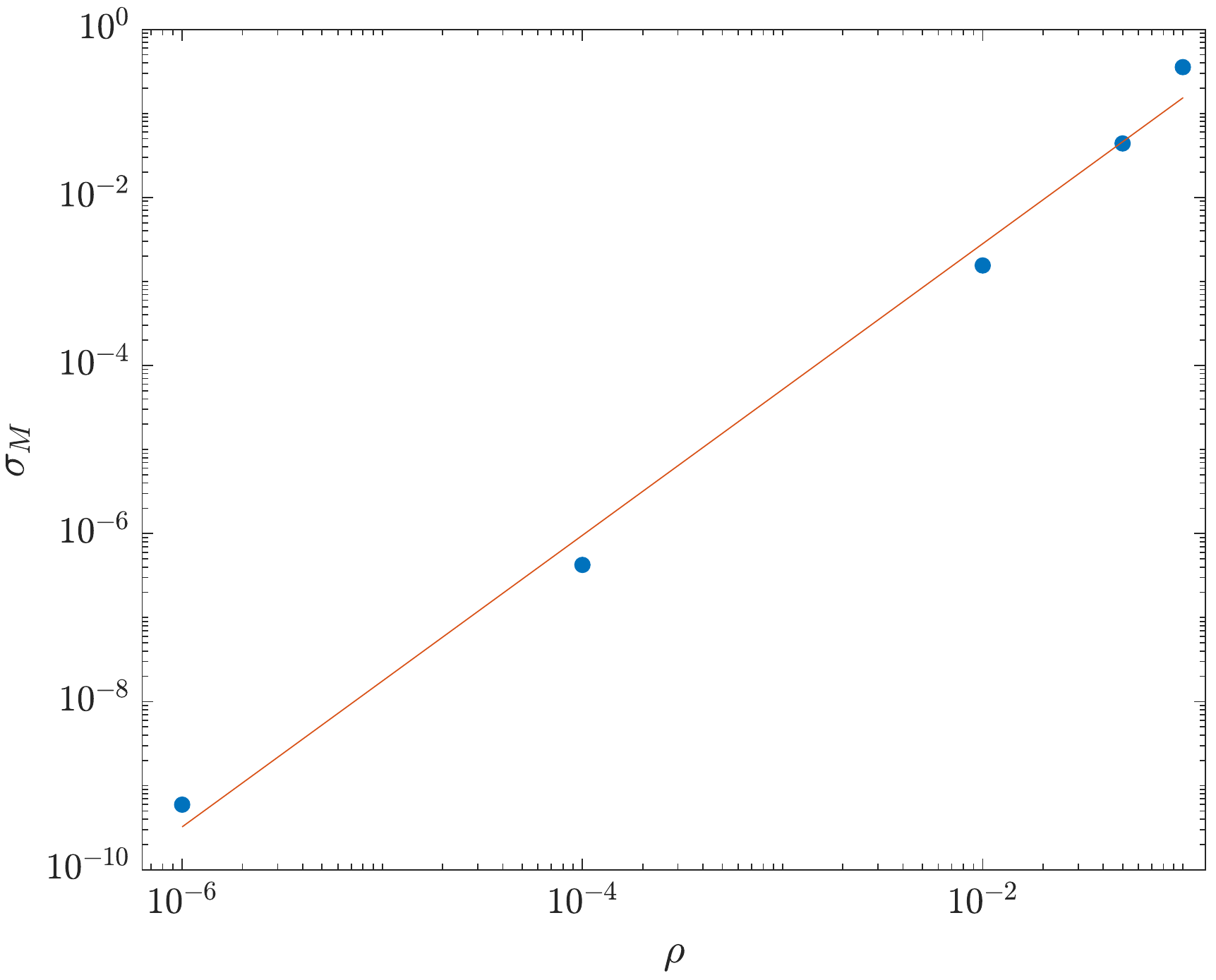}}}
\mbox{
\subfigure[]{\includegraphics[width=0.48\textwidth]{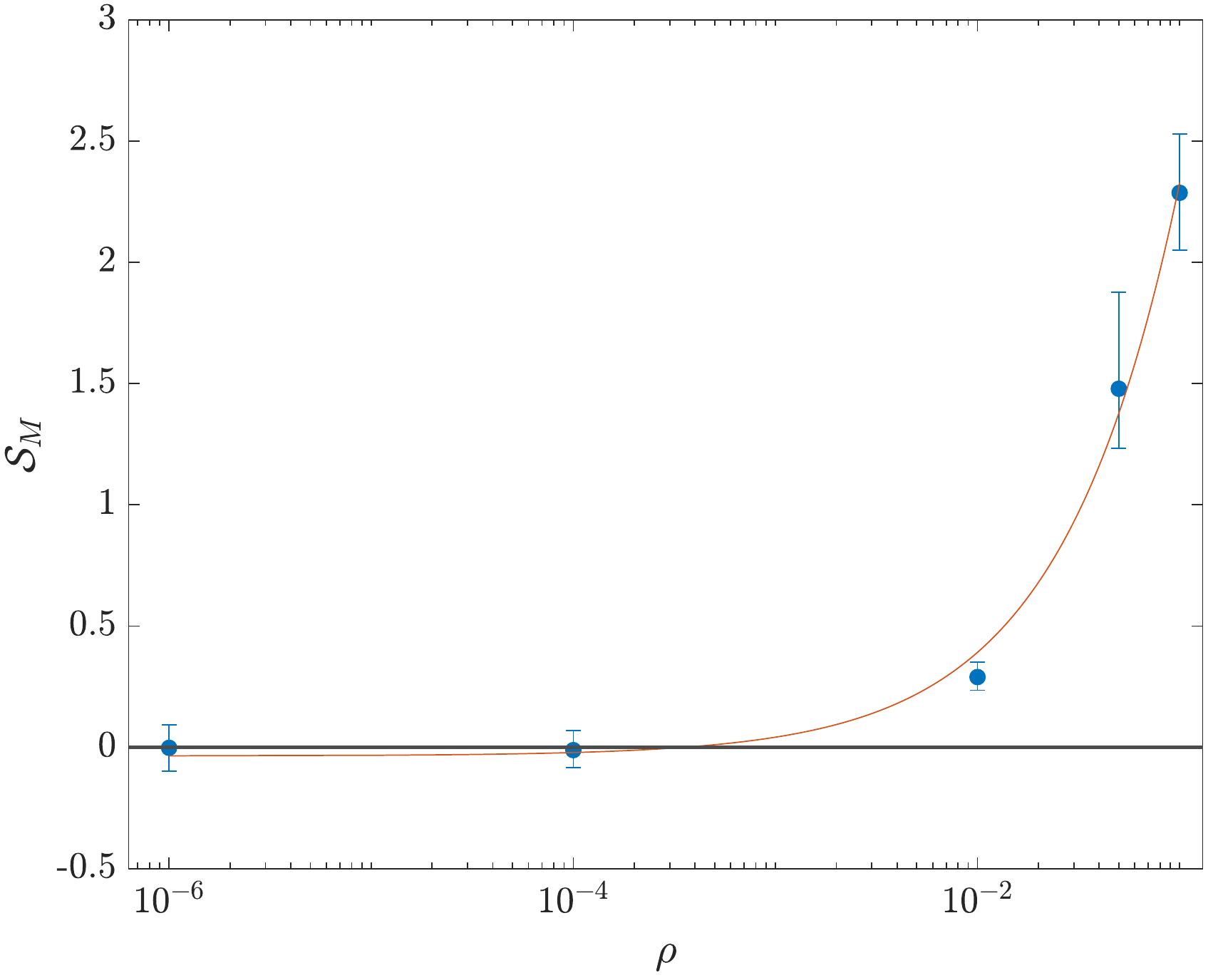}}
\subfigure[]{\includegraphics[width=0.48\textwidth]{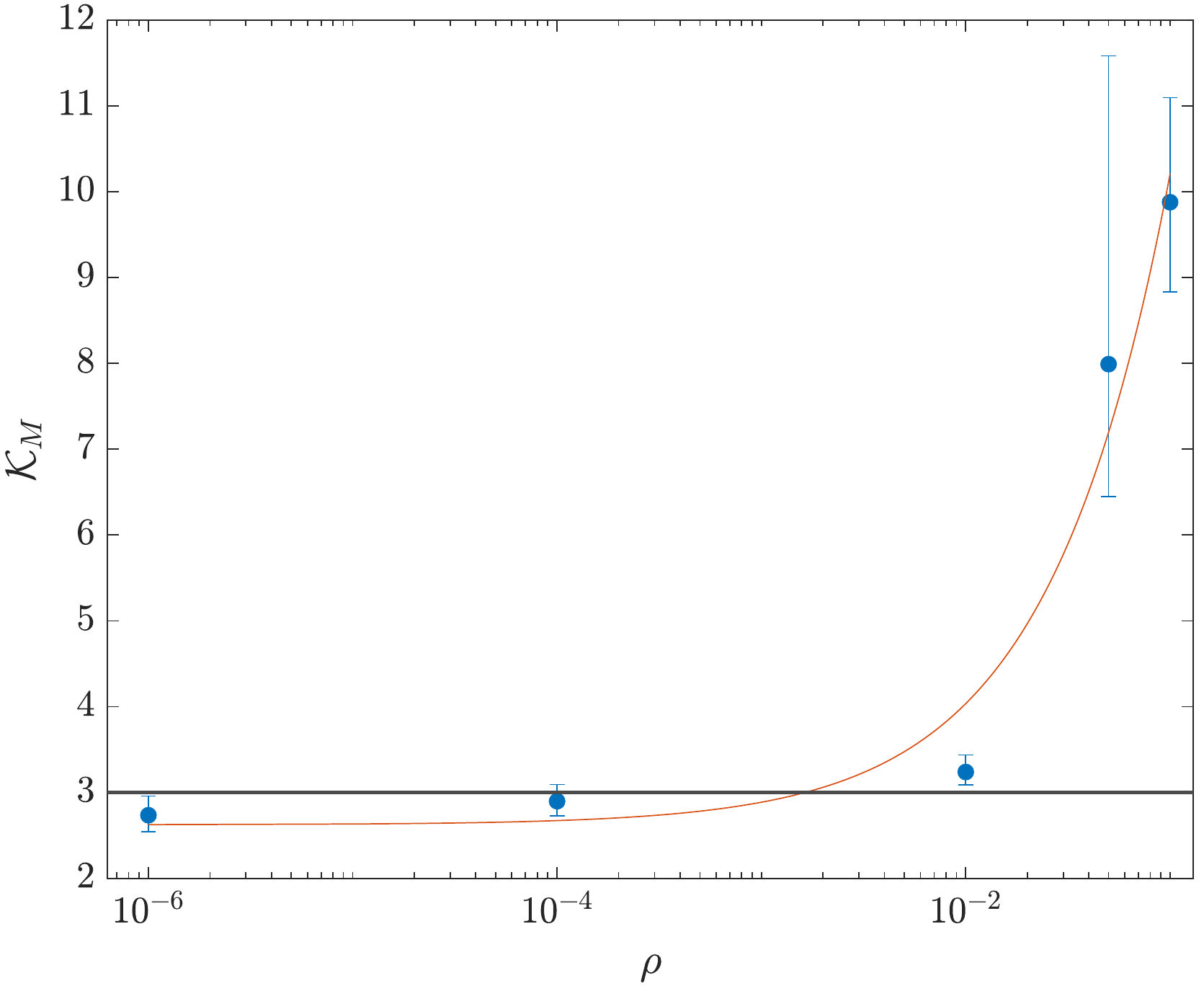}}}
\caption{Estimated (a) mean time $\Tm$ when the maximum enstrophy is
  attained, (b) its standard deviation, (c) skewness and (d) kurtosis
  as functions of the noise amplitude $\rho$. The horizontal line in
  (a) corresponds to the time where the maximum enstrophy is attained
  in the deterministic case, whereas the horizontal lines in (c) and
  (d) represent the values characterizing the Gaussian distribution.
  The red curves correspond to the fits obtained using the
    power-law relation \eqref{eq:fit} with the parameters reported in
    Table \ref{tab:FitTmsig}.}
\label{fig:Tmsig}
\end{figure}

In order to shed further light on the trends observed in Figures
\ref{fig:JPDF}a--d, the PDFs of $\Em$ and $\Tm$ for two selected
values of $\rho$ are shown in Figures \ref{fig:subPDF}a--b. It is
evident that as the noise amplitude becomes larger, the PDFs of both
$\Em$ and $\Tm$ become increasingly non-Gaussian, with this trend
appearing more pronounced in the case of $\Tm$. In particular, the
PDFs become strongly skewed with heavy tails towards large values of
$\Em$ and $\Tm$. On the other hand, values of these quantities
significantly smaller than in the deterministic case are much less
likely to occur. To analyze these trends in quantitative terms, the
four first statistical moments of the distributions of $\Em$ and
$\Tm$, defined analogously to \eqref{eq:mu}--\eqref{eq:K}, are plotted
in Figures \ref{fig:Emsig} and \ref{fig:Tmsig} as functions of the
noise amplitude $\rho$. We remark that the dependence of the
estimates of these moments on the number of samples used in their
evaluation is similar to what was shown in Figures
\ref{fig:m1}--\ref{fig:m2}, and these results are omitted here for
brevity. 

In Figures \ref{fig:Emsig}a we see that the mean value of $\Em$
increases with $\rho$ which is consistent with the trend observed in
Figure \ref{fig:subEt}c for the expected value of the enstrophy
$\EE[\E(u(t))]$, although the mean values of $\Em$ exceed the values
of $\max_{t>0} \EE[\E(u(t))]$ for all $\rho$. On the other hand, in
Figures \ref{fig:Tmsig}a we see that the mean time $\Tm$ when the
maximum enstrophy is attained also increases with $\rho$ which should
be contrasted with Figure \ref{fig:subEt}c showing that the maxima of
the expected value of the enstrophy $\EE[\E(u(t))]$ occur at earlier
times as the noise amplitude increases. The reason for this somewhat
surprising effect is that the operation of computing the maximum of a
function, or the argument corresponding to the maximum, is not linear
and hence does not commute with the operation of computing the
average.

In order to quantify the dependence of the statistical moments
of $\Em$ and $\Tm$ on the noise magnitude $\rho$, Figures
\ref{fig:Emsig}a--d and \ref{fig:Tmsig}a--d also include fits
obtained using the power-law relation \eqref{eq:fit} with parameters
reported in Tables \ref{tab:FitEmsig} and \ref{tab:FitTmsig}. As was
the case for the moments of the blow-up time in the supercritical
regime, cf.~Figure \ref{fig:msig}, we see that while the dependence
of the variance, skewness and kurtosis of the distributions of $\Em$
and $\Tm$ on $\rho$ is represented rather well by relation
\eqref{eq:fit}, this is not the case for the mean values of $\Em$
and $\Tm$. Interestingly, the values of the exponent $b$ describing
the dependence of the variance, skewness and kurtosis of the
distributions of $\Em$ on $\rho$, cf.~Table \ref{tab:FitEmsig}, are
similar to those reported in Table \ref{tab:mf} for the blow-up time
in the supercritical regime. Other than this, there seems to be no
obvious pattern in how the different statistical moments depend on
the noise amplitude $\rho$.

\begin{table}[t]
\centering
\begin{tabular}{|c|c|c|c|}
\hline
Statistical Moment&a&b&c \\ \hline
$\mu_M$ & 175.027 & 1.42&0\\ \hline
$\sigma_M$ & $2.18\cdot10^{6}$&2.05&0\\\hline
$\mathcal{S}_M$ & 21.33&1.13&-0.041\\ \hline
$\mathcal{K}_M$ & 868.48&1.95&2.99\\\hline
\end{tabular}
\caption{Parameters of the fits to the data shown in Figures 
\ref{fig:Emsig}a--d using the power-law relation  \eqref{eq:fit}.}
\label{tab:FitEmsig}
\bigskip
\begin{tabular}{|c|c|c|c|}
\hline
Statistical Moment&a&b&c \\ \hline
$\mu_M$ & 0.62 & 1.03&0\\ \hline
$\sigma_M$ & 8.307&1.73&0\\\hline
$\mathcal{S}_M$ & 13.07&0.74&-0.035\\ \hline
$\mathcal{K}_M$ & 40.72&0.72&2.62\\\hline
\end{tabular}
\caption{Parameters of the fits to the data shown in Figures 
\ref{fig:Tmsig}a--d using the power-law relation \eqref{eq:fit}.}
\label{tab:FitTmsig}
\end{table}

\FloatBarrier

\section{Summary and Conclusions}
\label{sec:final}

This study is motivated by the question of how singularity formation
and other forms of extreme behavior in nonlinear dissipative partial
differential equations are affected by stochastic excitations
\cite{f15}. As a model problem, we have considered the 1D fractional
Burgers equation \eqref{FBE} with additive colored noise. This system
is interesting, because in the deterministic setting it exhibits
finite-time blow-up or a globally well-posed behavior depending on the
value of the fractional dissipation exponent $\alpha$. The question we
are interested in is addressed by performing a series of very accurate
numerical computations combining spectrally-accurate spatial
discretizations with a Monte-Carlo approach where convergence of all
approximations was checked very carefully. It should be emphasized
that even though the problem is formulated in 1D, these computations
are in fact quite challenging since resolving an emerging singularity
requires refined numerical resolutions both in space and in time, and
this difficulty is compounded by a slow convergence of the Monte-Carlo
approach.

As the first main contribution, we carefully documented the
singularity formation in the deterministic system in the supercritical
regime. It was shown that for a fixed initial condition the blow-up
time $T^*$ is a decreasing function of the fractional dissipation
exponent $\alpha \in [0,1/2)$ and tends to the expected limiting
values as $\alpha \rightarrow 0$ and $\nu \rightarrow 0$, cf.~Figures
\ref{fig:Ts}b and \ref{fig:Tnu0}. When $\alpha \rightarrow (1/2)^-$
the blow-up time approaches a finite value.

Our second main finding is that while there is no evidence for noise
to prevent singularity formation in the supercritical regime, the
distribution of the blow-up times $T^*$ (understood as a stochastic
variable) becomes increasingly non-Gaussian as the noise amplitude
$\rho$ increases, cf.~Figure \ref{fig:supPDF}. Interestingly, as
$\rho$ grows, the mean blow-up time is reduced but at the same time
significantly delayed blow-up times also become more likely.  This is
because the PDFs of $T^*$ become more skewed and develop heavy tails
towards large blow-up times as $\rho$ increases. These observations
are complemented by our findings for the subcritical case where we
noted that the expected value of the enstrophy $\EE[\E(t)]$ of the
stochastic solution at any fixed time $t$ increases with the noise
magnitude $\rho$ (Figure \ref{fig:subEt}c). Interestingly, while its
maxima are attained at earlier times as compared to the deterministic
case, the mean time when the enstrophy maxima are attained in
individual stochastic realizations is in fact larger than in the
deterministic case and increases with $\rho$, cf.~Figure
\ref{fig:Tmsig}a.  Increasingly non-Gaussian behavior of the PDFs of
$\Em$ and $\Tm$ is evident as $\rho$ becomes large in the subcritical
case we well (Figures \ref{fig:subPDF}a--d).

To conclude, as is evident from the discussion above, the answer to
the question about the effect of stochastic excitations on extreme and
singular behavior in fractional Burgers flows is rather nuanced.  It
is clear, however, that there is no evidence for the noise to
regularize the evolution by suppressing blow-up in the supercritical
regime, or for the noise to trigger blow-up in the subcritical regime.
It has to be recognized that due to the computational cost the results
reported here are restricted to small and intermediate noise
magnitudes only.  Whether or not the trends reported in this study
will hold also for much larger noise amplitudes is an interesting open
question.

\section*{Acknowledgments}

The authors acknowledge funding from McMaster University and through
an NSERC (Canada) Discovery Grant.  Computational resources were
provided by Compute Canada under its Resource Allocation Competition.




\end{document}